\documentclass[a4paper]{amsart}
\usepackage{graphicx}
\usepackage{amsmath}
\usepackage{amssymb}
\usepackage{amsfonts}
\usepackage{amsthm}
\usepackage{eucal}
\usepackage{tikz}
\usetikzlibrary{trees, positioning, calc, cd}
\usepackage{stmaryrd}

\theoremstyle{plain}
\newtheorem{theorem}{Theorem}[section]

\newtheorem{lemma}[theorem]{Lemma}
\newtheorem{proposition}[theorem]{Proposition}

\newtheorem{corollary}[theorem]{Corollary}

\newtheorem{conjecture}[theorem]{Conjecture}

\newtheorem*{theorem*}{Theorem}
\newtheorem*{proposition*}{Proposition}
\newtheorem*{corollary*}{Corollary}

\theoremstyle{definition}
\newtheorem{definition}[theorem]{Definition}

\newtheorem{remark}[theorem]{Remark}

\makeatletter
\@addtoreset{theorem}{section}
\@addtoreset{subsection}{section}
\makeatother

\begin{document}

\title{Koszul duality and a conjecture of Francis--Gaitsgory}

\author{Gijs Heuts}

\date{}

\begin{abstract}
Koszul duality is a fundamental correspondence between algebras for an operad $\mathcal{O}$ and coalgebras for its dual cooperad $B\mathcal{O}$, built from $\mathcal{O}$ using the bar construction. Francis--Gaitsgory proposed a conjecture about the general behavior of this duality. The main result of this paper, roughly speaking, is that Koszul duality provides an equivalence between the subcategories of nilcomplete algebras and conilcomplete coalgebras and that these are the largest possible subcategories for which such a result holds. This disproves Francis--Gaitsgory's prediction, but does provide an adequate replacement. We show that many previously known partial results about Koszul duality can be deduced from our results.
\end{abstract}

\maketitle

\setcounter{tocdepth}{1}
\tableofcontents

\section{Introduction}

In his work on rational homotopy theory Quillen \cite{rationalhomotopy} exhibited a close relationship between Lie algebras and (coaugmented) cocommutative coalgebras over the rational numbers. This relationship can be interpreted as an adjoint pair of functors between $\infty$-categories
\[
\begin{tikzcd}
\mathrm{Lie}(\mathrm{Ch}_{\mathbb{Q}}) \ar[shift left]{r}{\mathrm{CE}} & \mathrm{coCAlg}^{\mathrm{aug}}(\mathrm{Ch}_{\mathbb{Q}}) \ar[shift left]{l}{\mathrm{prim}},
\end{tikzcd}
\]
where the left-hand (resp. right-hand) side is the $\infty$-category of differential graded Lie algebras (resp. coaugmented cocommutative coalgebras) over $\mathbb{Q}$. The left adjoint is a version of the complex computing Chevalley--Eilenberg homology, whereas the right adjoint takes the (derived) primitives of a coalgebra. Moreover, upon restricting to connected Lie algebras on the left and simply-connected coalgebras on the right, this pair becomes an adjoint equivalence. Another duality of a similar spirit arises in the work of Moore \cite{moore}, who established bar and cobar constructions relating associative algebras and coalgebras. In modern language (e.g. as in \cite[Section 5.2.2]{higheralgebra}), these constructions form an adjoint pair between $\infty$-categories
\[
\begin{tikzcd}
\mathrm{Alg}^{\mathrm{aug}}(\mathcal{C}) \ar[shift left]{r}{\mathrm{Bar}} & \mathrm{coAlg}^{\mathrm{aug}}(\mathcal{C})  \ar[shift left]{l}{\mathrm{Cobar}},
\end{tikzcd}
\]
where $\mathcal{C}$ is a monoidal $\infty$-category satisfying only minimal hypotheses. In classical examples $\mathcal{C}$ would be the category of chain complexes over a field, but the setup works in great generality.

Both of these dualities are instances of the general phenomenon of \emph{Koszul duality} for operads, which was first identified in the works of Ginzburg--Kapranov \cite{ginzburgkapranov} and Getzler--Jones \cite{getzlerjones}. In its general form, Koszul duality (or bar-cobar duality) gives a relation between algebras for an operad $\mathcal{O}$ and coalgebras for its bar construction (or `Koszul dual') $B\mathcal{O}$, which is a cooperad. The examples of the previous paragraph are special instances where the operad $\mathcal{O}$ is taken to be the Lie operad or the associative operad. Koszul duality plays a central role in modern treatments of formal deformation theory \cite{brantnermathew,calaquecamposnuiten,dagx,pridham} and in the theory of factorization homology (or chiral homology) \cite{ayalafrancis,francisgaitsgory}. 

Francis and Gaitsgory \cite{francisgaitsgory} give a description of operadic Koszul duality as an adjoint pair between the $\infty$-category of $\mathcal{O}$-algebras and that of $B\mathcal{O}$-coalgebras `with divided powers' (a notion we describe in detail later). They conjecture that this duality adjunction should restrict to an equivalence between the subcategories of \emph{pronilpotent} $\mathcal{O}$-algebras and \emph{ind-conilpotent} $B\mathcal{O}$-coalgebras with divided powers (cf. Conjecture 3.4.5 of \cite{francisgaitsgory}). Several special cases of this conjecture have been proved, e.g. in the works of Ching--Harper \cite{chingharper} and Brantner--Mathew \cite{brantnermathew}. In this paper we prove the following:

\begin{theorem*}[Theorem \ref{thm:completeness}]
Koszul duality implements an equivalence between the $\infty$-category of \emph{nilcomplete} $\mathcal{O}$-algebras and that of \emph{conilcomplete} $B\mathcal{O}$-coalgebras with divided powers. These subcategories are the largest possible for which such an equivalence is achieved: on any other object, the (co)unit of the Koszul duality adjunction is not an isomorphism.
\end{theorem*}

This result is different from Francis--Gaitsgory's prediction; we will use it to produce counterexamples to their conjecture. However, we do demonstrate that many special cases of the conjecture of Francis--Gaitsgory treated in the literature can be seen as instances of our result.

\subsection*{Acknowledgments} 
I thank Max Blans, Lukas Brantner, Hadrian Heine, Jacob Lurie, and Yuqing Shi for interesting conversations about Koszul duality. I thank the Institute for Advanced Study, where part of this work was done, for its hospitality. This work was supported by an ERC Starting Grant (no.\ 950048) and an NWO Vidi grant (no.\ 223.093).

\section{Main results}
\label{sec:mainresults}

Throughout this paper we let $\mathcal{C}$ be a presentable, stable, symmetric monoidal $\infty$-category for which the tensor product commutes with colimits in each variable separately (see Remark \ref{rmk:hypotheses}). We consider operads $\mathcal{O}$ in $\mathcal{C}$ for which $\mathcal{O}(0) = 0$ (i.e. $\mathcal{O}$ is nonunital) and $\mathcal{O}(1)$ is the tensor unit of $\mathcal{C}$. We will review a version of the theory of operads and cooperads we need here in Section \ref{sec:truncations}.

An $\mathcal{O}$-algebra is an object $X$ of $\mathcal{C}$ equipped with a system of `multiplication maps'
\begin{equation*}
(\mathcal{O}(n) \otimes X^{\otimes n})_{h\Sigma_n} \rightarrow X
\end{equation*}
for each $n \geq 1$, and a further coherent system of homotopies expressing the associativity and unitality of these multiplications. We will discuss the $\infty$-category $\mathrm{Alg}_{\mathcal{O}}(\mathcal{C})$ of $\mathcal{O}$-algebras in $\mathcal{C}$ in more detail in Section \ref{sec:coalgebras}. It is related to $\mathcal{C}$ by two different adjoint pairs:
\[
\begin{tikzcd}[column sep = large]
\mathcal{C} \ar[shift left]{r}{\mathrm{free}_{\mathcal{O}}} & \mathrm{Alg}_{\mathcal{O}}(\mathcal{C}) \ar[shift left]{l}{\mathrm{forget}_{\mathcal{O}}} \ar[shift left]{r}{\mathrm{cot}_{\mathcal{O}}} & \mathcal{C}. \ar[shift left]{l}{\mathrm{triv}_{\mathcal{O}}}
\end{tikzcd}
\]
Here the left adjoints are indicated on top. The pair on the left is the usual free-forgetful adjunction. The functor $\mathrm{triv}_{\mathcal{O}}$ assigns to an object of $\mathcal{C}$ the trivial $\mathcal{O}$-algebra structure on that object, for which all multiplication maps (except for the unit) are null. It admits a left adjoint $\mathrm{cot}_{\mathcal{O}}$ which is often referred to as the \emph{indecomposables}, or \emph{cotangent fiber}, or \emph{topological Andr\'{e}--Quillen homology}. It can be thought of as the universal homology theory for $\mathcal{O}$-algebras (cf. Remark \ref{rmk:homology}) and will be discussed in Section \ref{sec:conilpotentcoalgebras}.

The cotangent fiber $\mathrm{cot}_{\mathcal{O}}(A)$ is more than just an object of $\mathcal{C}$. Indeed, it can be upgraded to a coalgebra for the cooperad $B\mathcal{O}$, which is obtained as the \emph{bar construction} of $\mathcal{O}$ (to be discussed in Section \ref{sec:truncations}). Generally, for a cooperad $\mathcal{Q}$ in $\mathcal{C}$, a $\mathcal{Q}$-coalgebra is an object $X$ of $\mathcal{C}$ equipped with a collection of `comultiplication maps'
\begin{equation*}
X \rightarrow (\mathcal{Q}(n) \otimes X^{\otimes n})^{h\Sigma_n}, \quad n \geq 1,
\end{equation*}
together with a coherent system of homotopies expressing the associativity and unitality of this structure. Furthermore, a \emph{coalgebra with divided powers} is an object $X$ equipped with comultiplication maps
\begin{equation*}
X \rightarrow (\mathcal{Q}(n) \otimes X^{\otimes n})_{h\Sigma_n}, \quad n \geq 1,
\end{equation*}
and homotopies expressing their associativity and unitality. Note that we have replaced the homotopy invariants with homotopy coinvariants in this formula. These divided powers should be regarded as extra structure on a coalgebra. Indeed, using the norm maps
\begin{equation*}
\mathrm{Nm}_{\Sigma_n}\colon (\mathcal{Q}(n) \otimes X^{\otimes n})_{h\Sigma_n} \rightarrow (\mathcal{Q}(n) \otimes X^{\otimes n})^{h\Sigma_n}
\end{equation*}
associated with the symmetric groups $\Sigma_n$, any coalgebra with divided powers in particular has an underlying coalgebra. We will discuss the relevant $\infty$-categories $\mathrm{coAlg}_{\mathcal{Q}}(\mathcal{C})$ and $\mathrm{coAlg}_{\mathcal{Q}}^{\mathrm{dp}}(\mathcal{C})$ of (divided power) $\mathcal{Q}$-coalgebras in Section \ref{sec:coalgebras} and Appendix \ref{app:coalgebras}. There we will also discuss how the cotangent fiber naturally carries the structure of a divided power $B\mathcal{O}$-coalgebra, giving a factorization of $\mathrm{cot}_{\mathcal{O}}$ as
\begin{equation*}
\mathrm{Alg}_{\mathcal{O}} \xrightarrow{\mathrm{indec}_{\mathcal{O}}} \mathrm{coAlg}^{\mathrm{dp}}_{B\mathcal{O}}(\mathcal{C}) \xrightarrow{\mathrm{forget}_{B\mathcal{O}}} \mathcal{C},
\end{equation*}
where we will refer to the first functor as the \emph{indecomposables}. This functor is part of an adjoint pair 
\[
\begin{tikzcd}
\mathrm{Alg}_{\mathcal{O}} \ar[shift left]{r}{\mathrm{indec}_{\mathcal{O}}} & \mathrm{coAlg}^{\mathrm{dp}}_{B\mathcal{O}}(\mathcal{C}). \ar[shift left]{l}{\mathrm{prim}_{B\mathcal{O}}}
\end{tikzcd}
\]
The right adjoint $\mathrm{prim}_{B\mathcal{O}}$ takes the \emph{primitives} of a divided power $B\mathcal{O}$-coalgebra. The basic question this paper aims to address is the following: on what subcategories does this adjoint pair restrict to an equivalence of $\infty$-categories?

In Section \ref{sec:truncations} we will discuss the notions of \emph{nilpotent completion} of $\mathcal{O}$-algebras and \emph{conilpotent cocompletion} of $B\mathcal{O}$-coalgebras. Briefly, every $\mathcal{O}$-algebra $X$ gives rise to a tower of algebras
\begin{equation*}
X \rightarrow \cdots \rightarrow t_n X \rightarrow t_{n-1} X \rightarrow \cdots \rightarrow t_1 X,
\end{equation*}
where $t_n X$ is obtained by `killing all $k$-fold multiplications' for $k > n$. We call $\varprojlim_n t_n X$ the \emph{nilpotent completion} of $X$ and say $X$ is \emph{nilcomplete} if the map
\begin{equation*}
X \rightarrow \varprojlim_n t_n X
\end{equation*}
is an isomorphism. There are dual notions for (divided power) $B\mathcal{O}$-coalgebras. One of our main results states that Koszul duality restricts to an equivalence between the $\infty$-categories of nilcomplete algebras and conilcomplete divided power coalgebras:

\begin{theorem}
\label{thm:completeness}
The unit of the adjoint pair $(\mathrm{indec}_\mathcal{O}, \mathrm{prim}_{B\mathcal{O}})$ is naturally equivalent to the nilpotent completion of $\mathcal{O}$-algebras, whereas the counit is naturally equivalent to the conilpotent cocompletion of divided power $B\mathcal{O}$-coalgebras. As a consequence, these functors restrict to an adjoint equivalence of $\infty$-categories
\[
\begin{tikzcd}[column sep = large]
\mathrm{Alg}_{\mathcal{O}}(\mathcal{C})^{\mathrm{cpl}} \ar[shift left]{r}{\mathrm{indec}_{\mathcal{O}}} & \mathrm{coAlg}_{B\mathcal{O}}^{\mathrm{dp}}(\mathcal{C})^{\mathrm{cpl}} \ar[shift left]{l}{\mathrm{prim}_{B\mathcal{O}}},
\end{tikzcd}
\]
with the superscripts $\mathrm{cpl}$ denoting the full subcategory of nilcomplete algebras, respectively conilcomplete coalgebras.
\end{theorem}

We will now discuss the setting of the conjecture of Francis--Gaitsgory, which uses a slightly different $\infty$-category of coalgebras. The adjoint pair $(\mathrm{cot}_{\mathcal{O}}, \mathrm{triv}_{\mathcal{O}})$ induces a comonad $\mathrm{cot}_{\mathcal{O}} \circ \mathrm{triv}_{\mathcal{O}}$ on $\mathcal{C}$, which can be described by the formula (cf. Proposition \ref{prop:stabcomonad})
\begin{equation*}
\mathrm{cot}_{\mathcal{O}}\mathrm{triv}_{\mathcal{O}}(X) \cong \bigoplus_{n \geq 1} (B\mathcal{O}(n) \otimes X^{\otimes n})_{h\Sigma_n}.
\end{equation*}
Francis--Gaitsgory \cite{francisgaitsgory} define the $\infty$-category of \emph{ind-conilpotent divided power $B\mathcal{O}$-coalgebras} to be the $\infty$-category of coalgebras for this comonad:
\begin{equation*}
\mathrm{coAlg}^{\mathrm{dp},\mathrm{nil}}_{B\mathcal{O}}(\mathcal{C}) := \mathrm{coAlg}_{\mathrm{cot}_{\mathcal{O}}\mathrm{triv}_{\mathcal{O}}}(\mathcal{C}).
\end{equation*}
The adjective \emph{ind-conilpotent} refers to the direct sum, whereas the term divided powers again refers to the fact that the formula above involves coinvariants, rather than invariants, for the symmetric groups. By construction, the functor $\mathrm{cot}_{\mathcal{O}}$ naturally factors through a functor $\mathrm{indec}_{\mathcal{O}}^{\mathrm{nil}}$ as follows:
\begin{equation*}
\mathrm{Alg}_{\mathcal{O}} \xrightarrow{\mathrm{indec}^{\mathrm{nil}}_{\mathcal{O}}} \mathrm{coAlg}^{\mathrm{dp},\mathrm{nil}}_{B\mathcal{O}}(\mathcal{C}) \xrightarrow{\mathrm{forget}^{\mathrm{nil}}_{B\mathcal{O}}} \mathcal{C}.
\end{equation*}
We will see in Section \ref{sec:coalgebras} that there is a natural (left adjoint) comparison functor between $\infty$-categories
\begin{equation*}
\mathrm{coAlg}^{\mathrm{dp},\mathrm{nil}}_{B\mathcal{O}}(\mathcal{C}) \rightarrow \mathrm{coAlg}^{\mathrm{dp}}_{B\mathcal{O}}(\mathcal{C})
\end{equation*}
relating these ind-conilpotent coalgebras to our previous setup. The composition of this functor with $\mathrm{indec}^{\mathrm{nil}}$ then reproduces our previous indecomposables functor.

For formal reasons, the indecomposables functor $\mathrm{indec}_{\mathcal{O}}^{\mathrm{nil}}$ admits a right adjoint $\mathrm{prim}^{\mathrm{nil}}_{B\mathcal{O}}$ which we will refer to as \emph{primitives}. The conjecture of Francis--Gaitsgory concerns the class of \emph{pronilpotent} $\mathcal{O}$-algebras, which is the smallest class of $\mathcal{O}$-algebras  closed under limits and containing all nilpotent $\mathcal{O}$-algebras.

\begin{conjecture}[Francis--Gaitsgory \cite{francisgaitsgory}]
\label{conj:francisgaitsgory}
The adjoint functors $\mathrm{indec}_{\mathcal{O}}^{\mathrm{nil}}$ and $\mathrm{prim}^{\mathrm{nil}}_{B\mathcal{O}}$ restrict to an adjoint equivalence from the full subcategory of $\mathrm{Alg}_{\mathcal{O}}(\mathcal{C})$ on pronilpotent algebras to the $\infty$-category $\mathrm{coAlg}^{\mathrm{dp},\mathrm{nil}}_{B\mathcal{O}}(\mathcal{C})$ of conilpotent divided power coalgebras.
\end{conjecture}

We will demonstrate that in fact this conjecture fails in a rather basic case:

\begin{proposition}
\label{prop:counterexample}
Let $k$ be a field of characteristic zero and $\mathcal{C} = \mathrm{Mod}_k$ its (unbounded) derived $\infty$-category. Take $\mathcal{O}$ to be the nonunital commutative operad. Then there exists a pronilpotent $\mathcal{O}$-algebra $X$ such that the unit map
\begin{equation*}
X \rightarrow \mathrm{prim}_{B\mathcal{O}}^{\mathrm{nil}}\mathrm{indec}_{\mathcal{O}}^{\mathrm{nil}} X
\end{equation*}
is \emph{not} an isomorphism. Thus, the indecomposables functor $\mathrm{indec}_{\mathcal{O}}^{\mathrm{nil}}$ is not fully faithful when restricted to the full subcategory of pronilpotent $\mathcal{O}$-algebras and Conjecture \ref{conj:francisgaitsgory} does not hold.
\end{proposition}

An important difference between the contexts of Theorem \ref{thm:completeness} and Conjecture \ref{conj:francisgaitsgory} is the distinction between the two $\infty$-categories of coalgebras featuring in their statements. We can summarize the pros and cons of the two $\infty$-categories $\mathrm{coAlg}^{\mathrm{dp},\mathrm{nil}}_{B\mathcal{O}}(\mathcal{C})$ and $\mathrm{coAlg}^{\mathrm{dp}}_{B\mathcal{O}}(\mathcal{C})$ as follows:
\begin{itemize}
\item[(1)] For $\mathrm{coAlg}^{\mathrm{dp},\mathrm{nil}}_{B\mathcal{O}}(\mathcal{C})$, the cofree coalgebra functor is easily understood explicitly. Indeed, by construction, it is the functor $\mathrm{cot}_\mathcal{O} \circ \mathrm{triv}_{\mathcal{O}}$ for which we gave a formula above. However, in the case where $B\mathcal{O}$ is the commutative (or associative) operad, the $\infty$-category $\mathrm{coAlg}^{\mathrm{dp},\mathrm{nil}}_{B\mathcal{O}}(\mathcal{C})$ is not the usual $\infty$-category of commutative (resp. associative) coalgebras in $\mathcal{C}$, defined as the opposite of the $\infty$-category of commutative (resp. associative) algebras in $\mathcal{C}^{\mathrm{op}}$. More precisely, take $\mathcal{C} = \mathrm{Mod}_k$ for a field $k$ of characteristic zero (so that one can ignore divided powers) and $B\mathcal{O}$ the commutative cooperad. Then the $\infty$-categories $\mathrm{coAlg}^{\mathrm{dp},\mathrm{nil}}_{B\mathcal{O}}(\mathcal{C})$ and $\mathrm{coCAlg}(\mathcal{C})$ are not equivalent; it is a priori not even clear whether the first is (equivalent to) a full subcategory of the second (cf. Remark \ref{rmk:comparingcoalgs}).
\item[(2)] The $\infty$-categories $\mathrm{coAlg}_{B\mathcal{O}}(\mathcal{C})$ and $\mathrm{coAlg}^{\mathrm{dp}}_{B\mathcal{O}}(\mathcal{C})$ reproduce familiar examples for specific choices of $B\mathcal{O}$. Indeed, if $B\mathcal{O}$ is the commutative cooperad, one has an equivalence
\begin{equation*}
\mathrm{coAlg}_{B\mathcal{O}}(\mathcal{C}) \simeq \mathrm{CAlg}(\mathcal{C}^{\mathrm{op}})^{\mathrm{op}}
\end{equation*}
as expected. The downside of these $\infty$-categories is that the corresponding cofree coalgebra functors are produced via the adjoint functor theorem, but cannot be described by an explicit and simple formula (like that for $\mathrm{cot}_\mathcal{O} \circ \mathrm{triv}_{\mathcal{O}}$) in general.
\end{itemize}

We will now explain the relation between Theorem \ref{thm:completeness} and Conjecture \ref{conj:francisgaitsgory}. For an $\mathcal{O}$-algebra $X$, we have two notions of `completion' at our disposal. The first is the nilcompletion, which by Theorem \ref{thm:completeness} can be interpreted as the unit
\begin{equation*}
X \rightarrow \mathrm{prim}_{B\mathcal{O}}\mathrm{indec}_{\mathcal{O}}(X).
\end{equation*}
Recall that $X$ is nilcomplete precisely if this map is an isomorphism. The second notion of completion is the unit
\begin{equation*}
X \rightarrow \mathrm{prim}^{\mathrm{nil}}_{B\mathcal{O}}\mathrm{indec}^{\mathrm{nil}}_{\mathcal{O}}(X)
\end{equation*}
of the adjunction relating $\mathcal{O}$-algebras to conilpotent divided power $B\mathcal{O}$-coalgebras. We say $X$ is \emph{homologically complete} if this map is an isomorphism. Note that the left adjoint functor $\mathrm{coAlg}^{\mathrm{dp},\mathrm{nil}}_{B\mathcal{O}} \rightarrow \mathrm{coAlg}^{\mathrm{dp}}_{B\mathcal{O}}(\mathcal{C})$ induces a natural transformation
\begin{equation*}
\mathrm{prim}^{\mathrm{nil}}_{B\mathcal{O}}\mathrm{indec}^{\mathrm{nil}}_{\mathcal{O}} \rightarrow \mathrm{prim}_{B\mathcal{O}}\mathrm{indec}_{\mathcal{O}}.
\end{equation*}

\begin{remark}
\label{rmk:BKcompletion}
The reason for the term `homologically complete' is that the object $\mathrm{prim}^{\mathrm{nil}}_{B\mathcal{O}}\mathrm{indec}^{\mathrm{nil}}_{\mathcal{O}}(X)$ may be computed as the totalization of the cosimplicial object
\[
\begin{tikzcd}
\mathrm{triv}_\mathcal{O} \mathrm{cot}_{\mathcal{O}}(X) \ar[shift left = .5ex]{r} \ar[shift right = .5ex]{r} & (\mathrm{triv}_\mathcal{O} \mathrm{cot}_{\mathcal{O}})^2(X) \ar{l} \ar[shift left = 1ex]{r} \ar{r} \ar[shift right = 1ex]{r} & \ar[shift left = .5ex]{l} \ar[shift right = .5ex]{l} (\mathrm{triv}_\mathcal{O} \mathrm{cot}_{\mathcal{O}})^3(X) \ar[shift left = 1.5ex]{r}\ar[shift left = .5ex]{r}\ar[shift right = .5ex]{r}\ar[shift right = 1.5ex]{r} & \cdots. \ar[shift left = 1ex]{l} \ar{l} \ar[shift right = 1ex]{l}
\end{tikzcd}
\]
This is exactly the completion in the sense of Bousfield--Kan \cite{bousfieldkan} of $X$ with respect to the homology theory $\mathrm{cot}_{\mathcal{O}}$ (cf. Remark \ref{rmk:homology}). This formula follows from Lemma \ref{lem:primnilresolution}.
\end{remark}

Now consider a trivial algebra $X = \mathrm{triv}_{\mathcal{O}}(M)$, for some $M \in \mathcal{C}$. Then $X$ is homologically complete for formal reasons; the cosimplicial object described above admits a contracting homotopy. However, it is generally not at all clear whether $X$ is nilcomplete. We introduce the following terminology:

\begin{definition}
An operad $\mathcal{O}$ has \emph{good completion} if every trivial $\mathcal{O}$-algebra is nilcomplete.
\end{definition}

When $\mathcal{O}$ satisfies this condition, the two notions of completion discussed above will agree. Indeed, we prove the following in Section \ref{sec:completions}:

\begin{proposition}
\label{prop:twocompletions}
If the operad $\mathcal{O}$ has good completion, then the natural map
\begin{equation*}
\mathrm{prim}^{\mathrm{nil}}_{B\mathcal{O}}\mathrm{indec}^{\mathrm{nil}}_{\mathcal{O}}(X) \rightarrow \mathrm{prim}_{B\mathcal{O}}\mathrm{indec}_{\mathcal{O}}(X)
\end{equation*}
is an equivalence for any $\mathcal{O}$-algebra $X$. In particular, $X$ is homologically complete if and only if it is nilcomplete.
\end{proposition}


Theorem \ref{thm:completeness} yields in particular that the functor $\mathrm{indec}_{\mathcal{O}}$ is fully faithful on the subcategory of $\mathrm{Alg}_{\mathcal{O}}(\mathcal{C})$ consisting of the nilcomplete $\mathcal{O}$-algebras. It follows from Proposition \ref{prop:twocompletions} that if $\mathcal{O}$ has good completion, then this is also the largest full subcategory on which the functor $\mathrm{indec}^{\mathrm{nil}}_{\mathcal{O}}$ is fully faithful. 

Suppose $\mathcal{O}$ has good completion (we will produce a good supply of such $\mathcal{O}$ in Section \ref{sec:goodcompletion}). Then Conjecture \ref{conj:francisgaitsgory}, combined with Proposition \ref{prop:twocompletions}, predicts that every pronilpotent $\mathcal{O}$-algebra $X$ is nilcomplete. We will provide an explicit counterexample to this statement in Section \ref{sec:FGconjecture} and use it to prove Proposition \ref{prop:counterexample}. 

We view Theorem \ref{thm:completeness} as a replacement for Conjecture \ref{conj:francisgaitsgory}. For given $\mathcal{C}$ and $\mathcal{O}$, it reduces the study of Koszul duality for $\mathcal{O}$-algebras to identifying the subcategories of nilcomplete $\mathcal{O}$-algebras and conilcomplete $B\mathcal{O}$-coalgebras. Of course this can still be a challenging problem in practice. In Section \ref{sec:examples} we review some cases of Conjecture \ref{conj:francisgaitsgory} treated in the literature and see that they fit the setting of Theorem \ref{thm:completeness}. In particular, they can all be deduced by identifying certain classes of complete $\mathcal{O}$-algebras and cocomplete $B\mathcal{O}$-coalgebras in specific examples.


\begin{remark}
\label{rmk:hypotheses}
We briefly remark on our hypotheses on the $\infty$-category $\mathcal{C}$. The assumption that the tensor product commutes with colimits in each variable is rather harmless, since it is satisfied in all of the examples relevant to us, but might not quite be necessary. In our arguments about coalgebras we will not have the luxury of a tensor product commuting with (cosifted) limits and often have to introduce more elaborate methods; it is plausible that these could be dualized to the setting of algebras in order to weaken the hypothesis on colimits and tensor products, only assuming that the tensor product is exact in each variable. The presentability of $\mathcal{C}$ will be relevant when applying the adjoint functor theorem in various places; e.g., we only see how to obtain the construction of cofree coalgebras from such an abstract result. Thus we do not see how to remove the hypothesis that $\mathcal{C}$ is presentable.
\end{remark}

\section{Bar-cobar duality for operads and cooperads}
\label{sec:barcobar}

In this section we briefly review the notions of operad and cooperad in $\mathcal{C}$, as well as the bar and cobar constructions. The main result is Theorem \ref{thm:barcobaroperads}, stating that these constructions implement an equivalence between the $\infty$-categories of reduced operads and cooperads in $\mathcal{C}$.

Write $\mathrm{Fin}^{\cong}$ for the $\infty$-category of non-empty finite sets and isomorphisms between them, and
\begin{equation*}
\mathrm{SSeq}(\mathcal{C}) := \mathrm{Fun}(\mathrm{Fin}^{\cong}, \mathcal{C})
\end{equation*}
for the $\infty$-category of symmetric sequences in $\mathcal{C}$. This $\infty$-category has a monoidal structure given by the \emph{composition product}, which is loosely described by the following formula, for symmetric sequences $A$ and $B$ and a finite set $I$:
\begin{equation*}
(A \circ B)(I) \cong \bigoplus_{E \in \mathrm{Equiv}(I)} A(I/E) \otimes \bigl(\otimes_{J \in I/E} B(J)\bigr).
\end{equation*}
Here the sum is over all equivalence relations $E$ on $I$ and $I/E$ denotes the set of $E$-equivalence classes in $I$. A precise construction of the composition product is given in \cite{haugsengsymseq} or \cite[Section 4.1.2]{brantnerthesis}. To any symmetric sequence $A$ one can associate a functor
\begin{equation*}
\mathrm{Sym}_A \colon \mathcal{C} \rightarrow \mathcal{C}\colon X \mapsto \bigoplus_{n \geq 1} (A(n) \otimes X^{\otimes n})_{h\Sigma_n}.
\end{equation*}
Here $A(n)$ should be read as $A(\{1, \ldots, n\})$. The crucial property of the composition product (and essentially the motivation for its definition) is that the assignment
\begin{equation*}
\mathrm{SSeq}(\mathcal{C}) \rightarrow \mathrm{Fun}(\mathcal{C},\mathcal{C}) \colon A \mapsto \mathrm{Sym}_A
\end{equation*}
can be made into a monoidal functor. In other words, it takes the composition product of symmetric sequences $A \circ B$ to the composition of functors $\mathrm{Sym}_A \circ \mathrm{Sym}_B$ in a coherent way.

\begin{remark}
\label{rmk:compositionproduct}
An explicit construction of the composition product (and the one described in \cite{brantnerthesis}) is as follows. The Day convolution tensor product equips the $\infty$-category $\mathrm{SSeq}(\mathrm{Sp})$ of symmetric sequences of spectra with a symmetric monoidal structure (different from the composition product); for symmetric sequences $X$ and $Y$ and a finite set $I$, a formula for their tensor product is
\begin{equation*}
(X \otimes Y)(I) \cong \bigoplus_{I = I_1 \amalg I_2} X(I_1) \otimes Y(I_2),
\end{equation*}
where the right-hand side uses the symmetric monoidal structure of $\mathrm{Sp}$. With this tensor product, the $\infty$-category $\mathrm{SSeq}(\mathrm{Sp})$ is the universal presentably symmetric monoidal stable $\infty$-category on a single object. More precisely, for any other such $\mathcal{C}$, evaluation at the one-element set determines an equivalence of $\infty$-categories
\begin{equation*}
\mathrm{Fun}^{\otimes,L}(\mathrm{SSeq}(\mathrm{Sp}), \mathcal{C}) \cong \mathcal{C}.
\end{equation*}
The notation $\mathrm{Fun}^{\otimes,L}$ indicates symmetric monoidal functors preserving all colimits. In particular, we find an equivalence
\begin{equation*}
\mathrm{Fun}^{\otimes,L}(\mathrm{SSeq}(\mathrm{Sp}), \mathrm{SSeq}(\mathrm{Sp})) \cong \mathrm{SSeq}(\mathrm{Sp}).
\end{equation*}
Composition of functors clearly makes the left-hand side into a monoidal $\infty$-category. The corresponding monoidal structure on the right-hand side is the opposite of the composition product of symmetric sequences of spectra. In other words, the symmetric monoidal endofunctor of $\mathrm{SSeq}(\mathrm{Sp})$ associated to a symmetric sequence $X$ is precisely $- \circ X$. To get the composition product on symmetric sequences in a general presentably symmetric monoidal stable $\infty$-category $\mathcal{C}$ one follows the same construction, but working in $\mathcal{C}$-linear presentably symmetric monoidal $\infty$-categories, rather than just stable ones.
\end{remark}

With the composition product in place, we can define operads (resp. cooperads) as monoids in the $\infty$-category of symmetric sequences:

\begin{definition}
The $\infty$-category of operads (resp. cooperads) in $\mathcal{C}$ is that of associative algebras (resp. coalgebras) in symmetric sequences with respect to the composition product:
\begin{equation*}
\mathrm{Op}(\mathcal{C}) := \mathrm{Alg}(\mathrm{SSeq}(\mathcal{C})) \quad\quad (\text{resp.} \quad \mathrm{coOp}(\mathcal{C}) := \mathrm{coAlg}(\mathrm{SSeq}(\mathcal{C}))).
\end{equation*}
\end{definition}

\begin{remark}
A comparison between a definition of operads using symmetric sequences and other definitions of (enriched) $\infty$-operads in the literature is given in \cite{haugsengsymseq}. With respect to the definition of the composition product we use here, a comparison will appear in \cite{brantnerheutsKD}.
\end{remark}

Generally, augmented associative algebras and coalgebras in a monoidal $\infty$-category are related by bar-cobar duality. A very general form of this principle is explained in \cite[Section 4.3]{dagx} and \cite[Section 3.4]{BCN}. Essentially, if $\mathcal{D}$ is a pointed monoidal $\infty$-category which admits geometric realizations and totalizations, then there is an adjoint pair
\[
\begin{tikzcd}
\mathrm{Alg}^{\mathrm{aug}}(\mathcal{D}) \ar[shift left]{r}{\mathrm{Bar}} & \mathrm{coAlg}^{\mathrm{aug}}(\mathcal{D}). \ar[shift left]{l}{\mathrm{Cobar}}
\end{tikzcd}
\]
For an augmented associative algebra $X \in \mathrm{Alg}(\mathcal{D})$, its bar construction $\mathrm{Bar}(X)$ is the relative tensor product $\mathbf{1} \otimes_X \mathbf{1}$, i.e., the geometric realization of the simplicial object
\[
\begin{tikzcd}
\cdots \ar[shift left = 1.5ex]{r}\ar[shift left = .5ex]{r}\ar[shift right = .5ex]{r}\ar[shift right = 1.5ex]{r} & X \otimes X \ar[shift left = 1ex]{r}\ar{r}\ar[shift right = 1ex]{r} \ar[shift left = 1ex]{l}\ar{l}\ar[shift right = 1ex]{l} & X \ar[shift left = .5ex]{r}\ar[shift right =.5ex]{r}  \ar[shift left = .5ex]{l}\ar[shift right =.5ex]{l} & \mathbf{1}. \ar{l}
\end{tikzcd}
\]
Here $\mathbf{1}$ denotes the monoidal unit of $\mathcal{D}$. The face maps in this simplicial object are built from the multiplication of $X$ and its augmentation, whereas the degeneracies use the unit. Dually, the cobar construction of a coalgebra is the totalization of a cosimplicial object constructed similarly. In the particular case where $\mathcal{D} = \mathrm{SSeq}(\mathcal{C})$ and the monoidal structure is the composition product, we obtain an adjunction
\[
\begin{tikzcd}
\mathrm{Op}^{\mathrm{aug}}(\mathcal{C}) \ar[shift left]{r}{B} & \mathrm{coOp}^{\mathrm{aug}}(\mathcal{C}) \ar[shift left]{l}{C}
\end{tikzcd}
\]
between the $\infty$-categories of augmented operads and cooperads, where we have abbreviated the notations $\mathrm{Bar}$ and $\mathrm{Cobar}$ to $B$ and $C$ respectively. In this context $\mathbf{1}$ is the unit symmetric sequence, which has the monoidal unit of $\mathcal{C}$ in degree 1 and 0 everywhere else. All of the operads (and cooperads) we will deal with are \emph{reduced}, meaning they have $\mathcal{O}(1) = \mathbf{1}$ (and $\mathcal{Q}(1) = \mathbf{1}$). We will implicitly take their augmentations to be the obvious projections to the degree 1 term.

The fact that $\mathcal{C}$ is stable turns out to imply that bar-cobar duality for operads and cooperads in $\mathcal{C}$ works particularly well. The following result is quoted without proof in \cite{francisgaitsgory} and proved in the special case $\mathcal{C} = \mathrm{Sp}$ by Ching in \cite{chingbar}. His argument also works in the general case to prove the following:

\begin{theorem}
\label{thm:barcobaroperads}
The bar-cobar adjunction restricts to an equivalence between the $\infty$-categories of reduced operads and cooperads in $\mathcal{C}$. In particular, for a reduced operad $\mathcal{O}$ and reduced cooperad $\mathcal{Q}$, the unit and counit maps
\begin{equation*}
\mathcal{O} \rightarrow CB\mathcal{O} \quad \text{and} \quad BC\mathcal{Q} \rightarrow \mathcal{Q}
\end{equation*}
are isomorphisms.
\end{theorem}
\begin{remark}
This result and its version for right modules will be discussed in more detail in \cite{brantnerheutsKD}. We include an outline of proof below for the reader's convenience and stress again that the essence of the argument is due to Ching \cite{chingbar}.
\end{remark}
\begin{proof}[Outline of proof]
We prove that the first map is an isomorphism, the argument for the second is analogous. For any $n \geq 1$, write $\tau_n \mathcal{O}$ for the symmetric sequence with $\tau_n\mathcal{O}(k) = \mathcal{O}(k)$ for $k \leq n$ and $\tau_n\mathcal{O}(k) = 0$ for $k > n$. (We will have more to say about these `truncations' below.) For the purposes of this proof, $\tau_n\mathcal{O}$ can be regarded as a right $\mathcal{O}$-module in the evident way. We will use the extension of Lurie's bar-cobar duality to right modules and comodules as explained in \cite[Theorem 3.26]{BCN} or \cite[Proposition 4.1.22]{blansblom}. For a right module $M$ we may form the bar construction $B(M,\mathcal{O},\mathbf{1})$, defined as the geometric realization of the simplicial object
\begin{equation*}
B_\bullet(M,\mathcal{O},\mathbf{1}) := M \circ \mathcal{O}^{\circ \bullet}.
\end{equation*}
Its face maps are determined by the multiplication maps of $\mathcal{O}$, the right $\mathcal{O}$-action on $M$, and the left $\mathcal{O}$-action on $\mathbf{1}$. This bar construction is naturally a right $B\mathcal{O}$-comodule. Similarly, for any right $B\mathcal{O}$-comodule $N$ one can form a cobar construction $C(N,B\mathcal{O},1)$, which is a right $\mathcal{O}$-module. This cobar construction for comodules is right adjoint to the bar construction for modules. We claim that the unit map
\begin{equation*}
\tau_n \mathcal{O} \rightarrow C(B(\tau_n \mathcal{O}, \mathcal{O}, \mathbf{1}), B\mathcal{O}, \mathbf{1})
\end{equation*}
is an isomorphism for every $n$. The conclusion follows from this by taking the inverse limit over $n$. Indeed, the left-hand side becomes $\mathcal{O}$, whereas on the right-hand side one has $\varprojlim_n B(\tau_n\mathcal{O},\mathcal{O},\mathbf{1}) = \mathbf{1}$, using that in every fixed degree of these symmetric sequences the limit is eventually constant. We establish our claim by induction on $n$. For $n=1$, we are considering the map
\begin{equation*}
\mathbf{1} = \tau_1\mathcal{O} \rightarrow C(B\mathcal{O}, B\mathcal{O},\mathbf{1}).
\end{equation*}
The cosimplicial object $C^{\bullet}(B\mathcal{O},B\mathcal{O},\mathbf{1})$ admits a contracting homotopy and hence the map is an isomorphism. For $n>1$, consider the diagram of fiber sequences
\[
\begin{tikzcd}
\mathcal{O}(n) \ar{d}\ar{r} & C(B(\mathcal{O}(n), \mathcal{O}, \mathbf{1}), B\mathcal{O},\mathbf{1}) \ar{d} \\
\tau_n\mathcal{O}\ar{d} \ar{r} & C(B(\tau_n \mathcal{O}, \mathcal{O}, \mathbf{1}), B\mathcal{O}, \mathbf{1}) \ar{d} \\
\tau_{n-1}\mathcal{O} \ar{r} & C(B(\tau_{n-1} \mathcal{O}, \mathcal{O}, \mathbf{1}), B\mathcal{O}, \mathbf{1}).
\end{tikzcd}
\]
Here $\mathcal{O}(n)$ is considered as a symmetric sequence concentrated in degree $n$. Its right $\mathcal{O}$-module structure is trivial (for evident degree reasons); hence we find
\begin{equation*}
B(\mathcal{O}(n), \mathcal{O}, \mathbf{1}) \cong \mathcal{O}(n) \circ B(\mathbf{1}, \mathcal{O},\mathbf{1}) = \mathcal{O}(n) \circ B\mathcal{O}.
\end{equation*}
But then another contracting homotopy shows that
\begin{equation*}
C(B(\mathcal{O}(n), \mathcal{O}, \mathbf{1}), B\mathcal{O},\mathbf{1}) \cong \mathcal{O}(n)
\end{equation*}
and the top horizontal arrow is an isomorphism. The bottom horizontal arrow is an isomorphism by the inductive hypothesis and hence the middle arrow is an isomorphism as well, completing the inductive step.
\end{proof}

\section{Truncations of operads and cooperads}
\label{sec:truncations}

We will discuss two natural ways of `approximating' an operad $\mathcal{O}$ from the left and from the right by certain filtrations (and similarly for a cooperad $\mathcal{Q}$). These approximations will be a very useful tool to analyze the $\infty$-categories of $\mathcal{O}$-algebras and $\mathcal{Q}$-coalgebras later on.

We write $\mathrm{Fin}^{\cong}_{\leq n}$ for the category of non-empty finite sets of cardinality at most $n$ and bijections between them. (We use non-empty sets, since we are only interested in reduced operads and cooperads.) We define the $\infty$-category of \emph{$n$-truncated symmetric sequences} in $\mathcal{C}$ to be 
\begin{equation*}
\mathrm{SSeq}_{\leq n}(\mathcal{C}) := \mathrm{Fun}(\mathrm{Fin}^{\cong}_{\leq n}, \mathcal{C}).
\end{equation*}
Like the $\infty$-category of all symmetric sequences it carries a composition product, making it a monoidal $\infty$-category. In fact, this structure is essentially uniquely determined by requiring the evident restriction functor
\begin{equation*}
\rho_n\colon \mathrm{SSeq}(\mathcal{C}) \rightarrow \mathrm{SSeq}_{\leq n}(\mathcal{C})
\end{equation*}
to be monoidal. It is clear that $\rho_n$ preserves limits and colimits. We define the $\infty$-category of (nonunital) \emph{$n$-truncated operads in $\mathcal{C}$} to be that of algebra objects in $n$-truncated symmetric sequences:
\begin{equation*}
\mathrm{Op}_{\leq n}(\mathcal{C}) := \mathrm{Alg}(\mathrm{SSeq}_{\leq n}(\mathcal{C})).
\end{equation*}
Similarly, we define
\begin{equation*}
\mathrm{coOp}_{\leq n}(\mathcal{C}) := \mathrm{coAlg}(\mathrm{SSeq}_{\leq n}(\mathcal{C}))
\end{equation*}
to be the $\infty$-category of $n$-truncated cooperads. The restriction $\rho_n$ induces corresponding functors
\begin{equation*}
\mathrm{Op}(\mathcal{C}) \rightarrow \mathrm{Op}_{\leq n}(\mathcal{C}) \quad \text{and} \quad \mathrm{coOp}(\mathcal{C}) \rightarrow \mathrm{coOp}_{\leq n}(\mathcal{C}).
\end{equation*}
These two restrictions admit both a left and a right adjoint, as we will now explain. The functor
\begin{equation*}
\zeta_n\colon \mathrm{SymSeq}_{\leq n}(\mathcal{C}) \rightarrow \mathrm{SymSeq}(\mathcal{C})
\end{equation*}
that extends by zero above degree $n$ is both left and right adjoint to $\rho_n$. Since $\rho_n$ is monoidal, it follows that $\zeta_n$ admits both a lax monoidal and an oplax monoidal structure. The lax monoidal structure induces a functor on the level of algebras
\begin{equation*}
\tau_n\colon \mathrm{Op}_{\leq n}(\mathcal{C}) \rightarrow \mathrm{Op}(\mathcal{C})
\end{equation*}
which is right adjoint to the restriction $\mathrm{Op}(\mathcal{C}) \rightarrow \mathrm{Op}_{\leq n}(\mathcal{C})$, whereas the oplax monoidal structure gives a functor
\begin{equation*}
\tau^n\colon \mathrm{coOp}_{\leq n}(\mathcal{C}) \rightarrow \mathrm{coOp}(\mathcal{C})
\end{equation*}
which is left adjoint to the restriction $\mathrm{coOp}(\mathcal{C}) \rightarrow \mathrm{coOp}_{\leq n}(\mathcal{C})$. Thus for an operad $\mathcal{O}$, the operad $\tau_n \mathcal{O}$ (where we have suppressed the restriction $\rho_n$ from the notation) is its `$n$-truncation' in which all operations of arity greater than $n$ have been set to zero. It receives a map $\mathcal{O} \rightarrow \tau_n \mathcal{O}$. Dually, a cooperad $\mathcal{Q}$ has an $n$-truncation equipped with a natural map $\tau^n \mathcal{Q} \rightarrow \mathcal{Q}$.

On the other hand, since the restriction $\mathrm{Op}(\mathcal{C}) \rightarrow \mathrm{Op}_{\leq n}(\mathcal{C})$ preserves limits and filtered colimits, the adjoint functor theorem \cite[Corollary 5.5.2.9]{HTT} guarantees the existence of a left adjoint
\begin{equation*}
\varphi_n\colon \mathrm{Op}_{\leq n}(\mathcal{C}) \rightarrow \mathrm{Op}(\mathcal{C}).
\end{equation*}
The operad $\varphi_n\mathcal{O}$ (again suppressing $\rho_n$ from the notation) comes with a natural map to $\mathcal{O}$ which is an isomorphism in arities up to $n$. The operations of $\varphi_n\mathcal{O}$ in higher arities are not zero, but rather `freely generated' by those in arities up to $n$. Dually, the restriction $\mathrm{coOp}(\mathcal{C}) \rightarrow \mathrm{coOp}_{\leq n}(\mathcal{C})$ has a right adjoint
\begin{equation*}
\varphi^n\colon \mathrm{coOp}_{\leq n}(\mathcal{C}) \rightarrow \mathrm{coOp}(\mathcal{C}),
\end{equation*}
where a cooperad of the form $\varphi^n\mathcal{Q}$ has cooperations cofreely generated by those in arities up to $n$.

\begin{remark}
\label{rmk:freelygeneratedops}
A version of the left adjoint functor $\varphi_n$ is studied in \cite[Appendix C]{heutsgoodwillie}. Although we will not need them in this paper, it is not hard to give explicit formulas for the objects $(\varphi_n \mathcal{O})(k)$. To be precise, we have (cf. \cite[Lemma C.20]{heutsgoodwillie})
\begin{equation*}
(\varphi_n \mathcal{O})(k) \simeq \varinjlim_{T \in \mathrm{Dec}_n(C_k)^{\mathrm{op}}} \bigotimes_{v \in V(T)} \mathcal{O}(\mathrm{in}(v)).
\end{equation*}
Here $C_k$ denotes the $k$-corolla (the rooted tree with $k$ leaves and only one internal vertex) and the category $\mathrm{Dec}_n(C_k)$ is the diagram of all `$n$-decompositions' $C_k \rightarrow T$ of the $k$-corolla, meaning trees $T$ with $k$ leaves labelled by those of $C_k$, satisfying the condition that any internal vertex of $T$ has at most $n$ input edges. The tensor product on the right is over the internal vertices $v$ of $T$, where $\mathrm{in}(v)$ denotes the set of input edges of such a vertex. This formula makes precise the idea that the $k$-ary operations of $\varphi_n \mathcal{O}$ are built by considering all the ways in which a $k$-ary operation can be obtained as a composition of operations of $\mathcal{O}$ of arity $\leq n$.
\end{remark}

The bar and cobar constructions interact well with the truncations just discussed, in the sense of the proposition below. The first isomorphism also appears as Theorem 3.9 of \cite{amabel}. We remind the reader that we tacitly take all operads and cooperads to be reduced in what follows.

\begin{proposition}
\label{prop:bartruncation}
For an operad $\mathcal{O}$ there are natural isomorphisms
\begin{equation*}
B(\tau_n \mathcal{O}) \cong \varphi^n B\mathcal{O}, \quad B(\varphi_n\mathcal{O}) \cong \tau^nB\mathcal{O}
\end{equation*}
and dually, for a cooperad $\mathcal{Q}$,
\begin{equation*}
C(\tau^n\mathcal{Q}) \cong \varphi_n C\mathcal{Q}, \quad C(\varphi^n\mathcal{Q}) \cong \tau_nC\mathcal{Q}.
\end{equation*}
\end{proposition}
\begin{proof}
The first isomorphism (resp. the second) follows by taking right (resp. left) adjoints in the commutative square
\[
\begin{tikzcd}
\mathrm{Op}(\mathcal{C}) \ar{r}{\rho_n} & \mathrm{Op}_{\leq n}(\mathcal{C}) \\
\mathrm{coOp}(\mathcal{C}) \ar{r}{\rho_n}\ar{u}{C} & \mathrm{coOp}_{\leq n}(\mathcal{C}). \ar{u}{C} 
\end{tikzcd}
\]
Here we use that the cobar construction $C$ implements an equivalence between cooperads and operads in $\mathcal{C}$, so that its inverse $B$ serves both as a left and a right adjoint. The second set of isomorphisms follows similarly by taking adjoints in the commutative square
\[
\begin{tikzcd}
\mathrm{Op}(\mathcal{C}) \ar{r}{\rho_n}\ar{d}{B} & \mathrm{Op}_{\leq n}(\mathcal{C}) \ar{d}{B} \\
\mathrm{coOp}(\mathcal{C}) \ar{r}{\rho_n} & \mathrm{coOp}_{\leq n}(\mathcal{C}). 
\end{tikzcd}
\] 
\end{proof}

\begin{remark}
\label{rmk:formulabarconstr}
It is a useful consequence of Proposition \ref{prop:bartruncation} that one can write the following formula for the terms of the bar construction of $\mathcal{O}$:
\begin{equation*}
B\mathcal{O}(n) \cong \Sigma \mathrm{cof}\bigl((\varphi_{n-1}\mathcal{O})(n) \rightarrow \mathcal{O}(n)\bigr).
\end{equation*}
This formula follows by considering the pushout square of $n$-truncated operads
\[
\begin{tikzcd}
\rho_n(\varphi_{n-1}\mathcal{O})(n) \ar{r}\ar{d} & \rho_n\varphi_{n-1}\mathcal{O} \ar{d} \\
\mathcal{O}(n) \ar{r} & \rho_n\mathcal{O}.
\end{tikzcd}
\]
Here $\mathcal{O}(n)$ and $(\varphi_{n-1}\mathcal{O})(n)$ denote the corresponding symmetric sequences concentrated in degree $n$, considered as operads. (As $n$-truncated operads, they are both free and trivial. Also, to really make them operads one should add a unit term $\mathbf{1}$ in arity 1, but it will not matter for our argument here.) The bar constructions on these free operads are easily seen to be $\Sigma\mathcal{O}(n)$ and $(\varphi_{n-1}\mathcal{O})(n)$, interpreted as (trivial) cooperads concentrated in degree $n$. (More generally, applying the bar construction to a free operad $\mathrm{free}(A)$ gives the trivial cooperad on the suspended symmetric sequence $\Sigma A$, cf. \cite[Remark 5.2.2.14]{higheralgebra}.) Hence applying $B$ to the pushout square above gives a pushout of cooperads
\[
\begin{tikzcd}
\Sigma\bigl((\varphi_{n-1}\mathcal{O})(n)\bigr) \ar{r}\ar{d} & \tau^{n-1}B\mathcal{O} \ar{d} \\
\Sigma\mathcal{O}(n) \ar{r} & \tau^n B\mathcal{O}.
\end{tikzcd}
\]
Evaluating in degree $n$ and using $(\tau^{n-1}B\mathcal{O})(n) = 0$ gives the claimed formula for $B\mathcal{O}(n)$.
\end{remark}

\section{Algebras}
\label{sec:algebras}

In this section we review the $\infty$-category of $\mathcal{O}$-algebras and provide a basic structural result, Theorem  \ref{thm:algdecomposition}, which essentially provides an obstruction theory for $\mathcal{O}$-algebra structures. It is different from (and in some sense Koszul dual to) the usual `arity truncation' used in many other places. 


Recall that the functor $\mathrm{SSeq}(\mathcal{C}) \to \mathrm{End}(\mathcal{C})\colon A \mapsto \mathrm{Sym}_A$ is monoidal. In particular, it sends an operad $\mathcal{O}$ to a monad $\mathrm{Sym}_{\mathcal{O}}$, so that we can make the following definition:

\begin{definition}
Let $\mathcal{O}$ be an operad in $\mathcal{C}$. Then an \emph{$\mathcal{O}$-algebra} is an algebra for the monad $\mathrm{Sym}_{\mathcal{O}}$. We write $\mathrm{Alg}_{\mathcal{O}}(\mathcal{C}) := \mathrm{Alg}_{\mathrm{Sym}_{\mathcal{O}}}(\mathcal{C})$ for the $\infty$-category of $\mathcal{O}$-algebras.
\end{definition}

For a symmetric sequence $A$, we define the following `extended power' functors on $\mathcal{C}$:
\begin{equation*}
D_n^A(X) := (A(n) \otimes X^{\otimes n})_{h\Sigma_n}, \quad D^n_A(X) := (A(n) \otimes X^{\otimes n})^{h\Sigma_n}.
\end{equation*}
For a functor $F\colon \mathcal{C} \rightarrow \mathcal{C}$, we define the $\infty$-category $\mathrm{alg}_F(\mathcal{C})$ of \emph{$F$-algebras} to be the $\infty$-category of objects $X$ equipped with a map $F(X) \rightarrow X$. To be precise, it is defined by the following pullback square of simplicial sets:
\[
\begin{tikzcd}
\mathrm{alg}_F(\mathcal{C}) \ar{r}\ar{d} & \mathcal{C}^{\Delta^1} \ar{d}{(\mathrm{ev}_0,\mathrm{ev}_1)} \\
\mathcal{C} \ar{r}{(F, \mathrm{id})} & \mathcal{C} \times \mathcal{C}. 
\end{tikzcd}
\]
Note that we use lower-case notation $\mathrm{alg}_F(\mathcal{C})$ when speaking about algebras for a functor and upper-case notation $\mathrm{Alg}_{\mathcal{O}}(\mathcal{C})$ for algebras for an operad (or monad). 

Observe that if $X$ has the structure of an $\mathcal{O}$-algebra, then in particular it comes equipped with a map
\begin{equation*}
D_n^{\mathcal{O}}(X) \rightarrow X
\end{equation*}
for each $n \geq 2$. These maps are used to define the horizontal arrows in the square of the following theorem, which provides a very useful decomposition of the $\infty$-category of $\mathcal{O}$-algebras:

\begin{theorem}
\label{thm:algdecomposition}
For each $n \geq 2$, the commutative square of $\infty$-categories
\[
\begin{tikzcd}
\mathrm{Alg}_{\varphi_n\mathcal{O}}(\mathcal{C}) \ar{r}\ar{d} & \mathrm{alg}_{D_n^{\mathcal{O}}}(\mathcal{C}) \ar{d} \\
\mathrm{Alg}_{\varphi_{n-1}\mathcal{O}}(\mathcal{C}) \ar{r} & \mathrm{alg}_{D_n^{\varphi_{n-1}\mathcal{O}}}(\mathcal{C})
\end{tikzcd}
\]
is a pullback. Furthermore, the natural map
\begin{equation*}
\mathrm{Alg}_{\mathcal{O}}(\mathcal{C}) \rightarrow \varprojlim_n \mathrm{Alg}_{\varphi_n\mathcal{O}}(\mathcal{C})
\end{equation*}
is an equivalence of $\infty$-categories.
\end{theorem}

\begin{remark}
\label{rmk:decomposition}
Let us describe the contents of the theorem a bit more informally. First of all, the pullback square of the theorem shows that to upgrade a $\varphi_{n-1}\mathcal{O}$-algebra $X$ to a $\varphi_n\mathcal{O}$-algebra, it suffices to equip $X$ with a multiplication map
\begin{equation*}
\mu_n\colon (\mathcal{O}(n) \otimes X^{\otimes n})_{h\Sigma_n} \rightarrow X
\end{equation*}
which is compatible with the (already specified) map 
\begin{equation*}
\widetilde{\mu}_n\colon ((\varphi_{n-1}\mathcal{O})(n) \otimes X^{\otimes n})_{h\Sigma_n} \rightarrow X.
\end{equation*}
In view of Remark \ref{rmk:freelygeneratedops}, this compatibility is precisely that of $\mu_n$ with all the `lower' multiplication maps $\mu_k$ for $k<n$. The last part of the theorem can then be interpreted as saying that to equip an object $X$ of $\mathcal{C}$ with an $\mathcal{O}$-algebra structure amounts to inductively equipping it with multiplication maps $\mu_n$ for $n \geq 2$, each one compatible with the previously specified $\mu_k$ for $k<n$.
\end{remark}

\begin{proof}[Proof of Theorem \ref{thm:algdecomposition}]
The theorem is a direct consequence of the following three claims:
\begin{itemize}
\item[(a)] The functor
\begin{equation*}
\mathrm{Op}(\mathcal{C})^{\mathrm{op}} \rightarrow (\mathrm{Cat}_\infty)_{/\mathcal{C}}\colon \mathcal{O} \mapsto \mathrm{Alg}_{\mathcal{O}}(\mathcal{C})
\end{equation*}
preserves limits. Here $\mathrm{Alg}_{\mathcal{O}}(\mathcal{C})$ maps to $\mathcal{C}$ by the forgetful functor.
\item[(b)] The square
\[
\begin{tikzcd}
\mathrm{free}\bigl((\varphi_{n-1}\mathcal{O})(n)\bigr) \ar{r}\ar{d} & \varphi_{n-1} \mathcal{O} \ar{d} \\
\mathrm{free}(\mathcal{O}(n)) \ar{r} & \varphi_n\mathcal{O}
\end{tikzcd}
\]
is a pushout in the $\infty$-category $\mathrm{Op}(\mathcal{C})$. (Here $\mathrm{free}(\mathcal{O}(n))$ denotes the free operad generated by $\mathcal{O}(n)$, the latter thought of as a symmetric sequence concentrated in degree $n$.) Furthermore, the natural map
\begin{equation*}
\varinjlim_n \varphi_n\mathcal{O} \rightarrow \mathcal{O}
\end{equation*}
is an isomorphism in $\mathrm{Op}(\mathcal{C})$.
\item[(c)] For a symmetric sequence $A$, generating a free operad $\mathrm{free}(A)$, the natural functor
\begin{equation*}
\mathrm{Alg}_{\mathrm{free}(A)}(\mathcal{C}) \rightarrow \mathrm{alg}_{\mathrm{Sym}_A}(\mathcal{C})
\end{equation*}
from algebras for the operad $\mathrm{free}(A)$ to algebras for the functor $\mathrm{Sym}_A$ is an equivalence. Recall that $\mathrm{Sym}_A$ denotes the functor
\begin{equation*}
\mathrm{Sym}_A(X) = \bigoplus_{n \geq 1} D_n^A(X).
\end{equation*}
\end{itemize}

Let us begin with claim (c). We will use the fact that $\mathrm{Sym}\colon \mathrm{SSeq}(\mathcal{C}) \to \mathrm{End}(\mathcal{C})$ preserves colimits and is monoidal. First we observe that the monad $\mathrm{Sym}_{\mathrm{free}(A)}$ is in fact the free monad generated by the functor $\mathrm{Sym}_A$. Indeed, this can be seen from the explicit construction of free algebras in \cite[Appendix B]{BCN}, which applies to the $\infty$-categories $\mathrm{SSeq}(\mathcal{C})$ and the $\infty$-category $\mathrm{End}^{\mathrm{filt}}(\mathcal{C})$ of endofunctors preserving filtered colimits. From \cite[Remark B.4]{BCN} one sees that algebras for the free monad on the functor $\mathrm{Sym}_A$ are indeed the same thing as algebras for the functor $\mathrm{Sym}_A$.

We now prove claim (a). Write $\mathrm{Mnd}^{\mathrm{acc}}(\mathcal{C})$ for the $\infty$-category of accessible monads on $\mathcal{C}$. The assignment
\begin{equation*}
\mathrm{Mnd}^{\mathrm{acc}}(\mathcal{C})^{\mathrm{op}} \rightarrow (\mathrm{Cat}_\infty)_{/\mathcal{C}}\colon T \rightarrow \bigl(\mathrm{Alg}_T(\mathcal{C}) \xrightarrow{\mathrm{forget}} \mathcal{C}\bigr),
\end{equation*}
sending a monad on $\mathcal{C}$ to its $\infty$-category of algebras, preserves limits. Indeed, it can be factored as 
\begin{equation*}
\mathrm{Mnd}^{\mathrm{acc}}(\mathcal{C})^{\mathrm{op}} \rightarrow (\mathcal{P}\mathrm{r}_\infty^{R})_{/\mathcal{C}} \rightarrow (\mathrm{Cat}_\infty)_{/\mathcal{C}}
\end{equation*}
where $\mathcal{P}\mathrm{r}_\infty^{R}$ is the subcategory of $\mathrm{Cat}_\infty$ of presentable $\infty$-categories and right adjoint functors. (Right adjoints between presentable $\infty$-categories are precisely the accessible functors that preserve limits.) The second arrow is the inclusion, which preserves limits by \cite[Theorem 5.5.3.18]{HTT}. The first functor has a left adjoint, namely the functor that sends a right adjoint $R\colon \mathcal{D} \rightarrow \mathcal{C}$ to the accessible monad $RL$ on $\mathcal{C}$, with $L$ denoting the left adjoint of $R$ (see \cite[Theorem 1.3]{heine}). 
Thus, it will suffice to show that the assignment
\begin{equation*}
\mathrm{Op}(\mathcal{C}) \rightarrow \mathrm{Mnd}(\mathcal{C})\colon \mathcal{O} \mapsto \mathrm{Sym}_{\mathcal{O}}
\end{equation*}
preserves colimits. It will be convenient to further restrict our attention to functors and monads preserving sifted colimits; we add a superscript $\Sigma$ to indicate these subcategories, e.g., we write $\mathrm{Mnd}^{\Sigma}(\mathcal{C})$ for the $\infty$-category of monads on $\mathcal{C}$ preserving sifted colimits. In the commutative square
\[
\begin{tikzcd}
\mathrm{Op}(\mathcal{C}) \ar{r}{\mathrm{Sym}}\ar{d}{\mathrm{forget}} & \mathrm{Mnd}^{\Sigma}(\mathcal{C}) \ar{d}{\mathrm{forget}} \\ 
\mathrm{SSeq}(\mathcal{C}) \ar{r}{\mathrm{Sym}} & \mathrm{Fun}^{\Sigma}(\mathcal{C},\mathcal{C}), 
\end{tikzcd}
\]
the vertical arrows preserve and create sifted colimits. The bottom horizontal arrow preserves arbitrary colimits of symmetric sequences; we deduce that the top horizontal arrow preserves sifted colimits. It remains to show that it preserves finite coproducts. Now any operad $\mathcal{O}$ is canonically a colimit of a simplicial object consisting of free operads (in particular, this is a sifted colimit). Hence it suffices to show that the assignment $\mathcal{O} \mapsto \mathrm{Sym}_{\mathcal{O}}$ preserves coproducts of free operads. For symmetric sequences $A$ and $B$, we have $\mathrm{free}(A) \amalg \mathrm{free}(B) \cong \mathrm{free}(A \oplus B)$. Using claim (c) we can identify the canonical functor
\begin{equation*}
\mathrm{Alg}_{\mathrm{free}(A\oplus B)}(\mathcal{C}) \rightarrow \mathrm{Alg}_{\mathrm{free}(A)}(\mathcal{C}) \times_{\mathcal{C}} \mathrm{Alg}_{\mathrm{free}(B)}(\mathcal{C}) 
\end{equation*}
with the functor
\begin{equation*}
\mathrm{alg}_{\mathrm{Sym}_{A\oplus B}}(\mathcal{C}) \rightarrow \mathrm{alg}_{\mathrm{Sym}_A}(\mathcal{C}) \times_{\mathcal{C}} \mathrm{alg}_{\mathrm{Sym}_B}(\mathcal{C}),
\end{equation*}
where now we are considering algebras for the \emph{functors} $\mathrm{Sym}_A$, $\mathrm{Sym}_B$, and $\mathrm{Sym}_{A \oplus B}$. Clearly this last map is indeed an equivalence. 

It remains to prove claim (b). The fact that 
\begin{equation*}
\varinjlim_n \varphi_n\mathcal{O} \rightarrow \mathcal{O}
\end{equation*}
is an equivalence is clear; for every fixed degree $k$, the directed system $\{(\varphi_n\mathcal{O})(k)\}_n$ becomes constant with value $\mathcal{O}(k)$ for $n \geq k$. The square of (b) arises by applying the left adjoint $\mathrm{Op}_{\leq n}(\mathcal{C}) \rightarrow \mathrm{Op}(\mathcal{C})$ of the restriction functor $\rho_n$ to the following square of $n$-truncated operads:
\[
\begin{tikzcd}
(\varphi_{n-1}\mathcal{O})(n) \ar{r}\ar{d} & \rho_n\varphi_{n-1} \mathcal{O} \ar{d} \\
\mathcal{O}(n) \ar{r} & \rho_n\varphi_n\mathcal{O}.
\end{tikzcd}
\]
It is straightforward to see that this is a pushout by observing that the right vertical map is an equivalence in every degree $k < n$, whereas in the top degree $n$ it agrees with the left vertical map.
\end{proof}

\begin{remark}
If we were working in a formalism where the $\infty$-category of $\mathcal{O}$-algebras in $\mathcal{C}$ would be encoded as some kind of $\infty$-category of morphisms $\mathcal{O} \rightarrow \mathcal{C}^{\otimes}$ (as is the case in the 1-categorical setting), then claim (a) of the proof above would be immediate.
\end{remark}

\section{The cotangent fiber and ind-conilpotent coalgebras}
\label{sec:conilpotentcoalgebras}

Let $\mathcal{O}$ be an operad that is reduced, i.e., satisfies $\mathcal{O}(1) = \mathbf{1}$. The evident (unique) morphisms of operads
\begin{equation*}
\mathbf{1} \xrightarrow{i} \mathcal{O} \xrightarrow{t} \mathbf{1}
\end{equation*}
induce adjunctions on $\infty$-categories of algebras as follows, where we make the identification $\mathrm{Alg}_{\mathbf{1}}(\mathcal{C}) \cong \mathcal{C}$:
\[
\begin{tikzcd}
\mathcal{C} \ar[shift left]{r}{i_!} & \mathrm{Alg}_{\mathcal{O}}(\mathcal{C}) \ar[shift left]{r}{t_!} \ar[shift left]{l}{i^*} & \mathcal{C}. \ar[shift left]{l}{t^*}
\end{tikzcd}
\] 
The functor $i^*$ is clearly the forgetful functor, so that its left adjoint $i_!$ may be identified with the free algebra functor $\mathrm{free}_{\mathcal{O}}$. Similarly, the functor $t^*$ equips an object $X$ of $\mathcal{C}$ with the trivial $\mathcal{O}$-algebra structure, so that its left adjoint $t_!$ is precisely the cotangent fiber $\mathrm{cot}_{\mathcal{O}}$ described in Section \ref{sec:mainresults}. Thus, we can restate the previous diagram with more descriptive notation as follows:
\[
\begin{tikzcd}
\mathcal{C} \ar[shift left]{r}{\mathrm{free}_{\mathcal{O}}} & \mathrm{Alg}_{\mathcal{O}}(\mathcal{C}) \ar[shift left]{r}{\mathrm{cot}_{\mathcal{O}}} \ar[shift left]{l}{\mathrm{forget}_{\mathcal{O}}} & \mathcal{C}. \ar[shift left]{l}{\mathrm{triv}_{\mathcal{O}}}
\end{tikzcd}
\] 
Since $t \circ i \cong \mathrm{id}$, the horizontal composites in both directions are the identity of $\mathcal{C}$. In particular, we find the formula
\begin{equation*}
\mathrm{cot}_{\mathcal{O}} \circ \mathrm{free}_{\mathcal{O}} \cong \mathrm{id}_{\mathcal{C}}.
\end{equation*}

\begin{remark}
It would also have been natural to call the functor $\mathrm{cot}_{\mathcal{O}}$ the $\mathcal{O}$-\emph{indecomposables}, but we reserve that term for another (very closely related) functor. The term cotangent fiber is inspired by the special case where $\mathcal{O}$ is the nonunital commutative operad. If $\mathcal{C} = \mathrm{Mod}_k$ is the derived $\infty$-category of a commutative ring $k$ and one makes the usual identification between nonunital commutative $k$-algebras and commutative $k$-algebras augmented over $k$, then for an augmented $k$-algebra $A$ there is the formula
\begin{equation*}
\mathrm{cot}_{\mathcal{O}}(A) \cong \mathbf{L}_{A/k} \otimes_A k.
\end{equation*}
Here $\mathbf{L}_{A/k}$ denotes the topological cotangent complex of $A$ over $k$. Thus, geometrically, $\mathrm{cot}_{\mathcal{O}}(A)$ is the fiber of the cotangent complex at the $k$-point of $\mathrm{Spec}(A)$ determined by the augmentation of $A$.
\end{remark}

Now any $\mathcal{O}$-algebra $X$ is canonically the geometric realization of a simplicial object consisting of free algebras:
\begin{equation*}
X \cong |(\mathrm{free}_{\mathcal{O}}\mathrm{forget}_{\mathcal{O}})^{\bullet+1} X|.
\end{equation*}
In the notation of bar constructions, one can write
\begin{equation*}
X \cong |B(\mathcal{O},\mathcal{O},X)_\bullet|,
\end{equation*}
where the action of $\mathcal{O}$ on $X$ should be interpreted as sending $X$ to $\mathrm{Sym}_{\mathcal{O}}(X)$. Applying $\mathrm{cot}_{\mathcal{O}}$ therefore yields the following:

\begin{proposition}
\label{prop:cotangentfiber}
The cotangent fiber of an $\mathcal{O}$-algebra $X$ may be computed by the formula
\begin{equation*}
\mathrm{cot}_{\mathcal{O}}(X) \cong |B(\mathbf{1},\mathcal{O},X)_\bullet|.
\end{equation*}
\end{proposition}

In fact, the cotangent fiber functor $\mathrm{cot}_{\mathcal{O}}$ exhibits $\mathcal{C}$ as the stabilization of the $\infty$-category $\mathrm{Alg}_{\mathcal{O}}(\mathcal{C})$. The following is a version of a theorem of Basterra--Mandell \cite{basterramandell} in the context of operads in spectra (see also \cite[Theorem 7.3.4.7]{higheralgebra}):

\begin{proposition}
\label{prop:stabAlgO}
The linearization of the adjoint pair $(\mathrm{cot}_{\mathcal{O}}, \mathrm{triv}_{\mathcal{O}})$ induces an equivalence of stable $\infty$-categories
\[
\begin{tikzcd}
\mathrm{Sp}(\mathrm{Alg}_{\mathcal{O}}(\mathcal{C})) \ar[shift left]{r}{\partial \mathrm{cot}_{\mathcal{O}}} & \mathcal{C}. \ar[shift left]{l}{\partial \mathrm{triv}_{\mathcal{O}}}
\end{tikzcd}
\]
\end{proposition}
\begin{proof}
By \cite[Corollary 6.2.2.16]{higheralgebra} and the fact that $\mathrm{Alg}_{\mathcal{O}}(\mathcal{C})$ is monadic over $\mathcal{C}$ (for the monad $\mathrm{Sym}_\mathcal{O}$), the stabilization $\mathrm{Sp}(\mathrm{Alg}_{\mathcal{O}}(\mathcal{C}))$ is monadic over $\mathcal{C}$ with respect to the linearization of the monad $\mathrm{Sym}_\mathcal{O}$, which is clearly the identity functor of $\mathcal{C}$ (cf. the proof of \cite[Proposition 7.3.4.10]{higheralgebra} for a version of this argument). The conclusion follows from the fact that
\begin{equation*}
\partial \mathrm{cot}_{\mathcal{O}} \circ \partial \mathrm{free}_{\mathcal{O}} \cong \partial(\mathrm{cot}_{\mathcal{O}} \circ\mathrm{free}_{\mathcal{O}}) \cong \mathrm{id}_{\mathcal{C}}. 
\end{equation*}
\end{proof}

\begin{remark}
\label{rmk:homology}
Proposition \ref{prop:stabAlgO} demonstrates that $\mathrm{cot}_{\mathcal{O}}$ can be regarded as the universal homology theory for $\mathcal{O}$-algebras. For specific choices of $\mathcal{O}$ it reproduces familiar notions of homology, such as (topological) Andr\'{e}--Quillen homology of commutative algebras, Hochschild homology of associative algebras, and Chevalley--Eilenberg homology of Lie algebras.
\end{remark}

Any stabilization as above comes with a corresponding comonad $\Sigma^\infty\Omega^\infty$. According to Proposition \ref{prop:stabAlgO}, this comonad can in this case be identified with $\mathrm{cot}_{\mathcal{O}} \circ \mathrm{triv}_{\mathcal{O}}$. This admits an explicit description as follows (cf. \cite[Lemma 3.3.4]{francisgaitsgory}).

\begin{proposition}
\label{prop:stabcomonad}
The comonad $\mathrm{cot}_\mathcal{O}\mathrm{triv}_{\mathcal{O}}$ is naturally equivalent to the comonad $\mathrm{Sym}_{B\mathcal{O}}$ associated with the cooperad $B\mathcal{O}$.
\end{proposition}
\begin{proof}
For an object $X$ of $\mathcal{C}$, Proposition \ref{prop:cotangentfiber} gives the calculation
\begin{equation*}
\mathrm{cot}_\mathcal{O}\mathrm{triv}_{\mathcal{O}}(X) \cong |B(\mathbf{1}, \mathcal{O}, \mathrm{triv}_{\mathcal{O}} X)_\bullet| \cong \mathrm{Sym}_{|B(\mathbf{1},\mathcal{O},\mathbf{1})_\bullet|}(X).
\end{equation*}
The second isomorphism uses the fact that the functor $A \mapsto \mathrm{Sym}_A$ preserves colimits and is monoidal. The identification above is in fact one of comonads, see e.g. \cite[Proposition 3.28]{BCN}.
\end{proof}

The functor $\mathrm{cot}_{\mathcal{O}}$ canonically factors through the $\infty$-category of coalgebras for the comonad $\mathrm{cot}_\mathcal{O}\mathrm{triv}_{\mathcal{O}}$. By Proposition \ref{prop:stabcomonad} we may identify that comonad with $\mathrm{Sym}_{B\mathcal{O}}$. We denote the resulting factorization of $\mathrm{cot}_{\mathcal{O}}$ by
\begin{equation*}
\mathrm{Alg}_{\mathcal{O}}(\mathcal{C}) \xrightarrow{\mathrm{indec}_{\mathcal{O}}^{\mathrm{nil}}} \mathrm{coAlg}_{B\mathcal{O}}^{\mathrm{dp},\mathrm{nil}}(\mathcal{C}) \xrightarrow{\mathrm{forget}^{\mathrm{nil}}_{B\mathcal{O}}} \mathcal{C}.
\end{equation*}
Here, as in Section \ref{sec:mainresults}, the middle term denotes the $\infty$-category of $\mathrm{Sym}_{B\mathcal{O}}$-coalgebras, which we refer to as \emph{ind-conilpotent divided power $B\mathcal{O}$-coalgebras}. The functor $\mathrm{indec}_{\mathcal{O}}^{\mathrm{nil}}$ admits a right adjoint $\mathrm{prim}^{\mathrm{nil}}_{B\mathcal{O}}$, which may be described as follows:

\begin{lemma}
\label{lem:primnilresolution}
For $Y \in \mathrm{coAlg}_{B\mathcal{O}}^{\mathrm{dp},\mathrm{nil}}(\mathcal{C})$, the $\mathcal{O}$-algebra $\mathrm{prim}_{B\mathcal{O}}^{\mathrm{nil}}(Y)$ is equivalent to the totalization
\begin{equation*}
\mathrm{prim}_{B\mathcal{O}}^{\mathrm{nil}}(Y) \cong \mathrm{Tot}\bigl(\mathrm{triv}_{\mathcal{O}}(\mathrm{Sym}_{B\mathcal{O}})^\bullet \mathrm{forget}_{B\mathcal{O}}^{\mathrm{nil}} Y\bigr).
\end{equation*}
\end{lemma}
\begin{proof}
Write $\mathrm{cofree}_{B\mathcal{O}}^{\mathrm{nil}}$ for the right adjoint to the functor $\mathrm{forget}^{\mathrm{nil}}_{B\mathcal{O}}$. Then $Y$ admits a canonical cosimplicial resolution by cofree $\mathrm{Sym}_{B\mathcal{O}}$-coalgebras:
\begin{equation*}
Y \cong \mathrm{Tot}\bigl(\mathrm{cofree}_{B\mathcal{O}}^{\mathrm{nil}}(\mathrm{Sym}_{B\mathcal{O}})^\bullet \mathrm{forget}_{B\mathcal{O}}^{\mathrm{nil}} Y\bigr).
\end{equation*}
The formula of the lemma is immediate from the identification $\mathrm{prim}^{\mathrm{nil}}_{B\mathcal{O}} \circ \mathrm{cofree}_{B\mathcal{O}}^{\mathrm{nil}} \cong \mathrm{triv}_{\mathcal{O}}$, which follows from the corresponding one for left adjoints:
\begin{equation*}
\mathrm{forget}_{B\mathcal{O}}^{\mathrm{nil}} \circ \mathrm{indec}_{B\mathcal{O}}^{\mathrm{nil}} \cong \mathrm{cot}_{\mathcal{O}}.
\end{equation*}
\end{proof}

\section{Divided power coalgebras}
\label{sec:coalgebras}

The aim of this section is to discuss the $\infty$-category $\mathrm{coAlg}_{\mathcal{Q}}^{\mathrm{dp}}(\mathcal{C})$ of divided power coalgebras for a cooperad $\mathcal{Q}$, as opposed to the $\infty$-category $\mathrm{coAlg}_{\mathcal{Q}}^{\mathrm{dp,nil}}(\mathcal{C})$ of \emph{ind-conilpotent} such coalgebras introduced above. Clearly a reasonable definition of the $\infty$-category of divided power $\mathcal{Q}$-coalgebras should have a description dual to the description of $\mathcal{O}$-algebras given in Remark \ref{rmk:decomposition}. In other words, the definition of the $\infty$-category $\mathrm{coAlg}_{\mathcal{Q}}^{\mathrm{dp}}(\mathcal{C})$ should satisfy a dual version of Theorem \ref{thm:algdecomposition}.

We provided one such definition in \cite[Definition 4.14]{heutsgoodwillie} (for coalgebras) and \cite[Section 6.1]{heutsgoodwillie} (for divided power coalgebras). Below we will supply a different (but equivalent) approach that is more directly adapted to our current setup. The main result of this section is Theorem \ref{thm:coalgdpdecomposition}, giving an inductive decomposition of the $\infty$-categories of divided power coalgebras analogous to Theorem \ref{thm:algdecomposition}. Our construction of the $\infty$-category $\mathrm{coAlg}_{\mathcal{Q}}^{\mathrm{dp}}(\mathcal{C})$ includes some technicalities on Pro-objects, which we largely defer to Appendix \ref{app:coalgebras}.

A divided power coalgebra structure on an object $X \in \mathcal{C}$ should in particular induce a structure map
\[
\delta\colon X \to \prod_{n \geq 1} (\mathcal{Q}(n) \otimes X^{\otimes n})_{h\Sigma_n} =: \widehat{\mathrm{Sym}}_{\mathcal{Q}}(X)
\]
describing the collection of `comultiplication maps' $\delta_n\colon X \to (\mathcal{Q}(n) \otimes X^{\otimes n})_{h\Sigma_n}$. In addition, such a coalgebra $X$ should come equipped with a coherent system of homotopies relating the various compositions of the $\delta_n$. The issue in describing this system precisely is that the functor $\widehat{\mathrm{Sym}}_{\mathcal{Q}}$ is generally not a comonad on $\mathcal{C}$. Indeed, the assignment
\[
\mathrm{SSeq}(\mathcal{C}) \to \mathrm{Fun}(\mathcal{C},\mathcal{C})\colon A \mapsto \widehat{\mathrm{Sym}}_{A}
\]
is only lax monoidal; it will therefore send operads to monads, but not cooperads to comonads.

Our construction of $\mathrm{coAlg}_{\mathcal{Q}}^{\mathrm{dp}}(\mathcal{C})$ proceeds through the following steps, which we detail in Appendix \ref{app:coalgebras}:
\begin{itemize}
\item[(1)] For each $n \geq 1$, the functor $D_n^{\mathcal{Q}}\colon \mathcal{C} \to \mathcal{C}$ admits a natural extension to a functor $pD_n^{\mathcal{Q}}\colon \mathrm{Pro}(\mathcal{C}) \to \mathrm{Pro}(\mathcal{C})$ preserving filtered limits. We define an endofunctor of $\mathrm{Pro}(\mathcal{C})$ by the formula
\[
\widehat{p\mathrm{Sym}}_{\mathcal{Q}} := \prod_{n \geq 1} pD_n^{\mathcal{Q}}.
\]
\item[(2)] The functor $\widehat{p\mathrm{Sym}}_{\mathcal{Q}}$ is naturally a comonad on $\mathrm{Pro}(\mathcal{C})$ (cf. Lemma \ref{lem:pSymmon}).
\item[(3)] The $\infty$-category $\mathrm{coAlg}_{\mathcal{Q}}^{\mathrm{dp}}(\mathcal{C})$ is the full subcategory of the $\infty$-category of $\widehat{p\mathrm{Sym}}_{\mathcal{Q}}$-coalgebras spanned by coalgebras whose underlying object of $\mathrm{Pro}(\mathcal{C})$ is pro-constant (Definition \ref{def:dpQcoalgebras}).
\end{itemize}

The following is the main result of Appendix \ref{app:coalgebras}:

\begin{theorem}
\label{thm:coalgdpdecomposition}
For each $n \geq 2$, the evident commutative square of $\infty$-categories
\[
\begin{tikzcd}
\mathrm{coAlg}_{\varphi^n\mathcal{Q}}^{\mathrm{dp}}(\mathcal{C}) \ar{r}\ar{d} & \mathrm{coalg}_{D_n^{\mathcal{Q}}}(\mathcal{C}) \ar{d} \\
\mathrm{coAlg}_{\varphi^{n-1}\mathcal{Q}}^{\mathrm{dp}}(\mathcal{C}) \ar{r} & \mathrm{coalg}_{D_n^{\varphi^{n-1}\mathcal{Q}}}(\mathcal{C})
\end{tikzcd}
\]
is a pullback. Furthermore, the natural map
\begin{equation*}
\mathrm{coAlg}_{\mathcal{Q}}^{\mathrm{dp}}(\mathcal{C}) \rightarrow \varprojlim_n \mathrm{coAlg}_{\varphi^n\mathcal{Q}}^{\mathrm{dp}}(\mathcal{C})
\end{equation*}
is an equivalence of $\infty$-categories.
\end{theorem}

\begin{remark}
In fact, one could almost take the statement of Theorem \ref{thm:coalgdpdecomposition} as an inductive definition of the $\infty$-category of divided power $\mathcal{Q}$-coalgebras. All that would be needed to make this work is a construction of the bottom horizontal functor
\begin{equation*}
\mathrm{coAlg}_{\varphi^{n-1}\mathcal{Q}}^{\mathrm{dp}}(\mathcal{C}) \rightarrow \mathrm{coalg}_{D_n^{\varphi^{n-1}\mathcal{Q}}}(\mathcal{C}). 
\end{equation*}
Since Theorem \ref{thm:coalgdpdecomposition} is all that we will need of the $\infty$-category of divided power $\mathcal{Q}$-coalgebras, any definition that satisfies its conclusion will work for our purposes.
\end{remark}


In the statement of the theorem we have used the notation $\mathrm{coalg}_F(\mathcal{C})$ for the $\infty$-category of coalgebras for a functor $F\colon \mathcal{C} \to \mathcal{C}$. Its definition is dual to that of the $\infty$-category of $F$-algebras; to be precise, it is defined by the pullback square 
\[
\begin{tikzcd}
\mathrm{coalg}_F(\mathcal{C}) \ar{r}\ar{d} & \mathcal{C}^{\Delta^1} \ar{d}{(\mathrm{ev}_0,\mathrm{ev}_1)} \\
\mathcal{C} \ar{r}{(\mathrm{id},F)} & \mathcal{C} \times \mathcal{C}.
\end{tikzcd}
\]
If $F$ is accessible, then $\mathrm{coalg}_F(\mathcal{C})$ is presentable. Indeed, we know that all but the top left corner are presentable $\infty$-categories and both the right vertical and lower horizontal functor are accessible. Therefore $\mathrm{coalg}_F(\mathcal{C})$ is accessible by 
\cite[Proposition 5.4.6.6]{HTT}. Since it has all colimits (which are created by the left vertical functor to $\mathcal{C}$), it is also presentable.

\begin{proposition}
\label{prop:coalgebrascomonadic}
The forgetful functor
\begin{equation*}
\mathrm{coAlg}^{\mathrm{dp}}_{\mathcal{Q}}(\mathcal{C}) \rightarrow \mathcal{C}
\end{equation*}
is the left adjoint of a comonadic adjunction. Moreover, the $\infty$-category $\mathrm{coAlg}^{\mathrm{dp}}_{\mathcal{Q}}(\mathcal{C})$ is presentable.
\end{proposition}
\begin{proof}
First we show that if $F\colon \mathcal{C} \rightarrow \mathcal{C}$ is an accessible functor, then the forgetful functor
\begin{equation*}
U_F\colon \mathrm{coalg}_F(\mathcal{C}) \rightarrow \mathcal{C},
\end{equation*}
exhibits the $\infty$-category of $F$-coalgebras as comonadic over $\mathcal{C}$. (Note that $U_F$ admits a right adjoint since it preserves colimits and $\mathrm{coalg}_F(\mathcal{C})$ is presentable.) This follows rather straightforwardly by checking the conditions of the Barr--Beck--Lurie theorem \cite[Theorem 4.7.3.5]{higheralgebra}):
\begin{itemize}
\item[(1)] The forgetful functor $U_F$ is clearly conservative, i.e., it detects isomorphisms.
\item[(2)] Say 
\begin{equation*}
\mathcal{X}^\bullet\colon \mathbf{\Delta} \rightarrow \mathrm{coalg}_F(\mathcal{C})\colon [n] \mapsto (X^n \rightarrow F(X^n))
\end{equation*} 
is a cosimplicial object in $\mathrm{coalg}_F(\mathcal{C})$ that is $U_F$-split, i.e., admits an extra codegeneracy after applying the forgetful functor $U_F$. Then the totalization of $U_F\mathcal{X}^\bullet = X^\bullet$ is preserved by any functor. In particular, the morphism $F(\mathrm{Tot} X^\bullet) \rightarrow \mathrm{Tot} F(X^\bullet)$ is an isomorphism and using its inverse we find an $F$-coalgebra 
\begin{equation*}
\mathrm{Tot} X^\bullet \rightarrow \mathrm{Tot} F(X^\bullet) \simeq F(\mathrm{Tot} X^\bullet).
\end{equation*}
This is the limit of $\mathcal{X}^\bullet$ in $\mathrm{coalg}_F(\mathcal{C})$ and it is preserved by $U_F$.
\end{itemize}
To bootstrap this result to the $\infty$-category $\mathrm{coAlg}^{\mathrm{dp}}_{\mathcal{Q}}(\mathcal{C})$ we apply Theorem \ref{thm:coalgdpdecomposition} and use that (accessible) comonadic adjunctions over $\mathcal{C}$ are closed under limits. To be more precise, the functor 
\begin{equation*}
\mathrm{coMnd}^{\mathrm{acc}}(\mathcal{C}) \rightarrow \mathcal{P}\mathrm{r}^{L}_{/\mathcal{C}}\colon T \mapsto (\mathrm{coAlg}_T(\mathcal{C}) \xrightarrow{\mathrm{forget}} \mathcal{C})
\end{equation*}
admits a left adjoint, assigning to a left adjoint functor $\mathcal{D} \rightarrow \mathcal{C}$ the resulting comonad on $\mathcal{C}$ (cf. \cite[Theorem 1.3]{heine}). Moreover, the counit of this adjunction is an equivalence. This exhibits the $\infty$-category  $\mathrm{coMnd}^{\mathrm{acc}}(\mathcal{C})$ as a localization of $\mathcal{P}\mathrm{r}^L_{/\mathcal{C}}$. In particular, it makes $\mathrm{coMnd}^{\mathrm{acc}}(\mathcal{C})$ equivalent to a full subcategory of $\mathcal{P}\mathrm{r}^L_{/\mathcal{C}}$ that is closed under limits. The fact that limits in $\mathcal{P}\mathrm{r}^L_{/\mathcal{C}}$ agree with limits in $(\mathrm{Cat}_\infty)_{/\mathcal{C}}$ (cf. \cite[Proposition 5.5.3.13]{HTT}) then concludes the argument.
\end{proof}

As a consequence of Proposition \ref{prop:coalgebrascomonadic} we obtain cofree divided power coalgebra functors. Contrary to the case of algebras, there is generally no simple formula describing these, sometimes making arguments with coalgebras more subtle than their counterparts for algebras.

\begin{corollary}
The forgetful functor
\begin{equation*}
\mathrm{coAlg}^{\mathrm{dp}}_{\mathcal{Q}}(\mathcal{C}) \rightarrow \mathcal{C}
\end{equation*}
admits a right adjoint
\begin{equation*}
\mathcal{C} \xrightarrow{\mathrm{cofree}^{\mathrm{dp}}_{\mathcal{Q}}} \mathrm{coAlg}^{\mathrm{dp}}_{\mathcal{Q}}(\mathcal{C}).
\end{equation*}
\end{corollary}

The relation between the comonads $\mathrm{Sym}_{\mathcal{Q}}$ and $\widehat{p\mathrm{Sym}}_{\mathcal{Q}}$ produces a natural left adjoint comparison functor from the $\infty$-category $\mathrm{coAlg}^{\mathrm{dp},\mathrm{nil}}_{\mathcal{Q}}(\mathcal{C})$ of conilpotent divided power coalgebras to that of general divided power coalgebras. We give this construction in Section \ref{sec:appdpcoalgs}, right after Definition \ref{def:dpQcoalgebras}. There is also a way to construct the comparison functor inductively by means of Theorem \ref{thm:coalgdpdecomposition}. Indeed, an object $X$ of $\mathrm{coAlg}^{\mathrm{dp},\mathrm{nil}}_{\mathcal{Q}}(\mathcal{C})$ is precisely a coalgebra for the comonad $\mathrm{Sym}_{\mathcal{Q}}$. The map of comonads $\mathrm{Sym}_{\mathcal{Q}} \rightarrow \mathrm{Sym}_{\varphi^n\mathcal{Q}}$ then provides a commutative diagram as follows:
\[
\begin{tikzcd}
& \mathrm{Sym}_{\mathcal{Q}}(X) \ar{d}\ar{r} & D_n^{\mathcal{Q}}(X) \ar{d} \\
X \ar{ur}\ar{r} & \mathrm{Sym}_{\varphi^n\mathcal{Q}}(X) \ar{r} & D_n^{\varphi^n\mathcal{Q}}(X).
\end{tikzcd}
\]
Now suppose we have inductively constructed a functor $\mathrm{coAlg}^{\mathrm{dp},\mathrm{nil}}_{\mathcal{Q}}(\mathcal{C}) \rightarrow \mathrm{coAlg}^{\mathrm{dp}}_{\varphi^n\mathcal{Q}}(\mathcal{C})$, starting from the case $n=1$ given by 
\begin{equation*}
\mathrm{coAlg}^{\mathrm{dp},\mathrm{nil}}_{\mathcal{Q}}(\mathcal{C}) \xrightarrow{\mathrm{forget}^{\mathrm{nil}}_{\mathcal{Q}}} \mathcal{C}.
\end{equation*}
Then the diagram above, combined with Theorem \ref{thm:coalgdpdecomposition}, provides a lift to a functor
\begin{equation*}
\mathrm{coAlg}^{\mathrm{dp},\mathrm{nil}}_{\mathcal{Q}}(\mathcal{C}) \rightarrow \mathrm{coAlg}^{\mathrm{dp}}_{\varphi^{n+1}\mathcal{Q}}(\mathcal{C}).
\end{equation*}
Taking the limit over $n$ gives a left adjoint functor
\begin{equation*}
\mathrm{coAlg}^{\mathrm{dp},\mathrm{nil}}_{\mathcal{Q}}(\mathcal{C}) \rightarrow \mathrm{coAlg}^{\mathrm{dp}}_{\mathcal{Q}}(\mathcal{C})
\end{equation*}
as desired.

\begin{definition}
We define the indecomposables functor $\mathrm{indec}_{\mathcal{O}}$ to be the composition of the functors
\[
\mathrm{Alg}_{\mathcal{O}}(\mathcal{C}) \xrightarrow{\mathrm{indec}_{\mathcal{O}}^{\mathrm{nil}}} \mathrm{coAlg}^{\mathrm{dp},\mathrm{nil}}_{B\mathcal{O}}(\mathcal{C}) \to \mathrm{coAlg}^{\mathrm{dp}}_{B\mathcal{O}}(\mathcal{C}),
\]
with the second arrow being the comparison functor constructed above.
\end{definition}

Since $\mathrm{indec}_{\mathcal{O}}$ preserves colimits, the adjoint functor theorem guarantees the existence of a right adjoint, which we will denote by $\mathrm{prim}_{B\mathcal{O}}$. Most of our constructions so far can be conveniently summarized in the following commutative diagram of left adjoints
\[
\begin{tikzcd}
& \mathrm{Alg}_{\mathcal{O}}(\mathcal{C}) \ar{dd}{\mathrm{indec}_{\mathcal{O}}} \ar{dr}{\mathrm{cot}_{\mathcal{O}}} & \\
\mathcal{C} \ar{ur}{\mathrm{free}_{\mathcal{O}}}\ar{dr}[swap]{\mathrm{triv}_{B\mathcal{O}}} && \mathcal{C} \\
& \mathrm{coAlg}^{\mathrm{dp}}_{B\mathcal{O}}(\mathcal{C}) \ar{ur}[swap]{\mathrm{forget}_{B\mathcal{O}}} &
\end{tikzcd}
\]
or equivalently that of right adjoints
\[\begin{tikzcd}
& \mathrm{Alg}_{\mathcal{O}}(\mathcal{C}) \ar{dl}[swap]{\mathrm{forget}_{\mathcal{O}}} & \\
\mathcal{C} && \mathcal{C} \ar{ul}[swap]{\mathrm{triv}_{\mathcal{O}}}\ar{dl}{\mathrm{cofree}_{B\mathcal{O}}^{\mathrm{dp}}} \\
& \mathrm{coAlg}^{\mathrm{dp}}_{B\mathcal{O}}(\mathcal{C}). \ar{uu}[swap]{\mathrm{prim}_{B\mathcal{O}}} \ar{ul} &
\end{tikzcd}
\]
The composite functors $\mathcal{C} \rightarrow \mathcal{C}$ formed in these diagrams are all naturally isomorphic to the identity. The only functor we have not explicitly discussed yet is $\mathrm{triv}_{B\mathcal{O}}$. It can either be defined by first constructing trivial ind-conilpotent coalgebras in $\mathrm{coAlg}^{\mathrm{dp},\mathrm{nil}}_{B\mathcal{O}}(\mathcal{C})$ using the morphism of comonads
\begin{equation*}
\mathrm{id} = \mathrm{Sym}_{\mathbf{1}} \rightarrow \mathrm{Sym}_{B\mathcal{O}}
\end{equation*}
and applying the comparison functor $\mathrm{coAlg}^{\mathrm{dp},\mathrm{nil}}_{B\mathcal{O}}(\mathcal{C}) \rightarrow \mathrm{coAlg}^{\mathrm{dp}}_{B\mathcal{O}}(\mathcal{C})$, or by applying Theorem \ref{thm:coalgdpdecomposition} and using the zero maps $X \rightarrow D_n^{B\mathcal{O}}(X)$ for all $n \geq 2$ as coalgebra structure maps.

To conclude this section, we record an important observation on cofree divided power coalgebras. Although general formulas for such cofree coalgebras are hard to come by, we record one special instance where they \emph{can} be described explicitly, namely the case of a truncated cooperad. The following will be very useful to us later. We defer the proof to Appendix \ref{app:coalgebras}.

\begin{proposition}
\label{prop:cofreetruncated}
For a truncated cooperad $\mathcal{Q} = \tau^n\mathcal{Q}$, the comparison functor
\begin{equation*}
\mathrm{coAlg}^{\mathrm{dp},\mathrm{nil}}_{\mathcal{Q}}(\mathcal{C}) \rightarrow \mathrm{coAlg}^{\mathrm{dp}}_{\mathcal{Q}}(\mathcal{C})
\end{equation*}
is an equivalence of $\infty$-categories. In particular, for $X \in \mathcal{C}$ we have a natural isomorphism
\begin{equation*}
\mathrm{cofree}^{\mathrm{dp}}_{\mathcal{Q}}(X) \cong \mathrm{Sym}_{\mathcal{Q}}(X) = \bigoplus_{k=1}^n D_k^{\mathcal{Q}}(X).
\end{equation*}
\end{proposition}

\section{Graded and filtered coalgebras}
\label{sec:filteredcoalgebras}

In Section \ref{sec:truncationsalg} we will discuss Theorem \ref{thm:filtrations}, which provides several ways of decomposing algebras and coalgebras into simpler (meaning free or trivial) pieces. That theorem is the key ingredient in the proof of our main result, Theorem \ref{thm:completeness}, on Koszul duality between $\mathcal{O}$-algebras and $B\mathcal{O}$-coalgebras. To prepare for the proof of Theorem \ref{thm:filtrations} we develop some material on graded and filtered coalgebras in this section. Much of our discussion is dual to the discussion of filtered algebras in \cite{brantnermathew}, but it should be noted that our arguments here are a bit different due to the fact that we cannot assume that the tensor product in $\mathcal{C}$ commutes with cosifted limits. (When discussing filtered algebras, one often exploits the fact that the tensor product preserves sifted colimits in order to reduce proofs to the case of free algebras.)

\subsection{Graded coalgebras}

We write $\mathbb{Z}^{\delta}_{\geq 1}$ for the set of strictly positive natural numbers, regarded as a discrete category. We denote by
\begin{equation*}
\mathrm{gr}(\mathcal{C}) := \mathrm{Fun}(\mathbb{Z}^{\delta}_{\geq 1}) = \prod_{i=1}^\infty \mathcal{C}
\end{equation*}
the $\infty$-category of graded objects in $\mathcal{C}$. Often we write objects of $\mathrm{gr}(\mathcal{C})$ as $X = \{X(n)\}_{n \geq 1}$, with $X(n)$ the $n$th graded piece of $X$. 

The category $\mathbb{Z}^{\delta}_{\geq 1}$ has an evident nonunital symmetric monoidal structure given by addition. Applying Day convolution \cite[Section 2.2.6]{higheralgebra} this also equips $\mathrm{gr}(\mathcal{C})$ with a nonunital symmetric monoidal structure. For graded objects $X$ and $Y$, their tensor product is loosely described by the formula
\begin{equation*}
(X \otimes Y)(i) = \bigoplus_{a+b = i} X(a) \otimes Y(b).
\end{equation*}
Now let $\mathcal{Q}$ be a cooperad in $\mathcal{C}$. By regarding $\mathrm{gr}(\mathcal{C})$ as tensored over $\mathcal{C}$, we can speak of divided power $\mathcal{Q}$-coalgebras in $\mathrm{gr}(\mathcal{C})$. We will refer to these as \emph{graded divided power $\mathcal{Q}$-coalgebras} and write $\mathrm{coAlg}_{\mathcal{Q}}^{\mathrm{dp}}(\mathrm{gr}(\mathcal{C}))$ for the $\infty$-category of such. (Alternatively, one can regard $\mathrm{gr}(\mathcal{C})$ as a full subcategory of the symmetric monoidal $\infty$-category $\prod_{i=0}^\infty \mathcal{C}$ of nonnegatively graded objects. Then one can interpret $\mathcal{Q}$ as a cooperad concentrated in degree 0 in this $\infty$-category and take coalgebras for it that are concentrated in strictly positive degrees.)

Generally cofree coalgebras are difficult to describe, but in the graded case the `naive' formula is correct:

\begin{lemma}
\label{lem:cofreegradedcoalg}
Let $X \in \mathrm{gr}(\mathcal{C})$ be a graded object and $\mathcal{Q}$ a cooperad in $\mathcal{C}$. Then there is a natural isomorphism
\begin{equation*}
\mathrm{cofree}_{\mathcal{Q}}^{\mathrm{dp}}(X) \cong \mathrm{Sym}_{\mathcal{Q}}(X) = \bigoplus_{n=1}^\infty D^{\mathcal{Q}}_n(X).
\end{equation*}
\end{lemma}
\begin{proof}
Let $k \geq 1$ and write $\mathrm{gr}_{\leq k}(\mathcal{C}) = \prod_{n=1}^k \mathcal{C}$ for the $\infty$-category of `$k$-truncated' graded objects. Write 
\begin{equation*}
j\colon \mathrm{gr}_{\leq k}(\mathcal{C}) \rightarrow \mathrm{gr}(\mathcal{C}) 
\end{equation*}
for the `extension by zero' functor, sending a truncated object $\{X(n)\}_{1 \leq n \leq k}$ to the graded object
\begin{equation*}
j(X)(n) = \begin{cases} X(n) & \text{if } n \leq k \\ 0 & \text{if } n > k. \end{cases}
\end{equation*}
Then $j$ is fully faithful and it has a right adjoint $r$, simply restricting a graded object to degrees $\leq r$. (In fact $r$ is also left adjoint to $j$.) The $\infty$-category $\mathrm{gr}_{\leq k}(\mathcal{C})$ has an essentially unique nonunital symmetric monoidal structure making $r$ a nonunital symmetric monoidal functor; the tensor product of truncated graded objects is given by 
\begin{equation*}
X \otimes Y = r(j(X) \otimes j(Y)).
\end{equation*}
Since $j$ is left adjoint to $r$, it acquires an oplax nonunital symmetric monoidal structure. Because of this, and the fact that $j$ preserves colimits, it lifts to a functor defined on divided power coalgebras (see Appendix \ref{app:functoriality}) to give a commutative square
\[
\begin{tikzcd}
\mathrm{coAlg}_{\mathcal{Q}}^{\mathrm{dp}}(\mathrm{gr}_{\leq k}(\mathcal{C})) \ar{r}{j}\ar{d}{\mathrm{forget}} & \mathrm{coAlg}_{\mathcal{Q}}^{\mathrm{dp}}(\mathrm{gr}(\mathcal{C})) \ar{d}{\mathrm{forget}} \\
\mathrm{gr}_{\leq k}(\mathcal{C}) \ar{r}{j} & \mathrm{gr}(\mathcal{C}).
\end{tikzcd}
\]
Passing to right adjoints yields a commutative square
\[
\begin{tikzcd}
\mathrm{coAlg}_{\mathcal{Q}}^{\mathrm{dp}}(\mathrm{gr}_{\leq k}(\mathcal{C})) & \mathrm{coAlg}_{\mathcal{Q}}^{\mathrm{dp}}(\mathrm{gr}(\mathcal{C})) \ar[shift left]{l}[swap]{r} \\
\mathrm{gr}_{\leq k}(\mathcal{C}) \ar{u}[swap]{\mathrm{cofree}^{\mathrm{dp}}_{\mathcal{Q}}} & \mathrm{gr}(\mathcal{C}) \ar{u}[swap]{\mathrm{cofree}^{\mathrm{dp}}_{\mathcal{Q}}} \ar[shift left]{l}[swap]{r}.
\end{tikzcd}
\]
The tensor product on $\mathrm{gr}_{\leq k}(\mathcal{C})$ is nilpotent; more precisely, for any $m > k,$ the $m$-fold tensor product functor is identically zero. It is straightforward to deduce from this that
\begin{equation*}
\mathrm{coAlg}_{\mathcal{Q}}^{\mathrm{dp}}(\mathrm{gr}_{\leq k}(\mathcal{C})) \cong \mathrm{coAlg}_{\tau^k\mathcal{Q}}^{\mathrm{dp}}(\mathrm{gr}_{\leq k}(\mathcal{C}))
\end{equation*}
e.g., by applying Theorem \ref{thm:coalgdpdecomposition} for the cooperads $\mathcal{Q}$ and $\tau^k\mathcal{Q}$. Applying Proposition \ref{prop:cofreetruncated} to identify $\mathrm{cofree}_{\tau^k\mathcal{Q}}^{\mathrm{dp}}$ with $\mathrm{Sym}_{\tau^k\mathcal{Q}}$, we conclude that
\begin{equation*}
r \circ \mathrm{cofree}^{\mathrm{dp}}_{\mathcal{Q}} \cong \mathrm{Sym}_{\tau^k\mathcal{Q}} \circ r.
\end{equation*}
In different words, $\mathrm{cofree}^{\mathrm{dp}}_{\mathcal{Q}}$ agrees with $\mathrm{Sym}_{\tau^k\mathcal{Q}}$ in degrees $\leq k$. Since $k$ was arbitrary, this implies the lemma.
\end{proof}

\subsection{Filtered coalgebras}

Write $\mathbb{Z}_{\geq 1}$ for the linearly ordered set of positive integers, regarded as a category in which there is a unique morphism $i \rightarrow j$ whenever $i \leq j$. We write $\mathrm{fil}(\mathcal{C})$ for the $\infty$-category $\mathrm{Fun}(\mathbb{Z}_{\geq 1},\mathcal{C})$ and think of an object of this category as representing an object $X^1$ of $\mathcal{C}$ equipped with a sequence of `quotients' $X^n$ of it. We denote a general object of $\mathrm{fil}(\mathcal{C})$ by $X^\bullet$ or 
\begin{equation*}
X^1 \rightarrow X^2 \rightarrow \cdots.
\end{equation*}

The $\infty$-category $\mathrm{fil}(\mathcal{C})$ can be equipped with the (nonunital) symmetric monoidal structure of \emph{right} Day convolution. To be precise, consider $\mathbb{Z}_{\geq 1}$ as a nonunital symmetric monoidal category with tensor product given by $i \otimes j = i + j$. Then (cf. \cite[Section 2.2.6]{higheralgebra}) the opposite $\infty$-category
\begin{equation*}
\mathrm{fil}(\mathcal{C})^{\mathrm{op}} = \mathrm{Fun}(\mathbb{Z}_{\geq 1}^{\mathrm{op}},\mathcal{C}^{\mathrm{op}})
\end{equation*}
admits the usual (nonunital) symmetric monoidal structure given by (left) Day convolution, in which the tensor product of objects $X^\bullet$ and $Y^\bullet$ is the left Kan extension of the two-variable functor $(a,b) \mapsto X^a \otimes Y^b$ along the tensor product on $\mathbb{Z}_{\geq 1}^{\mathrm{op}}$. This also gives a nonunital symmetric monoidal structure on the opposite $\infty$-category $\mathrm{fil}(\mathcal{C})$. Informally, it is described by the formula
\begin{equation*}
(X \otimes Y)^n = \varprojlim_{n\leq a+b} X^a \otimes Y^b.
\end{equation*}

\begin{remark}
Note that \cite[Example 2.2.6.17]{higheralgebra} does not quite apply as stated, as it would require the tensor product on $\mathcal{C}^{\mathrm{op}}$ to preserve $\kappa$-small colimits in each variable separately, for $\kappa$ an uncountable regular cardinal. However, inspecting the proofs of Corollary 2.2.6.14 and Proposition 2.2.6.16 of loc. cit. it is clear that for a symmetric monoidal $\infty$-category $\mathcal{D}$ the Day convolution on a functor category $\mathrm{Fun}(I,\mathcal{D})$ will exist under the following more precise hypothesis:
\begin{itemize}
\item[(*)] The $\infty$-category $\mathcal{D}$ admits colimits indexed over the $\infty$-categories $I^n \times_{I} I/x$ parametrizing $n$-tuples $(y_1, \ldots, y_n)$ equipped with a map $y_1 \otimes \cdots \otimes y_n \rightarrow x$, for each $x \in I$ and each $n$. Moreover, such colimits are preserved in each variable separately by the tensor product of $\mathcal{D}$.
\end{itemize}
In the particular case $I = \mathbb{Z}_{\geq 1}^{\mathrm{op}}$, the category $I^n \times_{I} I/x$ may be described as the opposite of the poset of $n$-tuples $(k_1, \ldots, k_n)$ of positive natural numbers such that $k_1 + \cdots + k_n \geq x$. This poset contains as a cofinal subset the one where we add the restriction $k_i \leq x$ for all $i$. Observe that the latter poset is finite. Therefore the relevant colimits may be computed over finite diagrams; the fact that $\mathcal{C}^{\mathrm{op}}$ is stable and its tensor product is exact in each variable then guarantees that condition (*) is satisfied.
\end{remark}

\begin{remark}
In our discussion of graded coalgebras we did not emphasize right Day convolution. The reason is that the indexing category $\mathbb{Z}^{\delta}_{\geq 1}$ there is discrete, so that left and right Day convolution agree.
\end{remark}

Evaluation at $1$ defines a functor we denote by 
\begin{equation*}
u\colon \mathrm{fil}(\mathcal{C}) \rightarrow \mathcal{C}.
\end{equation*}
We think of $u(X^\bullet) = X^1$ as the `underlying object' of $X^\bullet$. This functor preserves small limits and colimits. In fact it has a left adjoint sending $X \in \mathcal{C}$ to the constant filtration $\mathrm{con}(X) = (X = X = \cdots)$ and a right adjoint sending $X$ to $(X \rightarrow 0 = \cdots)$. Clearly both these adjoints are fully faithful.

Note that
\begin{equation*}
u(X^\bullet \otimes Y^\bullet) = \varprojlim_{1\leq a+b} X^a \otimes Y^b \cong X^1 \otimes Y^1,
\end{equation*}
using that $(1,1)$ is initial among pairs $(a,b) \in \mathbb{Z}_{\geq 1} \times \mathbb{Z}_{\geq 1}$ with $1 \leq a+b$. More precisely, $u$ can be given the structure of a nonunital symmetric monoidal functor. Hence it induces a colimit-preserving functor
\begin{equation*}
u\colon \mathrm{coAlg}_{\mathcal{Q}}^{\mathrm{dp}}(\mathrm{fil}(\mathcal{C})) \rightarrow \mathrm{coAlg}_{\mathcal{Q}}^{\mathrm{dp}}(\mathcal{C}).
\end{equation*}
We may now apply the adjoint functor theorem to guarantee the existence of the following:

\begin{definition}
The \emph{adic filtration} of divided power $\mathcal{Q}$-coalgebras is the right adjoint of $u$:
\begin{equation*}
a^\bullet\colon \mathrm{coAlg}_{\mathcal{Q}}^{\mathrm{dp}}(\mathcal{C}) \rightarrow \mathrm{coAlg}_{\mathcal{Q}}^{\mathrm{dp}}(\mathrm{fil}(\mathcal{C}))\colon Y \mapsto (a^1 Y \rightarrow a^2 Y \rightarrow \cdots).
\end{equation*}
\end{definition}

Let us record the following simple observation:

\begin{lemma}
\label{lem:adicfullyfaithful}
The functor $a^\bullet$ is fully faithful or, equivalently, the counit map $u(a^\bullet Y) \rightarrow Y$ is an isomorphism for every $Y \in \mathrm{coAlg}_{\mathcal{Q}}^{\mathrm{dp}}(\mathcal{C})$. 
\end{lemma}
\begin{proof}
Since $u$ is nonunital symmetric monoidal, its left adjoint $\mathrm{con}\colon \mathcal{C} \rightarrow \mathrm{fil}(\mathcal{C})$ acquires a nonunital oplax symmetric monoidal structure. Since it also preserves colimits, $\mathrm{con}$ lifts to a functor
\begin{equation*}
\mathrm{con}\colon \mathrm{coAlg}_{\mathcal{Q}}^{\mathrm{dp}}(\mathcal{C}) \rightarrow \mathrm{coAlg}_{\mathcal{Q}}^{\mathrm{dp}}(\mathrm{fil}(\mathcal{C}))
\end{equation*}
defined on divided power coalgebras, which is still fully faithful. Now for $Y, Z \in \mathrm{coAlg}_{\mathcal{Q}}^{\mathrm{dp}}(\mathcal{C})$ we find natural isomorphisms
\begin{eqnarray*}
\mathrm{Map}(Z, u(a^\bullet Y)) & \cong & \mathrm{Map}(\mathrm{con}(Z), a^\bullet Y) \\
& \cong & \mathrm{Map}(u(\mathrm{con}(Z)), Y) \\
& \cong & \mathrm{Map}(Z, Y),
\end{eqnarray*}
proving the lemma.
\end{proof}

Although the cofree filtered coalgebra is as complicated to describe as the `usual' cofree coalgebra, we have a simple formula in the case of a truncated cooperad:

\begin{lemma}
\label{lem:cofreetruncatedfiltcoalg}
Let $n \geq 1$ and let $Y = \mathrm{cofree}^{\mathrm{dp}}_{\tau^n\mathcal{Q}}(V)$ for an object $V \in \mathcal{C}$. Then there are natural isomorphisms
\begin{equation*}
a^k Y \cong \bigoplus_{k \leq j \leq n} D_n^{\mathcal{Q}}(V).
\end{equation*}
\end{lemma}
\begin{proof}
The composition $a^{\bullet} \circ \mathrm{cofree}^{\mathrm{dp}}_{\tau^n\mathcal{Q}}$ is right adjoint to the functor assigning to $Z \in \mathrm{coAlg}^{\mathrm{dp}}_{\tau^n\mathcal{Q}}(\mathrm{fil}(\mathcal{C}))$ the object $Z^1 \in \mathcal{C}$. Hence $a^{\bullet} \circ \mathrm{cofree}^{\mathrm{dp}}_{\tau^n\mathcal{Q}}$ assigns to $V \in \mathcal{C}$ the cofree filtered divided power $\tau^n\mathcal{Q}$-coalgebra on the filtered object $[V]_1 = (V \rightarrow 0 \rightarrow \cdots)$. Proposition \ref{prop:cofreetruncated} then gives the formula we are after, observing that the $j$-fold tensor product of $[V]_1$ is the filtered object that equals $V^{\otimes j}$ in degrees $\leq j$ and vanishes in higher degrees.
\end{proof}

Generally there is no reason for the adic filtration $a^{\bullet} Y$ of a coalgebra $Y$ to be \emph{cocomplete}, in the sense that $\varinjlim_n a^n Y = 0$. (We will see later that this happens precisely when $Y$ is conilcomplete, cf. Remark \ref{rmk:adictn}.) However, this will always happen in the special case where $\mathcal{Q}$ is a truncated cooperad:  

\begin{lemma}
\label{lem:finiteadicfiltration}
For a truncated cooperad $\tau^n\mathcal{Q}$, the adic filtration
\begin{equation*}
a^{\bullet}\colon \mathrm{coAlg}_{\tau^n\mathcal{Q}}^{\mathrm{dp}}(\mathcal{C}) \rightarrow \mathrm{coAlg}_{\tau^n\mathcal{Q}}^{\mathrm{dp}}(\mathrm{fil}(\mathcal{C}))
\end{equation*}
satisfies $a^k \cong 0$ for $k > n$.
\end{lemma}
\begin{proof}
Write
\begin{equation*}
\mathrm{fil}^{\geq k}(\mathcal{C}) := \mathrm{Fun}(\mathbb{Z}_{\geq k}, \mathcal{C})
\end{equation*}
for the $\infty$-category of filtered objects $X^k \rightarrow X^{k+1} \rightarrow \cdots$ starting in degree $k$. The inclusion $\iota\colon \mathbb{Z}_{\geq k} \rightarrow \mathbb{Z}_{\geq 1}$ is nonunital symmetric monoidal. Hence the induced restriction functor
\begin{equation*}
\iota^*\colon \mathrm{fil}(\mathcal{C}) \rightarrow \mathrm{fil}^{\geq k}(\mathcal{C})
\end{equation*}
is nonunital symmetric monoidal as well. Its left adjoint $\iota_!$ therefore acquires a nonunital oplax symmetric monoidal structure and the pair $(\iota_!,\iota^*)$ lifts to an adjoint pair on divided power coalgebras as follows (see Appendix \ref{app:functoriality}):
\[
\begin{tikzcd}
\mathrm{coAlg}^{\mathrm{dp}}_{\tau^n\mathcal{Q}}(\mathrm{fil}^{\geq k}(\mathcal{C})) \ar[shift left]{r}{\iota_!} & \mathrm{coAlg}^{\mathrm{dp}}_{\tau^n\mathcal{Q}}(\mathrm{fil}(\mathcal{C})). \ar[shift left]{l}{\iota^*}
\end{tikzcd}
\]
Here we have also used that both $\iota_!$ and $\iota^*$ preserve colimits. Now observe that for a cofree coalgebra $Y = \mathrm{cofree}^{\mathrm{dp}}_{\tau^n\mathcal{Q}}(V)$, the adic filtration $a^j Y$ vanishes for $j > n$ by Lemma \ref{lem:cofreetruncatedfiltcoalg}. In other words, $\iota^*(a^\bullet Y) \cong 0$. Any $Z \in \mathrm{coAlg}^{\mathrm{dp}}_{\tau^n\mathcal{Q}}(\mathcal{C})$ can be written as a limit of cofree coalgebras; indeed, by Proposition \ref{prop:coalgebrascomonadic}, the $\infty$-category $\mathrm{coAlg}^{\mathrm{dp}}_{\tau^n\mathcal{Q}}(\mathcal{C})$ is comonadic over $\mathcal{C}$. The functor $\iota^* \circ a^{\bullet}$ is a right adjoint, hence preserves limits, and we conclude that $\iota^*Z \cong 0$ as desired.
\end{proof}

Denote by
\begin{equation*}
\mathrm{gr}\colon \mathrm{fil}(\mathcal{C}) \rightarrow \mathrm{gr}(\mathcal{C})
\end{equation*}
the functor taking a filtered object $X_\bullet$ to its associated graded $\mathrm{gr}(X_\bullet)$, of which the $i$th component is
\begin{equation*}
\mathrm{gr}(X_\bullet)(i) := \mathrm{fib}(X^i \rightarrow X^{i+1}).
\end{equation*}
The functor $\mathrm{gr}$ is nonunital symmetric monoidal (cf. \cite[Proposition 3.2.1]{lurieRot}). Since $\mathcal{C}$ is stable, the functor $\mathrm{gr}$ preserves all limits and colimits. It has a left adjoint
\begin{equation*}
\delta\colon \mathrm{gr}(\mathcal{C}) \rightarrow \mathrm{fil}(\mathcal{C})
\end{equation*}
assigning to a graded object $\{X(i)\}_{i \geq 1}$ the `discrete' filtered object
\begin{equation*}
\delta(X) = (X(1) \xrightarrow{0} X(2) \xrightarrow{0} X(3) \xrightarrow{0} \cdots).
\end{equation*}
Being left adjoint to the nonunital symmetric monoidal functor $\mathrm{gr}$, the functor $\delta$ acquires a nonunital oplax symmetric monoidal structure. In particular, the adjunction $(\delta,\mathrm{gr})$ lifts to divided power $\mathcal{Q}$-coalgebras to give a square
\[
\begin{tikzcd}
\mathrm{coAlg}_{\mathcal{Q}}^{\mathrm{dp}}(\mathrm{gr}(\mathcal{C})) \ar{d}[swap]{\mathrm{forget}}\ar[shift left]{r}{\delta} & \mathrm{coAlg}_{\mathcal{Q}}^{\mathrm{dp}}(\mathrm{fil}(\mathcal{C})) \ar{l}[shift left]{\mathrm{gr}}\ar{d}{\mathrm{forget}} \\
\mathrm{gr}(\mathcal{C}) \ar[shift left]{r}{\delta} & \mathrm{fil}(\mathcal{C}) \ar{l}[shift left]{\mathrm{gr}}
\end{tikzcd}
\]
where both the left and right adjoints commute with the forgetful functors down to $\mathcal{C}$.

The following describes the associated graded of the adic filtration. It can be compared to the dual result for algebras  \cite[Proposition 4.26]{brantnermathew}, although the proof in loc. cit. does not quite apply here.

\begin{lemma}
\label{lem:assocgradedadic}
Let $Y \in \mathrm{coAlg}_{\mathcal{Q}}^{\mathrm{dp}}(\mathcal{C})$. Then there is a natural isomorphism of graded divided power coalgebras
\begin{equation*}
\mathrm{gr}(a^{\bullet} Y) \cong \mathrm{Sym}_{\mathcal{Q}}([\mathrm{prim}_{\mathcal{Q}} Y]_1).
\end{equation*}
In particular, there are natural isomorphisms $\mathrm{gr}(a^{\bullet} Y)(n) \cong D_n^{\mathcal{Q}}(\mathrm{prim}_{\mathcal{Q}} Y)$. (Here $[\mathrm{prim}_{\mathcal{Q}} Y]_1$ denotes the graded object of $\mathcal{C}$ that equals $\mathrm{prim}_{\mathcal{Q}} Y$ in degree 1 and is zero in all other degrees.)
\end{lemma}
\begin{proof}
The functor 
\begin{equation*}
\mathrm{gr} \circ a^{\bullet}\colon \mathrm{coAlg}_{\mathcal{Q}}^{\mathrm{dp}}(\mathcal{C}) \rightarrow \mathrm{coAlg}_{\mathcal{Q}}^{\mathrm{dp}}(\mathrm{gr}(\mathcal{C}))
\end{equation*} 
is right adjoint to the functor
\begin{equation*}
u \circ \delta\colon \mathrm{coAlg}_{\mathcal{Q}}^{\mathrm{dp}}(\mathrm{gr}(\mathcal{C})) \rightarrow \mathrm{coAlg}_{\mathcal{Q}}^{\mathrm{dp}}(\mathcal{C})
\end{equation*}
sending a graded coalgebra $X = \{X(i)\}_{i \geq 1}$ to its degree 1 part $X(1)$. This is a trivial coalgebra for evident degree reasons; in other words, $u \circ \delta$ is isomorphic to the functor $X \mapsto \mathrm{triv}_{\mathcal{Q}}((\mathrm{forget}_{\mathcal{Q}} X)(1))$. The right adjoint of this functor clearly is the assignment
\begin{equation*}
Y \mapsto \mathrm{cofree}_{\mathcal{Q}}^{\mathrm{dp}}([\mathrm{prim}_{\mathcal{Q}} Y]_1).
\end{equation*}
Now Lemma \ref{lem:cofreegradedcoalg} gives the desired formula for the cofree graded coalgebra.
\end{proof}

\section{Truncations of algebras and coalgebras}
\label{sec:truncationsalg}

In this section we discuss how the `approximations' $\varphi_n\mathcal{O}$ and $\tau_n\mathcal{O}$ of an operad $\mathcal{O}$ (and similarly the approximations $\varphi^n\mathcal{Q}$ and $\tau^n\mathcal{Q}$ of a cooperad $\mathcal{Q}$) induce corresponding approximations of $\mathcal{O}$-algebras (or of divided power $\mathcal{Q}$-coalgebras), as well as how these various approximations relate. The result is summarized in Theorem \ref{thm:filtrations}.

Restriction along the morphism $\mathcal{O} \rightarrow \tau_n\mathcal{O}$ defines a functor
\begin{equation*}
\mathrm{Alg}_{\tau_n\mathcal{O}}(\mathcal{C}) \rightarrow \mathrm{Alg}_{\mathcal{O}}(\mathcal{C})
\end{equation*}
that admits a left adjoint. For an $\mathcal{O}$-algebra $X$, we denote the unit of the resulting adjunction by
\begin{equation*}
X \rightarrow t_n X
\end{equation*}
and refer to $t_n X$ as the \emph{$n$-truncation} of the algebra $X$. Informally speaking, it is the quotient of the algebra $X$ obtained by killing all $k$-ary operations on $X$ for $k>n$. Dually, the morphism $\varphi_n\mathcal{O} \rightarrow \mathcal{O}$ determines a functor
\begin{equation*}
\mathrm{Alg}_{\mathcal{O}}(\mathcal{C}) \rightarrow \mathrm{Alg}_{\varphi_n\mathcal{O}}(\mathcal{C})
\end{equation*}
again admitting a left adjoint. We write
\begin{equation*}
f_n X \rightarrow X
\end{equation*}
for the counit of this adjunction. The functors $t_n$ and $f_n$ give two different approximations of $X$. Indeed, letting $n$ vary we obtain a diagram
\[
\begin{tikzcd}
\vdots \ar{ddr} && \vdots \ar{d} \\
f_2 X \ar{u}\ar{dr} && t_2 X \ar{d} \\
f_1 X \ar{u}\ar{r} & X \ar{r}\ar{ur}\ar{uur} & t_1 X.
\end{tikzcd}
\]
Thus we can think of the system of algebras $\{f_n X\}_{n \geq 1}$ as a filtration of $X$ and the system $\{t_n X\}_{n \geq 1}$ as a tower of algebras under $X$.

We introduce some more notation to describe the dual situation for coalgebras. For $\mathcal{Q}$ a cooperad, the morphisms $\tau^n\mathcal{Q} \rightarrow \mathcal{Q}$ and $\mathcal{Q} \rightarrow \varphi^n\mathcal{Q}$ induce adjunctions on the corresponding $\infty$-categories of coalgebras. For $Y \in \mathrm{coAlg}_{\mathcal{Q}}^{\mathrm{dp}}(\mathcal{C})$, we write
\begin{equation*}
t^n Y \rightarrow Y
\end{equation*}
for the counit of the first and
\begin{equation*}
Y \rightarrow f^n Y
\end{equation*}
for the unit of the second. For a coalgebra $Y$ we therefore have the following diagram of approximations:
\[
\begin{tikzcd}
\vdots \ar{ddr} && \vdots \ar{d} \\
t^2 Y \ar{u}\ar{dr} && f^2 Y \ar{d} \\
t^1 Y \ar{u}\ar{r} & Y \ar{r}\ar{ur}\ar{uur} & f^1 Y.
\end{tikzcd}
\]

\begin{definition}
For $X \in \mathrm{Alg}_{\mathcal{O}}(\mathcal{C})$, its \emph{nilcompletion} is the map
\begin{equation*}
X \rightarrow \varprojlim_n t_n X.
\end{equation*}
We say $X$ is \emph{nilcomplete} if this map is an isomorphism. Dually, for $Y \in \mathrm{coAlg}_{\mathcal{Q}}^{\mathrm{dp}}(\mathcal{C})$, its \emph{conilcompletion} is the map
\begin{equation*}
\varinjlim_n t^n Y \rightarrow Y
\end{equation*}
and $Y$ is \emph{conilcomplete} if this is an isomorphism.
\end{definition}

\begin{remark}
Despite the terminology, the object $\varprojlim_n t_n X$ need generally \emph{not} be nilcomplete itself. In other words, nilcompletion need not give an idempotent functor on the $\infty$-category of $\mathcal{O}$-algebras. This is similar to the behavior of other familiar completion functors, such as Bousfield--Kan completion of spaces, or (derived) completion of modules over a ring $R$ with respect to an ideal $I$.
\end{remark}

We summarize the crucial properties of these approximations in the following:

\begin{theorem}
\label{thm:filtrations}
\begin{itemize}
\item[(1)] For an $\mathcal{O}$-algebra $X$, the tower of algebras $t_n X$ has associated graded
\begin{equation*}
\mathrm{fib}(t_n X \rightarrow t_{n-1} X) \cong \mathrm{triv}_{\mathcal{O}}\bigl((\mathcal{O}(n) \otimes (\mathrm{indec}_{\mathcal{O}} X)^{\otimes n})_{h\Sigma_n}\bigr).
\end{equation*}
\item[(2)] The filtration of the identity functor of $\mathrm{Alg}_\mathcal{O}(\mathcal{C})$ by the functors $f_n$ has associated graded 
\begin{equation*}
\mathrm{cof}(f_{n-1} X \rightarrow f_n X) \cong \mathrm{free}_{\mathcal{O}}\bigl((B\mathcal{O}(n) \otimes X^{\otimes n})_{h\Sigma_n}\bigr).
\end{equation*}
Moreover, this filtration is exhaustive in the sense that the natural map
\begin{equation*}
\varinjlim_n f_n X \rightarrow X
\end{equation*}
is an isomorphism.
\item[(3)] For a divided power $\mathcal{Q}$-coalgebra $Y$, the tower $f^n Y$ has associated graded
\begin{equation*}
\mathrm{fib}(f^n Y \rightarrow f^{n-1} Y) \cong \mathrm{cofree}_{\mathcal{Q}}\bigl((C\mathcal{Q}(n) \otimes Y^{\otimes n})_{h\Sigma_n} \bigr).
\end{equation*}
Moreover, this tower is complete in the sense that the natural map
\begin{equation*}
Y \rightarrow \varprojlim_n f^n Y
\end{equation*}
is an isomorphism.
\item[(4)] The filtration of the identity of $\mathrm{coAlg}_{\mathcal{Q}}^{\mathrm{dp}}(\mathcal{C})$ by the functors $t^n$ has associated graded
\begin{equation*}
\mathrm{cof}(t^{n-1} Y \rightarrow t^n Y) \cong \mathrm{triv}_{\mathcal{Q}}\bigl((\mathcal{Q}(n) \otimes (\mathrm{prim}_{\mathcal{Q}} Y)^{\otimes n})_{h\Sigma_n}\bigr).
\end{equation*}
\end{itemize}
\end{theorem}
\begin{proof}
(1) The fiber of $t_n X \rightarrow t_{n-1} X$ can explicitly be computed as the fiber of the map
\begin{equation*}
B(\tau_n\mathcal{O},\mathcal{O}, X) \rightarrow B(\tau_{n-1}\mathcal{O},\mathcal{O}, X).
\end{equation*}
Since the composition product is exact in the first variable, this fiber is the bar construction
\begin{equation*}
B(\mathcal{O}(n),\mathcal{O}, X),
\end{equation*}
where $\mathcal{O}(n)$ here should be read as a symmetric sequence concentrated in degree $n$. For evident degree reasons its $\mathcal{O}$-module structure is trivial, so this is equivalent to
\begin{equation*}
\bigl(\mathcal{O}(n) \otimes B(\mathbf{1}, \mathcal{O}, X)^{\otimes n}\bigr)_{h\Sigma_n} = D_n^{\mathcal{O}}(B(\mathbf{1}, \mathcal{O}, X)).
\end{equation*}
The term $B(\mathbf{1}, \mathcal{O}, X)$ is a formula for $\mathrm{indec}_{\mathcal{O}} X$. (Note that we have implicitly used that the functor $D_n^{\mathcal{O}}$ preserves the geometric realization computing the bar construction here.)


(2) Let $X$ and $Y$ be $\mathcal{O}$-algebras. To prove the first statement, it suffices to establish a natural isomorphism between the fiber $F$ of the map
\begin{equation*}
\mathrm{Map}_{\mathrm{Alg}_{\mathcal{O}}(\mathcal{C})}(f_n X, Y) \rightarrow \mathrm{Map}_{\mathrm{Alg}_{\mathcal{O}}(\mathcal{C})}(f_{n-1} X, Y)
\end{equation*}
and the space $\mathrm{Map}_{\mathcal{C}}((B\mathcal{O}(n) \otimes X^{\otimes n})_{h\Sigma_n}, Y)$. Theorem \ref{thm:algdecomposition} implies that $F$ is isomorphic to the fiber of the middle vertical map in the following diagram, in which all rows and columns are fiber sequences:
\[
\begin{tikzcd}
G \ar{d}\ar{r}{\cong} & F \ar{d}\ar{r} & * \ar{d} \\
\mathrm{Map}_{\mathcal{C}}(D_n^{\mathcal{O}}(X),\Omega Y) \ar{r}\ar{d} & \mathrm{Map}_{\mathrm{Alg}_{D_n^{\mathcal{O}}}}(X,Y) \ar{r}\ar{d} & \mathrm{Map}_{\mathcal{C}}(X,Y) \ar{d} \\
\mathrm{Map}_{\mathcal{C}}(D_n^{\varphi_{n-1}\mathcal{O}}(X), \Omega Y) \ar{r} & \mathrm{Map}_{\mathrm{Alg}_{D_n^{\varphi_{n-1}\mathcal{O}}}}(X,Y) \ar{r} & \mathrm{Map}_{\mathcal{C}}(X,Y).
\end{tikzcd}
\]
The fiber $G$ on the left may be identified as the space
\begin{equation*}
\mathrm{Map}_{\mathcal{C}}((\mathrm{cof}(\varphi_{n-1}\mathcal{O}(n) \rightarrow \mathcal{O}(n)) \otimes X^{\otimes n})_{h\Sigma_n}, \Omega Y) \cong \mathrm{Map}_{\mathcal{C}}((B\mathcal{O}(n) \otimes X^{\otimes n})_{h\Sigma_n}, Y),
\end{equation*}
where the expression on the right is deduced from Remark \ref{rmk:formulabarconstr} and the adjunction between $\Omega$ and $\Sigma$. This completes the argument. 

The fact that the filtration by the $f_n X$ is exhaustive is now a consequence of the second half of Theorem \ref{thm:algdecomposition}. Alternatively, one can directly observe the following natural isomorphisms:
\begin{eqnarray*}
\varinjlim_n f_n X & \cong & \varinjlim_n B(\mathcal{O}, \varphi_n\mathcal{O}, X) \\
& \cong & B(\mathcal{O}, \varinjlim_n \varphi_n\mathcal{O}, X) \\
& \cong & B(\mathcal{O}, \mathcal{O}, X) \\
& \cong & X.
\end{eqnarray*}

(3) The argument is entirely dual to that for (2), now relying on Theorem \ref{thm:coalgdpdecomposition}.

(4) We could start the proof in a style dual to the proof of (1), but we run into the issue of showing that the map
\begin{equation*}
D_n^{\mathcal{Q}}(\mathrm{prim}_{\mathcal{Q}} Y) \rightarrow \mathrm{Tot}(D_n^\mathcal{Q}(C^\bullet(\mathrm{id}_{\mathcal{C}}, \mathrm{cofree}_{\mathcal{Q}},Y)))
\end{equation*}
is an isomorphism. This is not immediately clear, since the functor $D_n^\mathcal{Q}$ need not preserve totalizations in general. To circumvent the issue we give an argument using filtered coalgebras. The adic filtration of divided power $\tau^n\mathcal{Q}$-coalgebras provides a filtration of the object $t^n Y$ of length $n$ (cf. Lemma \ref{lem:finiteadicfiltration}). Similarly, the adic filtration of $\tau^{n-1}\mathcal{Q}$-coalgebras gives $t^{n-1} Y$ a filtration of length $n-1$. On associated graded objects, Lemma \ref{lem:assocgradedadic} shows that the map $t^{n-1} Y \rightarrow t^n Y$ (now interpreted as a map of filtered coalgebras) gives the inclusion
\begin{equation*}
\mathrm{Sym}_{\tau^{n-1}\mathcal{Q}}(Y) \rightarrow \mathrm{Sym}_{\tau^n\mathcal{Q}}(Y).
\end{equation*}
We conclude that the filtered divided power $\mathcal{Q}$-coalgebra $\mathrm{cof}(t^{n-1} Y \rightarrow t^n Y)$ vanishes in degrees $>n$ and is constant with value $D_n^\mathcal{Q}(\mathrm{prim}_{\mathcal{Q}}(Y))$ in degrees $\leq n$. It is a trivial coalgebra for degree reasons. Taking underlying objects of these filtered coalgebras completes the proof.
\end{proof}

\begin{remark}
The filtration $f_n X$ of an $\mathcal{O}$-algebra $X$ should be thought of as the Koszul dual of the filtration $t^n Y$ of a $B\mathcal{O}$-coalgebra $Y$ (and similarly for $t_n X$ and $f^n Y$). Indeed, the $f_n X$ give a filtration $\mathrm{indec}_{\mathcal{O}}(f_n X)$ of the indecomposables $\mathrm{indec}_{\mathcal{O}} X$ with associated graded
\begin{equation*}
\mathrm{cof}(\mathrm{indec}_{\mathcal{O}}(f_{n-1} X) \rightarrow \mathrm{indec}_{\mathcal{O}}(f_n X)) \cong \mathrm{triv}_{B\mathcal{O}}\bigl((B\mathcal{O}(n) \otimes X^{\otimes n})_{h\Sigma_n}\bigr), 
\end{equation*}
where we have used Theorem \ref{thm:filtrations}(1) and the fact that $\mathrm{indec}_{\mathcal{O}} \circ \mathrm{free}_{\mathcal{O}}$ is equivalent to the identity functor of $\mathcal{C}$. We will study the interactions of these filtrations with Koszul duality in more detail in the next section. 
\end{remark}

\begin{remark}
In the specific case where $\mathcal{O}$ is the (nonunital) commutative operad, the functor $\mathrm{indec}_{\mathcal{O}}$ is more commonly known as topological Andr\'{e}--Quillen homology ($TAQ$). The filtration of the previous remark then specializes in the following way: for an augmented commutative ring spectrum $R$, the spectrum $TAQ(R)$ has a filtration with associated graded consisting of the terms $(\mathbf{L}(n)^{\vee} \otimes R^{\otimes n})_{h\Sigma_n}$. Here $\mathbf{L}(n)$ is the $n$th term of the spectral Lie operad and $\mathbf{L}(n)^{\vee}$ is its Spanier--Whitehead dual. This filtration of $TAQ$ has been studied by Kuhn \cite{kuhnmccord} and by Behrens--Rezk \cite[Section 4]{behrensrezk}, who refer to it as the \emph{Kuhn filtration}.
\end{remark}

\begin{remark}
\label{rmk:adictn}
It is not difficult to expand the argument proving Theorem \ref{thm:filtrations}(4) to show that for $Y \in \mathrm{coAlg}_{\mathcal{Q}}^{\mathrm{dp}}(\mathcal{C})$, there is a natural cofiber sequence $t^n Y \to Y \to a^n Y$ in $\mathcal{C}$. It follows that $Y$ is conilcomplete (i.e., $\varinjlim_n t^n Y \to Y$ is an isomorphism) if and only if its adic filtration is cocomplete (in the sense that $\varinjlim_n a^n Y \cong 0$). A dual remark applies to algebras, although we have not dealt with the dual case of adic filtrations for algebras in this paper (but see \cite{brantnermathew}).
\end{remark}

\section{Nilcompletion and conilcompletion}
\label{sec:proof}

The goal of this section is to prove Theorem \ref{thm:completeness}, stating that the adjoint pair
\[
\begin{tikzcd}[column sep = large]
\mathrm{Alg}_{\mathcal{O}}(\mathcal{C}) \ar[shift left]{r}{\mathrm{indec}_{\mathcal{O}}} & \mathrm{coAlg}^{\mathrm{dp}}_{B\mathcal{O}}(\mathcal{C}) \ar[shift left]{l}{\mathrm{prim}_{B\mathcal{O}}}
\end{tikzcd}
\]
restricts to an equivalence between nilcomplete algebras and conilcomplete coalgebras. Essentially, the proof will amount to showing that Koszul duality between $\mathcal{O}$-algebras and $B\mathcal{O}$-coalgebras intertwines the filtrations of Theorem \ref{thm:filtrations} in an appropriate sense, see Corollary \ref{cor:nilpotentunit}. Theorem \ref{thm:completeness} is then a direct consequence of Proposition \ref{prop:unit} below.

We begin by analyzing the truncated case:

\begin{proposition}
\label{prop:nilpotentunit}
For an operad $\mathcal{O}$ and any $n \geq 1$ the functor
\begin{equation*}
\mathrm{indec}_{\tau_n\mathcal{O}}\colon \mathrm{Alg}_{\tau_n\mathcal{O}}(\mathcal{C}) \rightarrow \mathrm{coAlg}_{B(\tau_n\mathcal{O})}^{\mathrm{dp}}(\mathcal{C})
\end{equation*}
is fully faithful. In other words, for any $X \in \mathrm{Alg}_{\tau_n\mathcal{O}}(\mathcal{C})$ the unit map
\begin{equation*}
X \rightarrow \mathrm{prim}_{B(\tau_n\mathcal{O})}\mathrm{indec}_{\tau_n\mathcal{O}}(X)
\end{equation*}
is an isomorphism. Dually, for a cooperad $\mathcal{Q}$ the functor
\begin{equation*}
\mathrm{prim}_{\tau^n\mathcal{Q}}\colon \mathrm{coAlg}^{\mathrm{dp}}_{\tau^n\mathcal{Q}}(\mathcal{C}) \rightarrow \mathrm{Alg}_{C(\tau^n\mathcal{Q})}(\mathcal{C})
\end{equation*}
is fully faithful. 
\end{proposition}
\begin{remark}
A version of this result (at least the first half for operads) is proved as Proposition 6.9 of \cite{heutsgoodwillie}.
\end{remark}
\begin{proof}
Given our work in Theorem \ref{thm:filtrations} the proofs of the two halves of the proposition are entirely dual, so we spell out the first case. Since the operad $\tau_n\mathcal{O}$ is $n$-truncated, any $X \in \mathrm{Alg}_{\tau_n\mathcal{O}}(\mathcal{C})$ has a \emph{finite} tower of `nilpotent approximations'
\begin{equation*}
X = t_n X \rightarrow t_{n-1} X \rightarrow \cdots \rightarrow t_1 X.
\end{equation*}
Theorem \ref{thm:filtrations}(1) states that the associated graded object is described (as an object of $\mathcal{C}$) by
\begin{equation*}
G(X) = \bigoplus_{k=1}^n \bigl(\mathcal{O}(k) \otimes (\mathrm{indec}_{\tau_n \mathcal{O}} X)^{\otimes k}\bigr)_{h\Sigma_k}.
\end{equation*}
For $k \leq n$, consider the commutative square of left adjoints
\[
\begin{tikzcd}[column sep = large]
\mathrm{Alg}_{\tau_n\mathcal{O}}(\mathcal{C}) \ar{r}{\mathrm{indec}_{\tau_n\mathcal{O}}}\ar{d} & \mathrm{coAlg}^{\mathrm{dp}}_{B(\tau_n\mathcal{O})}(\mathcal{C}) \ar{d} \\
\mathrm{Alg}_{\tau_k\mathcal{O}}(\mathcal{C}) \ar{r}{\mathrm{indec}_{\tau_k\mathcal{O}}} & \mathrm{coAlg}^{\mathrm{dp}}_{B(\tau_k\mathcal{O})}(\mathcal{C})
\end{tikzcd}
\]
with vertical functors induced by the map of operads $\tau_n\mathcal{O} \rightarrow \tau_k\mathcal{O}$ and the map of cooperads $B(\tau_n\mathcal{O}) \rightarrow B(\tau_k\mathcal{O})$. The map $X \rightarrow t_k X$ is the unit of the left vertical adjunction. The unit of the composite left adjoint
\begin{equation*}
\mathrm{Alg}_{\tau_n\mathcal{O}}(\mathcal{C}) \rightarrow \mathrm{coAlg}_{B(\tau_k\mathcal{O})}^{\mathrm{dp}}(\mathcal{C}) 
\end{equation*}
is the natural map
\begin{equation*}
X \rightarrow \mathrm{prim}_{B(\tau_n\mathcal{O})} f^k \mathrm{indec}_{\tau_n\mathcal{O}} X.
\end{equation*}
The commutativity of the square implies that the tower $\{t_k X\}_{k \geq 1}$ under $X$ naturally maps to the tower 
\begin{equation*}
\mathrm{prim}_{B(\tau_n\mathcal{O})} \mathrm{indec}_{\tau_n\mathcal{O}} X = \mathrm{prim}_{B(\tau_n\mathcal{O})} f^n \mathrm{indec}_{\tau_n\mathcal{O}} X \rightarrow \cdots \rightarrow \mathrm{prim}_{B(\tau_n\mathcal{O})} f^1 \mathrm{indec}_{\tau_n\mathcal{O}} X.
\end{equation*}
This second tower has associated graded
\begin{equation*}
G'(X) = \bigoplus_{k=1}^n \bigl(CB(\tau_n\mathcal{O})(k) \otimes( \mathrm{indec}_{\tau_n \mathcal{O}} X)^{\otimes k}\bigr)_{h\Sigma_k}
\end{equation*}
by Theorem \ref{thm:filtrations}(3) and the fact that $\mathrm{prim}_{B(\tau_n\mathcal{O})}\mathrm{cofree}_{B(\tau_n\mathcal{O})} \cong \mathrm{triv}_{\tau_n\mathcal{O}}$. The map $G \rightarrow G'$ is a sum of homogeneous pieces and hence induced by a map of symmetric sequences
\begin{equation*}
\tau_n\mathcal{O} \rightarrow CB(\tau_n\mathcal{O}).
\end{equation*}
We claim that it is the unit map of the bar-cobar adjunction for operads and hence an equivalence by Theorem \ref{thm:barcobaroperads}, completing the proof. To verify the claim, let us abbreviate notation by $\mathcal{P} = \tau_n\mathcal{O}$ and evaluate the natural map $G(X) \to G'(X)$ at a free $\mathcal{P}$-algebra $X = \mathrm{Sym}_{\mathcal{P}}(V)$. Then $G(X)$ arises as the associated graded of the evident finite tower $\{\tau_k \mathcal{P} \circ V\}_{k \leq n}$. Observing that $\mathrm{indec}_{\mathcal{P}} X \cong \mathrm{triv}_{B\mathcal{P}} V$, we see that $G'(X)$ arises from the tower of which the $k$th stage is
\[
C(\mathbf{1}, B\mathcal{P} ,C(B\mathcal{P}, \varphi^k B\mathcal{P}, \mathbf{1})) \circ V \cong C(\mathbf{1}, \varphi^k B\mathcal{P}, \mathbf{1}) \circ V.
\]
Inspecting the square diagram above and identifying $\varphi^k B\mathcal{P}$ with $B(\tau_k\mathcal{P})$, we see that $G \to G'$ indeed arises from the map
\[
\tau_K\mathcal{P} \to C(\mathbf{1}, B(\tau_k\mathcal{P}), \mathbf{1})
\]
considered in Theorem \ref{thm:barcobaroperads}.
\end{proof}

\begin{corollary}
\label{cor:nilpotentunit}
For an operad $\mathcal{O}$ and $X \in \mathrm{Alg}_{\mathcal{O}}(\mathcal{C})$ there are natural isomorphisms
\begin{equation*}
t_n X \xrightarrow{\cong} \mathrm{prim}_{B\mathcal{O}}(f^n \mathrm{indec}_{\mathcal{O}} X)
\end{equation*}
under $X$. Dually, for a cooperad $\mathcal{Q}$ and $Y \in \mathrm{coAlg}_{B\mathcal{O}}^{\mathrm{dp}}(\mathcal{C})$ there are natural isomorphisms
\begin{equation*}
t^n Y \xrightarrow{\cong} \mathrm{indec}_{C\mathcal{Q}}(f_n \mathrm{prim}_{\mathcal{Q}} Y).
\end{equation*}
over $Y$.
\end{corollary}
\begin{proof}
Again we only prove the first half, the second half being dual. Consider the square of left adjoints
\[
\begin{tikzcd}[column sep = large]
\mathrm{Alg}_{\mathcal{O}}(\mathcal{C}) \ar{r}{\mathrm{indec}_{\mathcal{O}}}\ar{d} & \mathrm{coAlg}^{\mathrm{dp}}_{B\mathcal{O}}(\mathcal{C}) \ar{d} \\
\mathrm{Alg}_{\tau_n\mathcal{O}}(\mathcal{C}) \ar{r}{\mathrm{indec}_{\tau_n\mathcal{O}}} & \mathrm{coAlg}^{\mathrm{dp}}_{B(\tau_n\mathcal{O})}(\mathcal{C}).
\end{tikzcd}
\]
The map $X \rightarrow t_n X$ is the unit for the left vertical adjunction. Since the unit for the bottom horizontal adjunction is an isomorphism by Proposition \ref{prop:nilpotentunit}, this map is isomorphic to the unit for the composite adjunction from top left to bottom right. Chasing around the square along the top and right, we find the unit map
\begin{equation*}
X \rightarrow  \mathrm{prim}_{B\mathcal{O}}(f^n \mathrm{indec}_{\mathcal{O}} X).
\end{equation*}
\end{proof}

\begin{proposition}
\label{prop:unit}
For $X \in \mathrm{Alg}_{\mathcal{O}}(\mathcal{C})$, the unit map $X \rightarrow \mathrm{prim}_{B\mathcal{O}}\mathrm{indec}_{\mathcal{O}} X$ is naturally isomorphic to the nilcompletion
\begin{equation*}
X \rightarrow \varprojlim_n t_n X.
\end{equation*}
Dually, for $Y \in \mathrm{coAlg}^{\mathrm{dp}}_{\mathcal{Q}}(\mathcal{C})$, the counit map $\mathrm{indec}_{C\mathcal{Q}}\mathrm{prim}_{\mathcal{Q}} Y \rightarrow Y$ is naturally isomorphic to the conilcompletion
\begin{equation*}
\varinjlim_n t^n Y \rightarrow Y.
\end{equation*}
\end{proposition}
\begin{proof}
Again we prove only the first half, the second being dual. We have natural isomorphisms
\begin{equation*}
\varprojlim_n t_n X \cong \varprojlim_n \mathrm{prim}_{B\mathcal{O}}(f^n\mathrm{indec}_{\mathcal{O}} X)
\end{equation*}
under $X$ by Corollary \ref{cor:nilpotentunit} and natural isomorphisms
\begin{equation*}
\varprojlim_n \mathrm{prim}_{B\mathcal{O}}(f^n\mathrm{indec}_{\mathcal{O}} X) \cong \mathrm{prim}_{B\mathcal{O}}(\varprojlim_n f^n\mathrm{indec}_{\mathcal{O}} X) \cong \mathrm{prim}_{B\mathcal{O}}(\mathrm{indec}_{\mathcal{O}} X)
\end{equation*}
since $\mathrm{prim}_{B\mathcal{O}}$ preserves limits (being a right adjoint) and from the completeness assertion of Theorem \ref{thm:filtrations}(3). Combining these yields the proposition.
\end{proof}

\section{Nilcompletion and homological completion}
\label{sec:completions}

Recall that we say an operad $\mathcal{O}$ has \emph{good completion} if every trivial $\mathcal{O}$-algebra is nilcomplete, i.e., if for every $X \in \mathcal{C}$ the map
\begin{equation*}
\mathrm{triv}_{\mathcal{O}} X \rightarrow \varprojlim_n t_n \mathrm{triv}_{\mathcal{O}} X
\end{equation*}
is an isomorphism. In this section we show that if $\mathcal{O}$ has good completion, then the notions of nilcompletion and homological completion of $\mathcal{O}$-algebras agree. In Section \ref{sec:FGconjecture} we will use this fact to show how our main result (Theorem \ref{thm:completeness}) can be used to find counterexamples to the conjecture of Francis--Gaitsgory. In the subsequent Section \ref{sec:goodcompletion} we provide a source of examples of operads with good completion.

The left adjoint functors
\begin{equation*}
\mathrm{Alg}_{\mathcal{O}}(\mathcal{C}) \xrightarrow{\mathrm{indec}^{\mathrm{nil}}_{\mathcal{O}}} \mathrm{coAlg}^{\mathrm{dp,nil}}_{B\mathcal{O}}(\mathcal{C}) \rightarrow \mathrm{coAlg}^{\mathrm{dp}}_{B\mathcal{O}}(\mathcal{C}) 
\end{equation*}
induce, for $X \in \mathrm{Alg}_{\mathcal{O}}(\mathcal{C})$, two unit maps
\begin{equation*}
X \rightarrow \mathrm{prim}^{\mathrm{nil}}_{B\mathcal{O}}\mathrm{indec}^{\mathrm{nil}}_{\mathcal{O}}(X) \rightarrow \mathrm{prim}_{B\mathcal{O}}\mathrm{indec}_{\mathcal{O}}(X).
\end{equation*}
By Theorem \ref{thm:completeness} the term $\mathrm{prim}_{B\mathcal{O}}\mathrm{indec}_{\mathcal{O}}(X)$ is naturally isomorphic to the nilcompletion of $X$, whereas we dubbed $\mathrm{prim}^{\mathrm{nil}}_{B\mathcal{O}}\mathrm{indec}^{\mathrm{nil}}_{\mathcal{O}}(X)$ the \emph{homological completion} of $X$. We promised the following result in Section \ref{sec:mainresults}:

\begin{proposition}
\label{prop:goodcompletion}
If the operad $\mathcal{O}$ has good completion, then the natural map
\begin{equation*}
\mathrm{prim}^{\mathrm{nil}}_{B\mathcal{O}}\mathrm{indec}^{\mathrm{nil}}_{\mathcal{O}}(X) \rightarrow \mathrm{prim}_{B\mathcal{O}}\mathrm{indec}_{\mathcal{O}}(X)
\end{equation*}
is an isomorphism for any $\mathcal{O}$-algebra $X$. In particular, $X$ is homologically complete if and only if it is nilcomplete.
\end{proposition}
\begin{proof}
By Lemma \ref{lem:primnilresolution} the homological completion $\mathrm{prim}^{\mathrm{nil}}_{B\mathcal{O}}\mathrm{indec}^{\mathrm{nil}}_{\mathcal{O}}(X)$ may be computed as the totalization of the cosimplicial object
\begin{equation*}
(\mathrm{triv}_{\mathcal{O}}\mathrm{cot}_{\mathcal{O}})^{\bullet + 1} X.
\end{equation*}
In other words, the homological completion is just the completion in the sense of Bousfield--Kan with respect to the homology theory $\mathrm{cot}_{\mathcal{O}}$. We claim that for any $n \geq 1$, the map
\begin{equation*}
t_n X \rightarrow \mathrm{Tot}\bigl(t_n(\mathrm{triv}_{\mathcal{O}}\mathrm{cot}_{\mathcal{O}})^{\bullet + 1}X\bigr)
\end{equation*}
is an equivalence. For $n=1$ we have $t_1 = \mathrm{triv}_{\mathcal{O}}\mathrm{cot}_{\mathcal{O}}$ and the cosimplicial object on the right admits a contracting codegeneracy onto $t_1 X$. For $n > 1$ we reason by induction, using the following diagram of fiber sequences (cf. Theorem \ref{thm:filtrations}(1)):
\[
\begin{tikzcd}
\mathrm{triv}_{\mathcal{O}}D_n^{\mathcal{O}}(\mathrm{cot}_{\mathcal{O}} X) \ar{r}\ar{d} & \mathrm{Tot}\bigl(\mathrm{triv}_{\mathcal{O}}D_n^{\mathcal{O}}(\mathrm{cot}_{\mathcal{O}}(\mathrm{triv}_{\mathcal{O}}\mathrm{cot}_{\mathcal{O}})^{\bullet + 1}X)\bigr) \ar{d} \\
t_n X \ar{r}\ar{d} & \mathrm{Tot}\bigl(t_n(\mathrm{triv}_{\mathcal{O}}\mathrm{cot}_{\mathcal{O}})^{\bullet + 1}X\bigr) \ar{d} \\
t_{n-1} X \ar{r} & \mathrm{Tot}\bigl(t_{n-1}(\mathrm{triv}_{\mathcal{O}}\mathrm{cot}_{\mathcal{O}})^{\bullet + 1}X\bigr).
\end{tikzcd}
\]
The bottom arrow is an isomorphism by induction, while the top is again an isomorphism by the fact that $\mathrm{cot}_{\mathcal{O}}(\mathrm{triv}_{\mathcal{O}}\mathrm{cot}_{\mathcal{O}})^{\bullet + 1}X$ admits a contracting codegeneracy to $\mathrm{cot}_{\mathcal{O}} X$.

Now consider the following commutative square:
\[
\begin{tikzcd}
X \ar{r}\ar{d} & \varprojlim_n t_n X \ar{d} \\
\mathrm{Tot}\bigl((\mathrm{triv}_{\mathcal{O}}\mathrm{cot}_{\mathcal{O}})^{\bullet + 1} X \bigr) \ar{r} & \varprojlim_n\mathrm{Tot} \bigl(t_n(\mathrm{triv}_{\mathcal{O}}\mathrm{cot}_{\mathcal{O}})^{\bullet + 1} X \bigr).
\end{tikzcd}
\]
We have assumed that $\mathrm{triv}_{\mathcal{O}} \cong \varprojlim_n t_n \mathrm{triv}_{\mathcal{O}}$, so the bottom horizontal arrow is an isomorphism. The vertical arrow on the right is an isomorphism by the first part of the proof. It now follows that the top horizontal map is an isomorphism if and only if the map on the left is an isomorphism; in other words, $X$ is nilcomplete if and only if it is homologically complete.
\end{proof}

\begin{remark}
\label{rmk:comparingcoalgs}
In general it is unclear to us to what extent the comparison functor
\begin{equation*}
\mathrm{coAlg}^{\mathrm{dp,nil}}_{B\mathcal{O}}(\mathcal{C}) \xrightarrow{\gamma} \mathrm{coAlg}^{\mathrm{dp}}_{B\mathcal{O}}(\mathcal{C})
\end{equation*}
from conilpotent divided power $B\mathcal{O}$-coalgebras to divided power $B\mathcal{O}$-coalgebras is fully faithful (but see Section \ref{sec:truncatedcompleteness} for some related observations). However, the hypothesis that $\mathcal{O}$ has good completion has some useful consequences. Indeed, under this hypothesis Proposition \ref{prop:goodcompletion} implies that $\gamma$ is fully faithful on the essential image of $\mathrm{indec}^{\mathrm{nil}}_{\mathcal{O}}$. To see this, consider the commutative square of mapping spaces\[
\begin{tikzcd}
\mathrm{Map}_{\mathrm{coAlg}^{\mathrm{dp,nil}}_{B\mathcal{O}}}(\mathrm{indec}^{\mathrm{nil}}_{\mathcal{O}} X, \mathrm{indec}^{\mathrm{nil}}_{\mathcal{O}} Y) \ar{r}{\gamma}\ar{d}{\cong} & \mathrm{Map}_{\mathrm{coAlg}^{\mathrm{dp}}_{B\mathcal{O}}}(\mathrm{indec}_{\mathcal{O}} X, \mathrm{indec}_{\mathcal{O}} Y) \ar{d}{\cong} \\
\mathrm{Map}_{\mathrm{Alg}(\mathcal{O})}(X, \mathrm{prim}^{\mathrm{nil}}_{B\mathcal{O}}\mathrm{indec}^{\mathrm{nil}}_{\mathcal{O}} Y) \ar{r} & \mathrm{Map}_{\mathrm{Alg}(\mathcal{O})}(X, \mathrm{prim}_{B\mathcal{O}}\mathrm{indec}_{\mathcal{O}} Y)
\end{tikzcd}
\]
and apply the proposition to conclude that the bottom horizontal map is an isomorphism. Note that the essential image of $\mathrm{indec}^{\mathrm{nil}}_{\mathcal{O}}$ includes both cofree and trivial conilpotent $B\mathcal{O}$-coalgebras, since $\mathrm{indec}^{\mathrm{nil}}_{\mathcal{O}} \circ \mathrm{triv}_{\mathcal{O}} \cong \mathrm{cofree}^{\mathrm{nil}}_{B\mathcal{O}}$ and $\mathrm{indec}^{\mathrm{nil}}_{\mathcal{O}} \circ \mathrm{free}_{\mathcal{O}} \cong \mathrm{triv}^{\mathrm{nil}}_{B\mathcal{O}}$.
\end{remark}

\section{The conjecture of Francis--Gaitsgory}
\label{sec:FGconjecture}

In this section we explain how our main result, Theorem \ref{thm:completeness}, can be used to find counterexamples to Conjecture \ref{conj:francisgaitsgory}. An $\mathcal{O}$-algebra $X$ is called \emph{nilpotent} if there exists an $n \geq 1$ so that $X$ is in the essential image of the restriction functor $\mathrm{Alg}_{\tau_n\mathcal{O}}(\mathcal{C}) \to \mathrm{Alg}_{\mathcal{O}}(\mathcal{C})$. An $\mathcal{O}$-algebra is called \emph{pronilpotent} if it can be obtained as a limit of nilpotent $\mathcal{O}$-algebras. Let us begin by observing the following straightforward consequence of Proposition \ref{prop:goodcompletion}.

\begin{corollary}
\label{cor:goodcompletion}
Let $\mathcal{O}$ be an operad with good completion. If $\mathcal{O}$ satisfies Conjecture \ref{conj:francisgaitsgory}, then every pronilpotent $\mathcal{O}$-algebra is nilcomplete.
\end{corollary}
\begin{proof}
Conjecture \ref{conj:francisgaitsgory} in particular implies that every pronilpotent $\mathcal{O}$-algebra is homologically complete. The statement then follows immediately from Proposition \ref{prop:goodcompletion}.
\end{proof}

Therefore, one may disprove Conjecture \ref{conj:francisgaitsgory} by finding an operad $\mathcal{O}$ with good completion and a pronilpotent $\mathcal{O}$-algebra that is \emph{not} nilcomplete. For an explicit such example, take $k$ to be a field of characteristic zero and $\mathcal{C}$ the $\infty$-category $\mathrm{Mod}_{k}$ of $k$-modules.  Take $\mathcal{O}$ to be the nonunital commutative operad in $\mathcal{C}$. We will demonstrate in the next section that this operad has good completion (see Theorem \ref{thm:goodcompletionchar0}). For ease of exposition, let us identify the $\infty$-category $\mathrm{Alg}_{\mathcal{O}}(\mathrm{Mod}_{k})$ of nonunital commutative $k$-algebras with that of augmented commutative $k$-algebras via the functor taking the augmentation ideal:
\begin{equation*}
\mathrm{CAlg}^{\mathrm{aug}}(\mathrm{Mod}_{k}) \xrightarrow{\cong} \mathrm{Alg}_{\mathcal{O}}(\mathrm{Mod}_{k})\colon (A \xrightarrow{\varepsilon} k) \mapsto \mathrm{fib}(\varepsilon).
\end{equation*}
The pronilpotent algebra we will use as a counterexample is a power series ring on a countable set of generators of degree 0:
\begin{equation*}
R := k\llbracket x_1, x_2, \ldots \rrbracket.
\end{equation*}
Its augmentation ideal consists of power series without constant term. Note that $R$ is indeed pronilpotent: it arises as the nilcompletion of the free commutative algebra $k[x_1, x_2, \ldots]$ and is therefore an inverse limit of nilpotent commutative algebras. To evaluate the nilcompletion of $R$ we need a basic observation:

\begin{lemma}
\label{lem:pi0nilcompletion}
Let $A$ be a connective augmented commutative $k$-algebra. Then for every $n$ there is a canonical isomorphism 
\begin{equation*}
\pi_0(t_n A) \cong \pi_0(A)/\mathfrak{m}^{n+1},
\end{equation*}
where $\mathfrak{m}$ denotes the augmentation ideal of $\pi_0(A)$.
\end{lemma}
\begin{proof}
The functor $t_n\colon \mathrm{CAlg}^{\mathrm{aug}}(\mathrm{Mod}_{k}) \rightarrow \mathrm{CAlg}^{\mathrm{aug}}(\mathrm{Mod}_{k})$ preserves sifted colimits and is therefore completely determined by its behavior on free algebras, i.e., algebras of the form $\mathrm{Sym}(V)$ for $V \in \mathrm{Mod}_k$. On such algebras it acts by
\begin{equation*}
t_n \mathrm{Sym}(V) \cong \bigoplus_{i \leq n} \mathrm{Sym}^i(V).
\end{equation*}
In particular, it is clear that $t_n$ preserves connectivity and hence restricts to a functor on connective algebras
\begin{equation*}
\mathrm{CAlg}^{\mathrm{aug}}(\mathrm{Mod}_{k, \geq 0}) \xrightarrow{t_n} \mathrm{CAlg}^{\mathrm{aug}}(\mathrm{Mod}_{k, \geq 0})
\end{equation*}
that still preserves sifted colimits. Writing $\mathrm{Vect}_k$ for the (ordinary) category of vector spaces over $k$, we may now postcompose with the symmetric monoidal colimit-preserving functor $\pi_0\colon \mathrm{Mod}_{k, \geq 0} \rightarrow \mathrm{Vect}_k$ to obtain a functor
\begin{equation*}
\mathrm{CAlg}^{\mathrm{aug}}(\mathrm{Mod}_{k, \geq 0}) \rightarrow \mathrm{CAlg}^{\mathrm{aug}}(\mathrm{Vect}_k)\colon A \mapsto \pi_0(t_n A),
\end{equation*}
again preserving sifted colimits. Since $A \mapsto \pi_0(A)/\mathfrak{m}^{n+1}$ preserves sifted colimits as well, it suffices to compare the two functors on free algebras $A = \mathrm{Sym}(V)$, for $V$ a connective $k$-module. But then we clearly have natural isomorphisms of (discrete) $k$-algebras
\begin{equation*}
\pi_0(t_n \mathrm{Sym}(V)) \cong \bigoplus_{i \leq n} \mathrm{Sym}^i(\pi_0 V) \cong \pi_0(\mathrm{Sym}(V))/\mathfrak{m}^{n+1}
\end{equation*}
and the conclusion follows.
\end{proof}

Proposition \ref{prop:counterexample} is now a consequence of the following exercise in commutative algebra: 

\begin{lemma}
The commutative algebra $R = k\llbracket x_1, x_2, \ldots \rrbracket$ is not nilcomplete.
\end{lemma}
\begin{proof}
The Milnor sequence in particular implies that the natural map 
\begin{equation*}
\pi_0(\varprojlim_n t_n R) \rightarrow \varprojlim_n \pi_0(t_n R)
\end{equation*}
is surjective. By Lemma \ref{lem:pi0nilcompletion} we may identify $\varprojlim_n \pi_0(t_n R)$ with the completion $\pi_0(R)^{\wedge}_{\mathfrak{m}}$ of $\pi_0(R)$ at its augmentation ideal. To prove that $\pi_0 R \rightarrow \pi_0(\varprojlim_n t_n R)$ is not an isomorphism it will thus suffice to show that $\pi_0 R \rightarrow \pi_0(R)^{\wedge}_{\mathfrak{m}}$ is not surjective. This is done in \cite[Section 05JA]{stacksproject}. (In fact, this $R$ is a standard example to show that completion at a non-finitely generated ideal need not yield a complete ring.)
\end{proof}

\section{Operads with good completion}
\label{sec:goodcompletion}

The aim of this section is to provide several examples of operads with good completion. These will include some of the standard `Koszul operads' in $\mathcal{C} = \mathrm{Mod}_k$ for $k$ a field of characteristic zero (Theorem \ref{thm:goodcompletionchar0}). In particular, this covers the nonunital commutative operad and the Lie operad. Also, we will show that the nonunital associative operad has good completion for \emph{any} $\mathcal{C}$, not necessarily linear over a field of characteristic zero (Proposition \ref{prop:assgoodcomletion}). It would be interesting to determine whether the $\mathbf{E}_n$-operads always have good completion; we do not attempt to address the issue here. Finally, it is straightforward to see that any truncated operad $\tau_n\mathcal{O}$ has good completion; we will conclude this section by deducing some useful consequences of this observation.

\subsection{A nilpotence criterion for good completion}

Recall that an operad $\mathcal{O}$ has good completion if for every $X \in \mathcal{C}$, the map
\begin{equation*}
\mathrm{triv}_{\mathcal{O}} X \rightarrow \varprojlim_n t_n \mathrm{triv}_{\mathcal{O}} X = \varprojlim_n B(\tau_n \mathcal{O}, \mathcal{O}, X)
\end{equation*}
is an isomorphism. Let us introduce the abbreviated notation
\begin{equation*}
K_n(d) := B(\tau_n \mathcal{O}, \mathcal{O}, \mathbf{1})(d).
\end{equation*}
Then the completion map above can, at the level of underlying objects of $\mathcal{C}$, explicitly be written as
\begin{equation*}
X \rightarrow \varprojlim_n \bigoplus_{d \geq 1} (K_n(d) \otimes X^{\otimes d})_{h\Sigma_d}.
\end{equation*}
If the limit were to commute past the sum, this map would always be an isomorphism; indeed, it is straightforward to see that $K_n(1) = \mathbf{1}$ and $K_n(d) \cong 0$ for $n > d > 1$. In particular, $\varprojlim_n (K_n(d) \otimes X^{\otimes d})_{h\Sigma_d} \cong 0$ for any fixed $d$. In general there seems to be no reason why the limit would commute past the sum, but it will if the following `uniform nilpotence' criterion is satisfied:

\begin{proposition}
\label{prop:uniformnilpotence}
Suppose that for any $n \geq 1$ there exists $N \geq 1$ such that for any $d >1$, the map $K_{n+N}(d) \rightarrow K_n(d)$ is $\Sigma_d$-equivariantly null. Then $\mathcal{O}$ has good completion.
\end{proposition}
\begin{proof}
The condition of the proposition implies that the map
\begin{equation*}
\bigoplus_{d \geq 2} (K_{n+N}(d) \otimes X^{\otimes d})_{h\Sigma_d} \rightarrow \bigoplus_{d \geq 2} (K_n(d) \otimes X^{\otimes d})_{h\Sigma_d}
\end{equation*}
is null and hence
\begin{equation*}
\varprojlim_n \bigoplus_{d \geq 2} (K_n(d) \otimes X^{\otimes d})_{h\Sigma_d} \cong 0.
\end{equation*}
\end{proof}
\begin{remark}
In ongoing joint work with Brantner, Hahn, and Yuan we show that the $p$-local spectral Lie operad has good completion by proving that it satisfies the criterion of Proposition \ref{prop:uniformnilpotence}.
\end{remark}

\subsection{Koszul operads over a field of characteristic zero}

For the time being we work in the $\infty$-category $\mathrm{Mod}_k$ of $k$-modules, with $k$ a field of characteristic zero, and we show that many of the `classical' operads have good completion. To be precise, we fix the following two hypotheses on an operad $\mathcal{O}$:
\begin{itemize}
\item[(1)] The homotopy groups $\pi_*\mathcal{O}(d)$ are concentrated in degree 0.
\item[(2)] The homotopy groups $\pi_*B\mathcal{O}(d)$ are concentrated in degree $d-1$.
\end{itemize}
Item (2) is essentially the assumption that $\mathcal{O}$ is `Koszul' in the sense of \cite{lodayvallette}. The (nonunital versions of the) operads $\mathbf{Lie}$, $\mathbf{Com}$, and $\mathbf{Ass}$ are standard examples. Indeed, there are isomorphisms
\begin{equation*}
B\mathbf{Lie} \cong s\mathbf{Com}^{\vee}, \quad B\mathbf{Com} \cong s\mathbf{Lie}^{\vee}, \quad B\mathbf{Ass} \cong s\mathbf{Ass}^{\vee}
\end{equation*}
where $s$ denotes the operadic suspension. It sends a symmetric sequence $\mathcal{O}$ to the symmetric sequence
\begin{equation*}
(s \mathcal{O})(d) = S^{\rho_d} \otimes \mathcal{O}(d), 
\end{equation*}
where $S^{\rho_d}$ is the one-point compactification of the reduced standard representation $\rho_d$ of $\Sigma_d$ (i.e., the quotient of the permutation representation on $\mathbb{R}^d$ by its diagonal). In particular, if $\mathcal{O}(d)$ is concentrated in degree 0, then $(s\mathcal{O})(d)$ is concentrated in degree $d-1$.

\begin{theorem}
\label{thm:goodcompletionchar0}
Let $k$ be a field of characteristic zero and $\mathcal{O}$ an operad in $\mathrm{Mod}_k$ satisfying assumptions (1) and (2) above. Then $\mathcal{O}$ has good completion.
\end{theorem}

\begin{remark}
A version of this result is \cite[Theorem 2.9.4]{gaitsgoryrozenblyum}. We offer a somewhat different argument (under slightly different hypotheses) here.
\end{remark}

Theorem \ref{thm:goodcompletionchar0} follows by combining Proposition \ref{prop:uniformnilpotence} with the following:

\begin{proposition}
\label{prop:Knd}
For any $n \geq 1$ and  $d > 1$, the map $K_{n+1}(d) \rightarrow K_n(d)$ is $\Sigma_d$-equivariantly null.
\end{proposition}

\begin{remark}
\label{rmk:rationalequiv}
Since we are working over a field of characteristic zero, the $\Sigma_d$-equivariance in the statement above is irrelevant: a map $M \rightarrow N$ of $\Sigma_d$-modules is null if and only if the underlying map of $k$-modules is null, since the module of equivariant maps is a summand of the module of $k$-linear maps $M \rightarrow N$.\end{remark}

We will deduce Proposition \ref{prop:Knd} from the following result, which does not need the assumption that $k$ is of characteristic zero.

\begin{proposition}
\label{prop:Knddegree}
For a field $k$ of arbitrary characteristic, the $k$-module $K_n(d)$ is concentrated in degree $d-n$.
\end{proposition}

We will prove Proposition \ref{prop:Knddegree} after the following:

\begin{proof}[Proof of Proposition \ref{prop:Knd}]
By Remark \ref{rmk:rationalequiv} it suffices to prove that the map of $k$-modules $K_{n+1}(d) \rightarrow K_n(d)$ is null. By Proposition \ref{prop:Knddegree} the modules $K_{n+1}(d)$ and $K_n(d)$ are concentrated in different degrees, from which the conclusion follows immediately.
\end{proof}

The proof of Proposition \ref{prop:Knddegree} consists of Lemmas \ref{lem:Knd1} and \ref{lem:Knd2} below.

\begin{lemma}
\label{lem:Knd1}
The module $K_n(d)$ is $(d-n)$-coconnective, i.e., concentrated in degrees $\leq d-n$.
\end{lemma}
\begin{remark}
The proof of this lemma relies only on assumption (1) above and does not need that $\mathcal{O}$ is Koszul.
\end{remark}
\begin{proof}
Write $\tau_{>n} \mathcal{O}$ for the fiber of $\mathcal{O} \rightarrow \tau_n\mathcal{O}$ and consider the fiber sequence of symmetric sequences
\begin{equation*}
B(\tau_{>n} \mathcal{O}, \mathcal{O},\mathbf{1}) \rightarrow B(\mathcal{O},\mathcal{O},\mathbf{1}) \rightarrow B(\tau_n \mathcal{O}, \mathcal{O},\mathbf{1}).
\end{equation*}
The middle term is isomorphic to the symmetric sequence $\mathbf{1}$, so that for $d>1$ the module $K_n(d)$ is the suspension of $B(\tau_{>n} \mathcal{O}, \mathcal{O},\mathbf{1})(d)$. Hence it suffices to show that the latter module is $(d-n-1)$-coconnective. Observe that the submodule of nondegenerate $p$-simplices inside
\begin{equation*}
B(\tau_{>n} \mathcal{O}, \mathcal{O},\mathbf{1})_p(d) = \tau_{>n} \mathcal{O} \circ \mathcal{O}^{\circ p}
\end{equation*}
is concentrated in arities $\geq n+p+1$. Moreover, this module is concentrated in homological degree 0 and therefore contributes to the homology of the realization $B(\tau_{>n} \mathcal{O}, \mathcal{O},\mathbf{1})$ in degree $p$. Thus in arity $d$ we only see contributions from terms with $n+p+1 \leq d$, or equivalently $p \leq d-n-1$, as desired.
\end{proof}

\begin{lemma}
\label{lem:Knd2}
The module $K_n(d)$ is $(d-n)$-connective, i.e., concentrated in degrees $\geq d-n$.
\end{lemma}
\begin{proof}
For $n=1$ we have $K_1(d) = B\mathcal{O}(d)$ and the conclusion follows from hypothesis (2) for the operad $\mathcal{O}$. For $n > 1$ we reason by induction on the fiber sequences
\begin{equation*}
B(\mathcal{O}(n), \mathcal{O}, \mathbf{1}) \rightarrow B(\tau_n \mathcal{O}, \mathcal{O}, \mathbf{1}) \rightarrow B(\tau_{n-1} \mathcal{O}, \mathcal{O}, \mathbf{1}).
\end{equation*}
In arity $d$, the inductive hypothesis states that the rightmost term has connectivity $d-(n-1) > d-n$. The leftmost term may be identified with the composition product $\mathcal{O}(n) \circ B\mathcal{O}$. In arity $d$ this is a colimit of terms of the form
\begin{equation*}
\mathcal{O}(n) \otimes B\mathcal{O}(k_1) \otimes \cdots \otimes B\mathcal{O}(k_n)
\end{equation*}
with $k_1 + \cdots + k_n = d$. By hypothesis (2) for the operad $\mathcal{O}$, the connectivity of this term is $\sum_i (k_i -1) = d-n$. We conclude that in arity $d$ the middle term is also $(d-n)$-connective.
\end{proof}

\subsection{The associative operad}

Now let us return to the case of a general $\mathcal{C}$ and prove that the nonunital associative operad $\mathbf{Ass}$ has good completion. This operad can be defined in the $\infty$-category $\mathrm{Sp}$ of spectra, but since any stable presentably symmetric monoidal $\infty$-category $\mathcal{C}$ is tensored over $\mathrm{Sp}$, the operad $\mathbf{Ass}$ can be interpreted in $\mathcal{C}$ as well.

\begin{proposition}
\label{prop:assgoodcomletion}
The nonunital associative operad has good completion.
\end{proposition}
\begin{proof}
It suffices to treat the case where $\mathcal{C} = \mathrm{Sp}$, since the general case is obtained by base change along the unique colimit-preserving symmetric monoidal functor $\mathrm{Sp} \rightarrow \mathcal{C}$. The symmetric groups act freely on the terms of $\mathbf{Ass}$. As a consequence, $\Sigma_d$ also acts freely on the spectrum  
\begin{equation*}
K_n(d) := B(\tau_n \mathbf{Ass}, \mathbf{Ass}, \mathbf{1})(d).
\end{equation*}
Hence it suffices to show that there exists $N$, independent of $d$, such that $K_{n+N}(d) \rightarrow K_n(d)$ is null as a map of spectra, without taking $\Sigma_d$-equivariance into account. Observe that each $K_n(d)$ is in particular a finite spectrum, hence has finitely generated integral homology $H_*K_n(d)$. Moreover, Proposition \ref{prop:Knddegree} implies that the homology $H_*K_n(d)$ is concentrated in degree $d-n$ and is torsionfree; indeed, if it had $p$-primary torsion for some $p$, then $H_*(K_n(d); \mathbb{F}_p)$ could not have been concentrated in a single degree. We conclude that $K_n(d)$ must be a finite direct sum of copies of the shifted sphere spectrum $\mathbb{S}^{d-n}$. Each of the maps $K_{n+1}(d) \rightarrow K_n(d)$ can therefore be written as a matrix of which the coefficients live in
\begin{equation*}
\pi_{d-n-1} \mathbb{S}^{d-n} \cong \pi_{-1} \mathbb{S}^0 \cong 0.
\end{equation*}
It follows that the inverse system $\{K_n(d)\}_{n \geq 1}$ satisfies the nilpotence hypothesis of Proposition  \ref{prop:uniformnilpotence}.
\end{proof}

\subsection{Truncated operads}
\label{sec:truncatedcompleteness}

Consider a truncated operad $\mathcal{O} = \tau_n\mathcal{O}$. Then clearly any $\mathcal{O}$-algebra $X$ is nilcomplete by construction. In particular, we have:

\begin{lemma}
\label{lem:truncatedgoodcompletion}
A truncated operad $\mathcal{O} = \tau_n\mathcal{O}$ has good completion.
\end{lemma}

This observation is trivial, but it allows us to get a much sharper formulation of Koszul duality, as well as of the relation between ind-conilpotent and general $B\mathcal{O}$-coalgebras, in this truncated case:

\begin{proposition}
\label{prop:dualitytruncatedO}
Let $\mathcal{O} = \tau_n\mathcal{O}$ be a truncated operad. Then the indecomposables functor
\begin{equation*}
\mathrm{indec}_{\mathcal{O}}^{\mathrm{nil}}\colon \mathrm{Alg}_{\mathcal{O}}(\mathcal{C}) \rightarrow \mathrm{coAlg}^{\mathrm{dp},\mathrm{nil}}_{B\mathcal{O}}(\mathcal{C})
\end{equation*}
is fully faithful. Moreover, the comparison functor
\begin{equation*}
\mathrm{coAlg}^{\mathrm{dp},\mathrm{nil}}_{B\mathcal{O}}(\mathcal{C}) \rightarrow \mathrm{coAlg}^{\mathrm{dp}}_{B\mathcal{O}}(\mathcal{C})
\end{equation*}
of Section \ref{sec:coalgebras} is fully faithful on the essential image of $\mathrm{indec}_{\mathcal{O}}^{\mathrm{nil}}$. The essential image of $\mathrm{indec}_{\mathcal{O}}$ in $\mathrm{coAlg}^{\mathrm{dp}}_{B\mathcal{O}}(\mathcal{C})$ can be characterized in the following two equivalent ways:
\begin{itemize}
\item[(1)] It consists of those $Y \in  \mathrm{coAlg}^{\mathrm{dp}}_{B\mathcal{O}}(\mathcal{C})$ that are conilcomplete.
\item[(2)] It is the full subcategory of $\mathrm{coAlg}^{\mathrm{dp}}_{B\mathcal{O}}(\mathcal{C})$ generated under colimits by trivial divided power $B\mathcal{O}$-coalgebras.
\end{itemize}
\end{proposition}
\begin{proof}
Any $\mathcal{O}$-algebra is nilcomplete and hence, by Proposition \ref{prop:goodcompletion} and Lemma \ref{lem:truncatedgoodcompletion}, also homologically complete. But this is equivalent to saying $\mathrm{indec}_{\mathcal{O}}^{\mathrm{nil}}$ is fully faithful, since homological completion is the unit of the adjoint pair resulting from that functor.

Also, combining the fact that every $\mathcal{O}$-algebra is nilcomplete with Theorem \ref{thm:completeness} guarantees that the horizontal arrow in
\[
\begin{tikzcd}
\mathrm{Alg}_{\mathcal{O}}(\mathcal{C}) \ar{rr}{\mathrm{indec}_{\mathcal{O}}} \ar{dr}[swap]{\mathrm{indec}^{\mathrm{nil}}_{\mathcal{O}}} && \mathrm{coAlg}^{\mathrm{dp}}_{B\mathcal{O}}(\mathcal{C}) \\
&\mathrm{coAlg}^{\mathrm{dp},\mathrm{nil}}_{B\mathcal{O}}(\mathcal{C}) \ar{ur} &
\end{tikzcd}
\]
is fully faithful and provides characterization (1) of its essential image. It only remains to establish description (2). Recall that $\mathrm{Alg}_{\mathcal{O}}(\mathcal{C})$ is generated under colimits by free $\mathcal{O}$-algebras. The functor $\mathrm{indec}_{\mathcal{O}}$ sends free algebras to trivial coalgebras, hence its essential image is generated by colimits of such.
\end{proof}


For general $\mathcal{O}$ it is unclear when one should expect
\begin{equation*}
\mathrm{coAlg}^{\mathrm{dp},\mathrm{nil}}_{B\mathcal{O}}(\mathcal{C}) \rightarrow \mathrm{coAlg}^{\mathrm{dp}}_{B\mathcal{O}}(\mathcal{C})
\end{equation*}
to be fully faithful (cf. Remark \ref{rmk:comparingcoalgs}). However, Proposition \ref{prop:dualitytruncatedO} suggests that it might be worth singling out the following subcategory of  $\mathrm{coAlg}^{\mathrm{dp}}_{B\mathcal{O}}(\mathcal{C})$ instead:

\begin{definition}
A coalgebra $Y \in \mathrm{coAlg}^{\mathrm{dp}}_{B\mathcal{O}}(\mathcal{C})$ is \emph{weakly conilpotent} if for every $n \geq 1$, its image in $\mathrm{coAlg}^{\mathrm{dp}}_{\varphi^n B\mathcal{O}}(\mathcal{C})$ is contained in the full subcategory described in points (1) and (2) of Proposition \ref{prop:dualitytruncatedO}. Write $\mathrm{coAlg}^{\mathrm{dp},\mathrm{wnil}}_{B\mathcal{O}}(\mathcal{C})$ for the full subcategory of $\mathrm{coAlg}^{\mathrm{dp}}_{B\mathcal{O}}(\mathcal{C})$ on the weakly conilpotent coalgebras.
\end{definition}
\begin{remark}
Note that we have tacitly made the identification $B(\tau_n\mathcal{O}) \cong \varphi^n B\mathcal{O}$. Also, despite the terminology, it is not clear that the image of a conilpotent coalgebra under the comparison functor $\mathrm{coAlg}^{\mathrm{dp},\mathrm{nil}}_{B\mathcal{O}}(\mathcal{C}) \rightarrow \mathrm{coAlg}^{\mathrm{dp}}_{B\mathcal{O}}(\mathcal{C})$ is also weakly conilpotent.
\end{remark}

Then Proposition \ref{prop:dualitytruncatedO} has the following consequence:

\begin{proposition}
\label{prop:weaklyconil}
The functors $\mathrm{indec}_{\tau_n\mathcal{O}}\colon \mathrm{Alg}_{\tau_n\mathcal{O}}(\mathcal{C}) \rightarrow \mathrm{coAlg}^{\mathrm{dp}}_{\varphi^nB\mathcal{O}}(\mathcal{C})$ induce an equivalence of $\infty$-categories
\begin{equation*}
\varprojlim_n \mathrm{Alg}_{\tau_n\mathcal{O}}(\mathcal{C}) \simeq \mathrm{coAlg}^{\mathrm{dp},\mathrm{wnil}}_{B\mathcal{O}}(\mathcal{C}).
\end{equation*}
\end{proposition}
\begin{proof}
From Proposition \ref{prop:dualitytruncatedO} we see that the functors $\mathrm{indec}_{\tau_n\mathcal{O}}$ give a fully faithful functor
\begin{equation*}
\varprojlim_n \mathrm{Alg}_{\tau_n\mathcal{O}}(\mathcal{C}) \to \varprojlim_n \mathrm{coAlg}_{\varphi^nB\mathcal{O}}^{\mathrm{dp}}(\mathcal{C}).
\end{equation*}
By Theorem \ref{thm:coalgdpdecomposition}, the inverse limit on the right is equivalent to $\mathrm{coAlg}_{B\mathcal{O}}^{\mathrm{dp}}(\mathcal{C})$ and the relevant full subcategory is, by construction, precisely that of weakly conilpotent coalgebras.
\end{proof}

We conclude this section by noting that the result of Proposition \ref{prop:dualitytruncatedO} can be dualized in the following way:

\begin{proposition}
Let $\mathcal{Q} = \tau^n\mathcal{Q}$ be a truncated cooperad. Then the adjoint pair
\[
\begin{tikzcd}
\mathcal{C} \ar[shift left]{r}{\mathrm{triv}_{\mathcal{Q}}} & \mathrm{coAlg}^{\mathrm{dp}}_{\mathcal{Q}}(\mathcal{C}) \ar[shift left]{l}{\mathrm{prim}_{\mathcal{Q}}}
\end{tikzcd}
\]
induces a fully faithful functor from $\mathrm{coAlg}^{\mathrm{dp}}_{\mathcal{Q}}(\mathcal{C})$ to the $\infty$-category of algebras for the monad on $\mathcal{C}$ described by
\begin{equation*}
X \mapsto \prod_{k=1}^\infty (C\mathcal{Q}(k) \otimes X^{\otimes k})_{h\Sigma_k}.
\end{equation*}
Furthermore, the functor
\begin{equation*}
\mathrm{prim}_{\mathcal{Q}}\colon \mathrm{coAlg}^{\mathrm{dp}}_{\mathcal{Q}}(\mathcal{C}) \rightarrow \mathrm{Alg}_{C\mathcal{Q}}(\mathcal{C})
\end{equation*}
is fully faithful. Its essential image can be described in the following two equivalent ways:
\begin{itemize}
\item[(1)] It consists of the nilcomplete $C\mathcal{Q}$-algebras.
\item[(2)] It is the smallest full subcategory of $\mathrm{Alg}_{C\mathcal{Q}}(\mathcal{C})$ closed under limits that contains the trivial $C\mathcal{Q}$-algebras.
\end{itemize}
\end{proposition}
\begin{proof}
All of this is dual to the arguments proving \ref{prop:dualitytruncatedO}, except perhaps for the explicit identification of the monad $\mathrm{prim}_{\mathcal{Q}} \circ \mathrm{triv}_{\mathcal{Q}}$. Note that $\mathrm{triv}_{\mathcal{Q}} \cong \mathrm{indec}_{C\mathcal{Q}} \circ \mathrm{free}_{C\mathcal{Q}}$, so that we may apply Theorem \ref{thm:completeness} to calculate
\begin{equation*}
\mathrm{prim}_{\mathcal{Q}}(\mathrm{triv}_{\mathcal{Q}} X) \cong \mathrm{prim}_{\mathcal{Q}}\mathrm{indec}_{C\mathcal{Q}}(\mathrm{free}_{C\mathcal{Q}} X) \cong \varprojlim_n t_n \mathrm{free}_{C\mathcal{Q}} X.
\end{equation*}
The last expression easily leads to the formula of the proposition. (We have implicitly used the equivalence $BC\mathcal{Q} \cong \mathcal{Q}$ throughout.)
\end{proof}

Similarly one can dualize Proposition \ref{prop:weaklyconil} to a statement for algebras. We leave the details to the interested reader.

\section{Examples}
\label{sec:examples}

In this section we sample some applications and special cases of Theorem \ref{thm:completeness} and relate it to existing results on the Francis--Gaitsgory conjecture.

\subsection{Connected algebras}

Imposing appropriate connectivity hypotheses often forces algebras to be complete. To be precise, suppose we have a $t$-structure on the stable $\infty$-category $\mathcal{C}$ that is compatible with its symmetric monoidal structure, meaning $X, Y \in \mathcal{C}_{\geq 0}$ implies $X \otimes Y \in \mathcal{C}_{\geq 0}$. Also, let us assume that this $t$-structure is left complete, meaning that each $X \in \mathcal{C}$ is the limit of its truncations:
\begin{equation*}
X \xrightarrow{\cong} \varprojlim_n X_{\leq n}.
\end{equation*}
Now suppose that $\mathcal{O}$ is a connective operad in $\mathcal{C}$, meaning $\mathcal{O}(n) \in \mathcal{C}_{\geq 0}$ for each $n$. The following was proved by Ching--Harper \cite{chingharper} in the specific case where $\mathcal{C} = \mathrm{Mod}_R$ for a connective commutative ring spectrum $R$:

\begin{theorem}
\label{thm:chingharper}
For $\mathcal{O}$ a connective operad, the adjoint pair $(\mathrm{indec}_{\mathcal{O}}, \mathrm{prim}_{B\mathcal{O}})$ restricts to an adjoint equivalence of $\infty$-categories
\[
\begin{tikzcd}
\mathrm{Alg}_{\mathcal{O}}(\mathcal{C}_{\geq 1}) \ar[shift left]{r}{\cong} & \mathrm{coAlg}_{B\mathcal{O}}^{\mathrm{dp}}(\mathcal{C}_{\geq 1}) \ar[shift left]{l}
\end{tikzcd}
\]
between connected $\mathcal{O}$-algebras and connected divided power $B\mathcal{O}$-coalgebras.
\end{theorem}

First of all, note that it is indeed clear that the functor $\mathrm{indec}_{\mathcal{O}}$ sends connected algebras to connected coalgebras, for example by inspecting the formula $\mathrm{indec}_{\mathcal{O}}(X) \cong B(\mathbf{1},\mathcal{O},X)$. We will see below (cf. Lemma \ref{lem:primconnected}) that $\mathrm{prim}_{B\mathcal{O}}$ preserves connectedness as well, so that indeed the adjunction between these functors restricts to connected objects as claimed above. To prove the theorem we begin by checking that the unit is an equivalence, which by Theorem \ref{thm:completeness} is the same as checking that every connected $\mathcal{O}$-algebra is nilcomplete:

\begin{lemma}
\label{lem:connectedcomplete}
For any $X \in \mathrm{Alg}_{\mathcal{O}}(\mathcal{C}_{\geq 1})$, the map $X \rightarrow \varprojlim_n t_n X$ is an isomorphism.
\end{lemma}
\begin{remark}
A version of this result is also already proved by Amabel \cite[Corollary 6.8]{amabel}.
\end{remark}
\begin{proof}
Write $r_n X$ for the fiber of $X \rightarrow t_n X$. Then $r_n$ is a functor from $\mathrm{Alg}_{\mathcal{O}}(\mathcal{C})$ to itself preserving sifted colimits. If $V \in \mathcal{C}_{\geq 1}$, then
\begin{equation*}
r_n(\mathrm{free}_{\mathcal{O}} V) \cong \bigoplus_{k > n} (\mathcal{O}(k) \otimes V^{\otimes k})_{h\Sigma_k} \in \mathcal{C}_{\geq n+1}.
\end{equation*}
Since a general connected $\mathcal{O}$-algebra $X$ is a sifted colimit of free connected $\mathcal{O}$-algebras, it follows that $r_n X \in \mathcal{C}_{\geq n+1}$ as well. We have assumed that the $t$-structure of $\mathcal{C}$ is left complete, so it follows that $\varprojlim_n r_n X \cong 0$, implying the result.
\end{proof}

To analyze the counit of the pair $(\mathrm{indec}_{\mathcal{O}}, \mathrm{prim}_{B\mathcal{O}})$ and prove Theorem \ref{thm:chingharper}, we first establish some basic properties of the primitives functor.

\begin{lemma}
\label{lem:primconnected}
For a connected divided power coalgebra $Y \in \mathrm{coAlg}_{B\mathcal{O}}^{\mathrm{dp}}(\mathcal{C}_{\geq 1})$, the primitives $\mathrm{prim}_{B\mathcal{O}}(Y)$ are connected as well. Moreover, the functor
\begin{equation*}
\mathrm{prim}_{B\mathcal{O}}\colon \mathrm{coAlg}_{B\mathcal{O}}^{\mathrm{dp}}(\mathcal{C}_{\geq 1}) \rightarrow \mathrm{Alg}_{\mathcal{O}}(\mathcal{C}_{\geq 1})
\end{equation*}
is conservative, i.e., detects isomorphisms.
\end{lemma}
\begin{proof}
As a consequence of Theorem \ref{thm:filtrations}(4), there is an equivalence
\begin{equation*}
\mathrm{prim}_{B\mathcal{O}}(Y) \cong \varprojlim_n \mathrm{prim}_{B\mathcal{O}}( \cdots \rightarrow f^n Y \rightarrow f^{n-1} Y \rightarrow \cdots)
\end{equation*}
and the associated graded of the tower on the right consists of the objects
\begin{equation*}
(\mathcal{O}(n) \otimes Y^{\otimes n})_{h\Sigma_n} = D_n^{\mathcal{O}}(Y)
\end{equation*}
for $n \geq 1$. These are $n$-connective, implying that $\mathrm{prim}_{B\mathcal{O}}(Y)$ is indeed connected. 

Now suppose that $f\colon Y \rightarrow Z$ is a map of connected divided power $B\mathcal{O}$-coalgebras such that $\mathrm{prim}_{B\mathcal{O}}(f)$ is an isomorphism. For the remainder of this proof, write $s(Y)$ for the fiber of the map
\begin{equation*}
\mathrm{prim}_{B\mathcal{O}}(Y) \rightarrow \mathrm{prim}_{B\mathcal{O}}(f^1Y) \cong \mathrm{forget}_{B\mathcal{O}}(Y)
\end{equation*}
and similarly for $Z$. By the above, $s(Y)$ can naturally be written as the limit of a tower with associated graded consisting of the terms $D_n^{\mathcal{O}}(Y)$ with $n \geq 2$. Consider the following diagram of fiber sequences:
\[
\begin{tikzcd}
s(Y) \ar{d}\ar{r} & \mathrm{prim}_{B\mathcal{O}}(Y) \ar{r}\ar{d} & \mathrm{forget}_{B\mathcal{O}}(Y) \ar{d} \\
s(Z) \ar{r} & \mathrm{prim}_{B\mathcal{O}}(Z) \ar{r} & \mathrm{forget}_{B\mathcal{O}}(Z).
\end{tikzcd}
\]
For any $n \geq 0$ we will prove the following two statements:
\begin{itemize}
\item[($a_n$)] The  map $s(Y) \rightarrow s(Z)$ induces isomorphisms on homotopy groups $\pi_k$ for $k \leq n$.
\item[($b_n$)] The map $Y \rightarrow Z$ induces isomorphisms on homotopy groups $\pi_k$ for $k \leq n$.
\end{itemize}
Note that ($b_0$) holds since $Y$ and $Z$ are connected. Moreover, if ($b_n$) holds, then the map $D_m^{\mathcal{O}}(Y) \rightarrow D_m^{\mathcal{O}}(Z)$ gives isomorphisms on $\pi_k$ for $k \leq n+m-1$. In particular, ($b_n$) implies ($a_{n+1}$). Furthermore, the long exact sequences associated with the fibrations defining $s(Y)$ and $s(Z)$ show that ($a_n$) implies ($b_n$). By induction, it follows that ($b_n$) holds for any $n$. Since we have assumed the $t$-structure on $\mathcal{C}$ to be left complete, we conclude that $Y \rightarrow Z$ is an isomorphism.
\end{proof}

\begin{proof}[Proof of Theorem \ref{thm:chingharper}]
It remains to check that the counit
\begin{equation*}
\varepsilon\colon \mathrm{indec}_{\mathcal{O}}\mathrm{prim}_{B\mathcal{O}}(Y) \rightarrow Y
\end{equation*}
is an isomorphism for any connected divided power $B\mathcal{O}$-coalgebra $Y$. By Lemma \ref{lem:primconnected} it will suffice to do this after applying the functor $\mathrm{prim}_{B\mathcal{O}}$ to this map. But then it follows from the commutative diagram
\[
\begin{tikzcd}[column sep = huge]
\mathrm{prim}_{B\mathcal{O}}(Y) \ar{r}{\eta \circ \mathrm{prim}_{B\mathcal{O}}}\ar[equal]{dr} & \mathrm{prim}_{B\mathcal{O}}\mathrm{indec}_{\mathcal{O}}\mathrm{prim}_{B\mathcal{O}}(Y) \ar{d}{\mathrm{prim}_{B\mathcal{O}} \circ \varepsilon} \\
& \mathrm{prim}_{B\mathcal{O}}(Y)
\end{tikzcd}
\]
since we already know (see Lemma \ref{lem:connectedcomplete}) that $\eta$ is an isomorphism.
\end{proof}

\begin{corollary}
    \label{cor:connectedcocomplete}
If $\mathcal{O}$ is a connective operad, then any connected divided power $B\mathcal{O}$-coalgebra $Y$ is conilcomplete.
\end{corollary}
\begin{proof}
This is one way of stating that the counit in the adjunction of Theorem \ref{thm:chingharper} is an isomorphism.
\end{proof}

Finally, let us observe the following variation on Theorem \ref{thm:chingharper}:

\begin{theorem}
    \label{thm:chingharpernil}
For $\mathcal{O}$ a connective operad, the functor 
\[
\mathrm{Alg}_{\mathcal{O}}(\mathcal{C}_{\geq 1}) \xrightarrow{\mathrm{indec}_{\mathcal{O}}^{\mathrm{nil}}} \mathrm{coAlg}_{B\mathcal{O}}^{\mathrm{nil},\mathrm{dp}}(\mathcal{C}_{\geq 1})
\]
is an equivalence of $\infty$-categories. The comparison functor
\[
\mathrm{coAlg}_{B\mathcal{O}}^{\mathrm{nil},\mathrm{dp}}(\mathcal{C}_{\geq 1}) \to \mathrm{coAlg}_{B\mathcal{O}}^{\mathrm{dp}}(\mathcal{C}_{\geq 1})
\]
between connected $\mathrm{Sym}_{B\mathcal{O}}$-coalgebras and connected divided power $B\mathcal{O}$-coalgebras is an equivalence of $\infty$-categories as well.
\end{theorem}
\begin{proof}
The second claim follows from the first by the fact that the composition of the comparison functor with $\mathrm{indec}_{\mathcal{O}}^{\mathrm{nil}}$ is the functor $\mathrm{indec}_{\mathcal{O}}$, which is an equivalence by Theorem \ref{thm:chingharper}. To establish the first claim, we first recall that the forgetful functor exhibits the $\infty$-category $\mathrm{coAlg}_{B\mathcal{O}}^{\mathrm{dp}}(\mathcal{C}_{\geq 1})$ as comonadic over $\mathcal{C}_{\geq 1}$. By Theorem \ref{thm:chingharper}, the $\infty$-category $\mathrm{Alg}_{\mathcal{O}}(\mathcal{C}_{\geq 1})$ must therefore be comonadic over $\mathcal{C}_{\geq 1}$ as well, via the functor $\mathrm{cot}_{\mathcal{O}}$. The resulting comonad is $\mathrm{Sym}_{B\mathcal{O}}$ by Proposition \ref{prop:stabcomonad} and the associated $\infty$-category of coalgebras is, by definition, $\mathrm{coAlg}_{B\mathcal{O}}^{\mathrm{nil},\mathrm{dp}}(\mathcal{C}_{\geq 1})$.
\end{proof}

\begin{remark}
It is certainly possible to improve on Theorem \ref{thm:chingharper} by allowing $\mathcal{O}$-algebras with nontrivial $\pi_0$, but one then has to impose further completeness conditions. For a concrete example, see Theorem \ref{thm:formaldef} below.
\end{remark}

\subsection{Formal deformation theory}

Koszul duality plays a central role in recent treatments of formal deformation theory. Lurie \cite{dagx} and Pridham \cite{pridham} showed that formal moduli problems over a field of characteristic zero are classified by differential graded Lie algebras. Brantner--Mathew \cite{brantnermathew} treated formal moduli problems in positive and mixed characteristic, in both the context of derived and of spectral algebraic geometry. They provided a classification now using the more sophisticated notion of a (spectral) partition Lie algebra \cite[Theorem 1.21]{brantnermathew}. 

A key technical ingredient in these classification results is a Koszul duality statement like Theorem \ref{thm:formaldef} below (which only addresses the case relevant for spectral algebraic geometry). In its statement $k$ denotes a field and $\mathrm{CAlg}^{\mathrm{aug}}_k$ is the $\infty$-category of augmented commutative $k$-algebras. An augmented commutative $k$-algebra $A$ is \emph{complete local Noetherian} if $A$ is connective, $\pi_0(A)$ is a complete local Noetherian $k$-algebra (in the classical sense), and $\pi_i(A)$ is a finitely generated $\pi_0(A)$-module for all $i > 0$. We denote by $\mathrm{CAlg}^{\mathrm{cN}}_k$ the full subcategory of $\mathrm{CAlg}^{\mathrm{aug}}_k$ on the complete local Noetherian $k$-algebras and by $\mathrm{Alg}_{\mathrm{Lie}^{\pi}_k}$ the $\infty$-category of spectral partition Lie algebras over $k$.

\begin{theorem}[Brantner--Mathew]
\label{thm:formaldef}
There is an equivalence of $\infty$-categories
\[
\mathfrak{D}\colon \mathrm{CAlg}^{\mathrm{cN}}_k \to (\mathrm{Alg}_{\mathrm{Lie}^{\pi}_k}^{\mathrm{ft}, \leq 0})^{\mathrm{op}}
\]
for which the underlying $k$-module of $\mathfrak{D}(A)$ is $\mathrm{cot}(A)^{\vee}$, the $k$-linear dual of the cotangent fiber of $A$. The right-hand side denotes the opposite of the full subcategory of $\mathrm{Alg}_{\mathrm{Lie}^{\pi}_k}$ on those $\mathfrak{g}$ that are of finite type (meaning that $\pi_i \mathfrak{g}$ is finite-dimensional for each $i$) and coconnective (i.e. $\pi_i \mathfrak{g} = 0$ for $i > 0$).
\end{theorem}

The analog of this result for simplicial commutative rings is \cite[Theorem 1.15]{brantnermathew} and the statement above is a formal consequence of \cite[Theorem 4.20, Proposition 5.31]{brantnermathew}. Our only aim here is to sketch how the methods of this paper, especially our Theorem \ref{thm:completeness}, provide another approach to obtaining results like Theorem \ref{thm:formaldef}, hopefully inspiring future work in this direction.

\begin{proof}[Sketch proof of Theorem \ref{thm:formaldef}]
As we have remarked before, the cooperad $B\mathbf{Com}$ is isomorphic to $s\mathbf{Lie}^{\vee}$, the operadic suspension of the linear dual of the Lie operad. It is fairly immediate from Brantner--Mathew's constructions that $k$-linear duality gives an equivalence between the $\infty$-category $(\mathrm{Alg}_{\mathrm{Lie}^{\pi}_k}^{\mathrm{ft}, \leq 0})^{\mathrm{op}}$ and the full subcategory of $\mathrm{coAlg}^{\mathrm{dp}}_{s\mathbf{Lie}^{\vee}}(\mathrm{Mod}_k)$ on divided power coalgebras whose underlying object is of finite type and connective. To prove the theorem it therefore suffices to show that the indecomposables functor for $\mathrm{CAlg}^{\mathrm{aug}}_k$ restricts to a fully faithful functor
\[
\mathrm{indec}\colon \mathrm{CAlg}^{\mathrm{cN}}_k \to \mathrm{coAlg}^{\mathrm{dp}}_{s\mathbf{Lie}^{\vee}}(\mathrm{Mod}_k)
\]
and that its essential image consists of the coalgebras that are of finite type and connective. By our Theorem \ref{thm:completeness}, the full faithfulness of this functor is equivalent to the statement that complete local Noetherian $k$-algebras are nilcomplete. This is the content of \cite[Proposition 5.29]{brantnermathew}. If $A$ is a connective augmented commutative $k$-algebra, then certainly its cotangent fiber $\mathrm{cot}(A)$ is also connective. If $A$ is moreover Noetherian, then its cotangent fiber is of finite type by the combination of Lurie's version of the Hilbert Basis Theorem \cite[Proposition 7.2.4.31]{higheralgebra} and \cite[Theorem 7.4.3.18]{higheralgebra}. Therefore it remains to show that any $\mathfrak{g} \in \mathrm{coAlg}^{\mathrm{dp}}_{s\mathbf{Lie}^{\vee}}(\mathrm{Mod}_k)$ that is connective and of finite type is in the essential image of the indecomposables functor restricted to $\mathrm{CAlg}^{\mathrm{cN}}_k$. If $V$ is a finite dimensional $k$-vector space and $n \geq 0$, then the completed free algebra $\widehat{\mathrm{Sym}}(V{[n]})$ is complete local Noetherian. Its indecomposables are the trivial coalgebra $\mathrm{triv}_{s\mathbf{Lie}^{\vee}}(V{[n]})$, which is therefore in the essential image we aim to identify. The coalgebra $\mathfrak{g}$ admits a filtration by `cellular approximation'; it is the colimit of a sequence of coalgebras
\[
\mathfrak{g}_{\leq 0} \to \mathfrak{g}_{\leq 1} \to \cdots \to \mathfrak{g}_{\leq n-1} \to \mathfrak{g}_{\leq n} \to \cdots
\]
where $\pi_i \mathfrak{g}_{\leq n} = \pi_i\mathfrak{g}$ for $i \leq n$ and $\pi_i \mathfrak{g}_{\leq n} = 0$ for $i > n$. Moreover, for each $n \geq 1$ there is a cofiber sequence
\[
\mathrm{triv}_{s\mathbf{Lie}^{\vee}}(\pi_n\mathfrak{g}{[n-1]}) \to \mathfrak{g}_{\leq n-1} \to \mathfrak{g}_{\leq n}.
\]
Since we have assumed that $\pi_n\mathfrak{g}$ is finite-dimensional, the first coalgebra in the sequence is in the essential image we are concerned with. A straightforward induction on $n$ now shows that each $\mathfrak{g}_{\leq n}$ and ultimately $\mathfrak{g}$ itself is in this essential image.
\end{proof}

\begin{remark}
Theorem \ref{thm:formaldef} is relevant in the context of formal deformations parametrized by $\mathbf{E}_\infty$-rings over a field $k$ (i.e., the setting of spectral algebraic geometry). A different version of formal deformation theory exists in the context of derived algebraic geometry, where the basic objects are parametrized by simplicial commutative (or animated) rings. To address this context, it would be interesting to develop versions of the main results of this paper that apply, for example, to the Koszul duality between simplicial commutative rings (or more generally derived rings) and partition Lie algebras.
\end{remark}

\subsection{Koszul duality in pronilpotent symmetric monoidal $\infty$-categories}

In their original paper, Francis--Gaitsgory proved their conjecture in the special case where the tensor product on $\mathcal{C}$ is nilpotent or pronilpotent \cite[Proposition 4.1.2]{francisgaitsgory}. In this section we will see how this case also follows from the methods we have developed here. The next section singles out some special cases that are useful in practice.

\begin{definition}
\label{def:pronilpotent}
Let $\mathcal{C}$ be a stable presentable $\infty$-category equipped with a nonunital symmetric monoidal structure for which the tensor product commutes with colimits in each variable separately. Then we say that $\mathcal{C}$ is \emph{pronilpotent} if it is equivalent to the inverse limit of a system
\[
\cdots \xrightarrow{f_{n+1}} \mathcal{C}_{n+1} \xrightarrow{f_n} \mathcal{C}_n \xrightarrow{f_{n-1}} \cdots \xrightarrow{f_1} \mathcal{C}_1
\]
satisfying the following:
\begin{itemize}
\item[(a)] Each $\mathcal{C}_n$ is a stable nonunital presentably symmetric monoidal $\infty$-category.
\item[(b)] Each of the functors $f_n$ is nonunital symmetric monoidal and preserves colimits. 
\item[(c)] The restriction of the tensor product functor on $\mathcal{C}_n$ to the subcategory $\mathrm{ker}(f_{n-1}) \times \mathcal{C}_n$ is (isomorphic to) the zero functor. For $n=1$ we require the tensor product on $\mathcal{C}_1$ to be identically zero.
\end{itemize}
\end{definition}

\begin{remark}
A straightforward inductive argument shows that in the setting of the definition above each $\mathcal{C}_n$ is $(n+1)$-nilpotent, in the sense that any $(n+1)$-fold tensor product in $\mathcal{C}_n$ is zero.
\end{remark}

\begin{remark}
Our hypotheses are slightly different from those in \cite[Definition 4.1.1]{francisgaitsgory}, in that we do not require the $f_n$ to preserve limits. This assumption is not necessary for the argument we present below.
\end{remark}

\begin{theorem}
\label{thm:pronilpotentKD}
Let $\mathcal{C}$ be pronilpotent and let $\mathcal{O}$ be an operad in $\mathcal{C}$. Then both of the functors
\[
\mathrm{Alg}_{\mathcal{O}}(\mathcal{C}) \xrightarrow{\mathrm{indec}_{\mathcal{O}}^{\mathrm{nil}}} \mathrm{coAlg}_{B\mathcal{O}}^{\mathrm{nil},\mathrm{dp}}(\mathcal{C}) \to \mathrm{coAlg}_{B\mathcal{O}}^{\mathrm{dp}}(\mathcal{C})
\]
are equivalences of $\infty$-categories, where the second one is the usual comparison between ind-conilpotent divided power coalgebras and general divided power coalgebras.
\end{theorem}
\begin{proof}
We will prove that the functor $\mathrm{indec}_{\mathcal{O}}\colon \mathrm{Alg}_{\mathcal{O}}(\mathcal{C}) \to \mathrm{coAlg}_{B\mathcal{O}}^{\mathrm{dp}}(\mathcal{C})$ is an equivalence of $\infty$-categories. After that, the statement that the comparison functor $\mathrm{coAlg}_{B\mathcal{O}}^{\mathrm{nil},\mathrm{dp}}(\mathcal{C}) \to \mathrm{coAlg}_{B\mathcal{O}}^{\mathrm{dp}}(\mathcal{C})$ is an equivalence is proved by the same comonadicity argument used in the proof of Theorem \ref{thm:chingharpernil}.

Since the functors $f_n$ preserve tensor products and colimits, they induce corresponding functors $\mathrm{Alg}_{\mathcal{O}}(\mathcal{C}) \to \mathrm{Alg}_{\mathcal{O}}(\mathcal{C}_n)$ and similarly $\mathrm{coAlg}_{B\mathcal{O}}^{\mathrm{dp}}(\mathcal{C}) \to \mathrm{coAlg}_{B\mathcal{O}}^{\mathrm{dp}}(\mathcal{C}_n)$. By Lemma \ref{lem:algcoalgpronil} below, the resulting functor $\mathrm{Alg}_{\mathcal{O}}(\mathcal{C}) \to \varprojlim_n \mathrm{Alg}_{\mathcal{O}}(\mathcal{C}_n)$ is an equivalence and the analogous statement holds for $\mathrm{coAlg}_{B\mathcal{O}}^{\mathrm{dp}}(\mathcal{C})$. Moreover, for each $n$ the diagram
\[
\begin{tikzcd}
\mathrm{Alg}_{\mathcal{O}}(\mathcal{C}) \ar{r}{\mathrm{indec}_{\mathcal{O}}}\ar{d}{f_n} & \mathrm{coAlg}_{B\mathcal{O}}^{\mathrm{dp}}(\mathcal{C}) \ar{d}{f_n} \\
\mathrm{Alg}_{\mathcal{O}}(\mathcal{C}_n) \ar{r}{\mathrm{indec}_{\mathcal{O}}} & \mathrm{coAlg}_{B\mathcal{O}}^{\mathrm{dp}}(\mathcal{C}_n) 
\end{tikzcd}
\]
commutes, again as a straightforward consequence of Definition \ref{def:pronilpotent}(b). Thus, it will suffice to prove that the functor on the bottom row is an equivalence.

The map of operads $\mathcal{O} \to \tau_n\mathcal{O}$ induces a functor $\mathrm{Alg}_{\mathcal{O}}(\mathcal{C}_n) \to \mathrm{Alg}_{\tau_n\mathcal{O}}(\mathcal{C}_n)$, which is easily seen to be an equivalence by Theorem \ref{thm:algdecomposition}, relying on the fact that $k$-fold tensor products vanish in $\mathcal{C}_n$ for $k > n$. It follows that every $X \in \mathrm{Alg}_{\mathcal{O}}(\mathcal{C}_n)$ is nilcomplete. The dual argument shows that every coalgebra $Y \in \mathrm{coAlg}_{B\mathcal{O}}^{\mathrm{dp}}(\mathcal{C}_n)$ is conilcomplete, now relying on the map of cooperads $\tau^n B\mathcal{O} \to B\mathcal{O}$ and Theorem \ref{thm:coalgdpdecomposition}. We conclude from Theorem \ref{thm:completeness} that $\mathrm{indec}_{\mathcal{O}}$ gives an equivalence between $\mathrm{Alg}_{\mathcal{O}}(\mathcal{C}_n)$ and $\mathrm{coAlg}_{B\mathcal{O}}^{\mathrm{dp}}(\mathcal{C}_n)$, as desired.
\end{proof}

\begin{lemma}
    \label{lem:algcoalgpronil}
Let $\mathcal{C} = \varprojlim_n \mathcal{C}_n$ be as in Definition \ref{def:pronilpotent} and let $\mathcal{O}$ be an operad in $\mathcal{C}$. Then the resulting functors $\mathrm{Alg}_{\mathcal{O}}(\mathcal{C}) \to \varprojlim_n \mathrm{Alg}_{\mathcal{O}}(\mathcal{C}_n)$ and $\mathrm{coAlg}^{\mathrm{dp}}_{B\mathcal{O}}(\mathcal{C}) \to \varprojlim_n \mathrm{coAlg}^{\mathrm{dp}}_{B\mathcal{O}}(\mathcal{C}_n)$ are equivalences.
\end{lemma}
\begin{proof}
We will give the proof for $\mathrm{coAlg}^{\mathrm{dp}}_{B\mathcal{O}}(\mathcal{C})$; the argument for $\mathrm{Alg}_{\mathcal{O}}(\mathcal{C})$ is entirely dual. The limit decomposition of $\mathrm{coAlg}^{\mathrm{dp}}_{B\mathcal{O}}(\mathcal{C})$ described in Theorem \ref{thm:coalgdpdecomposition} shows that it will suffice to show the following statement: for every symmetric sequence $A$, the induced functor
\[
f_\infty\colon \mathrm{coalg}_{D_k^A}(\mathcal{C}) \to \varprojlim_n \mathrm{coalg}_{D_k^A}(\mathcal{C}_n)
\]
is an equivalence of $\infty$-categories. (Recall that $D_k^A(X) = (A(k) \otimes X^{\otimes k})_{h\Sigma_k}$ and that the definition of $\mathrm{coalg}_F$ was given right after Theorem \ref{thm:coalgdpdecomposition}.) Note that the construction of the functor $f_\infty$ relies on the commutative squares 
\[
\begin{tikzcd}
\mathcal{C} \ar{d}{f_n}\ar{r}{D_k^A} & \mathcal{C} \ar{d}{f_n} \\
\mathcal{C}_n \ar{r}{D_k^A} & \mathcal{C}_n
\end{tikzcd}
\]
provided by the nonunital symmetric monoidal structure of $f_n$ and the fact that it preserves colimits. The fact that $f_\infty$ is indeed an equivalence now follows immediately from the construction of the $\infty$-categories $\mathrm{coalg}_{D_k^A}(\mathcal{C})$ as a finite limit of $\infty$-categories of the form $\mathcal{C}$ and $\mathcal{C}^{\Delta^1}$.
\end{proof}

\subsection{Koszul duality for graded (co)algebras and left (co)modules}

In this final section we sample some applications of Theorem \ref{thm:pronilpotentKD}. In particular, we show that Koszul duality works perfectly in the setting of algebras and coalgebras that are graded on the positive natural numbers as in Section \ref{sec:filteredcoalgebras} (see Theorem \ref{thm:KDgraded}) and in the setting of left modules over an operad $\mathcal{O}$ and left comodules over its bar construction $B\mathcal{O}$ (see Theorem \ref{thm:KDleftmodules}).

In the setting of algebras and coalgebras that are graded on the positive natural numbers as in Section \ref{sec:filteredcoalgebras}, Koszul duality works perfectly. Recall the $\infty$-category $\mathrm{gr}(\mathcal{C}) = \prod_{i \geq 1} \mathcal{C}$ with Day convolution as its (nonunital) symmetric monoidal structure. 

\begin{theorem}
\label{thm:KDgraded}
Let $\mathcal{O}$ be any operad in $\mathcal{C}$. Then both of the functors
\[
\mathrm{Alg}_{\mathcal{O}}(\mathrm{gr}(\mathcal{C})) \xrightarrow{\mathrm{indec}_{\mathcal{O}}^{\mathrm{nil}}} \mathrm{coAlg}_{B\mathcal{O}}^{\mathrm{nil},\mathrm{dp}}(\mathrm{gr}(\mathcal{C})) \to \mathrm{coAlg}_{B\mathcal{O}}^{\mathrm{dp}}(\mathrm{gr}(\mathcal{C}))
\]
are equivalences of $\infty$-categories.
\end{theorem}
\begin{proof}
The nonunital symmetric monoidal $\infty$-category $\mathrm{gr}(\mathcal{C})$ is easily seen to be pronilpotent by writing it as the limit over $n$ of the $\infty$-categories $\mathrm{gr}_{\leq n}(\mathcal{C}) = \prod_{i=1}^n \mathcal{C}$. The theorem is therefore a direct consequence of Theorem \ref{thm:pronilpotentKD}.
\end{proof}


We have defined (co)operads as (co)algebra objects in the monoidal $\infty$-category $\mathrm{SSeq}(\mathcal{C})$ of symmetric sequences in $\mathcal{C}$, so it makes sense to speak of left and right (co)modules over a (co)operad in $\mathrm{SSeq}(\mathcal{C})$. Recall that by convention our symmetric sequences start in arity 1 and our (co)operads are reduced (meaning $\mathcal{O}(1) = \mathbf{1} = \mathcal{Q}(1)$). The cases of left and right (co)modules behave quite differently, since the composition product preserves colimits in the first, but not in the second variable. The case of right (co)modules is the simplest:

\begin{lemma}
\label{lem:rightmodules}
For $\mathcal{P}$ an operad in $\mathcal{C}$ (or $\mathcal{Q}$ a cooperad in $\mathcal{C}$), the $\infty$-category $\mathrm{RMod}_{\mathcal{P}}(\mathrm{SSeq}(\mathcal{C}))$ of right $\mathcal{P}$-modules (resp. the $\infty$-category $\mathrm{RcoMod}_{\mathcal{Q}}(\mathrm{SSeq}(\mathcal{C}))$ of right $\mathcal{Q}$-comodules) is stable, presentable, and has an essentially unique symmetric monoidal structure for which the forgetful functor
\[
\mathrm{RMod}_{\mathcal{P}}(\mathrm{SSeq}(\mathcal{C})) \to \mathrm{SSeq}(\mathcal{C}) \quad\quad \text{(resp.} \quad \mathrm{RcoMod}_{\mathcal{Q}}(\mathrm{SSeq}(\mathcal{C})) \to \mathrm{SSeq}(\mathcal{C}) \text{)}
\]
is symmetric monoidal, where $\mathrm{SSeq}(\mathcal{C})$ is equipped with the Day convolution symmetric monoidal structure (see Remark \ref{rmk:compositionproduct}).
\end{lemma}
\begin{proof}
These statements are easily deduced from the fact that $- \circ \mathcal{P}$ (resp. $- \circ \mathcal{Q}$) is a monad (resp. a comonad) on $\mathrm{SSeq}(\mathcal{C})$ that preserves colimits and the Day convolution tensor product (see Remark \ref{rmk:compositionproduct}).
\end{proof}
\begin{remark}
Below we will discuss Koszul duality for left (co)modules. The case of right (co)modules is more straightforward. The relevant arguments are essentially those discussed in the proof of Theorem \ref{thm:barcobaroperads}; we will deal with this case in detail in \cite{brantnerheutsKD}.

\end{remark}

The $\infty$-category $\mathcal{C}$ can be considered as the full subcategory of $\mathrm{SSeq}(\mathcal{C})$ on symmetric sequences concentrated in arity 1. Thought of in this way, $\mathcal{C}$ acts on $\mathrm{SSeq}(\mathcal{C})$ via the composition product:
\[
\mathcal{C} \to \mathrm{End}(\mathrm{SSeq}(\mathcal{C}))\colon X \mapsto X \circ -.
\]
Alternatively, this action is simply the levelwise tensor product of a symmetric sequence with the object $X$. Now consider an operad $\mathcal{O}$ in $\mathcal{C}$. Through the action of $\mathcal{C}$ on $\mathrm{SSeq}(\mathcal{C})$ it makes sense to speak of $\mathcal{O}$-algebras in $\mathrm{SSeq}(\mathcal{C})$. However, these turn out to be exactly the same thing as left $\mathcal{O}$-modules. Indeed, first observe that $\mathcal{O} = \mathrm{Sym}_{\mathcal{O}}(\mathbf{1})$, where $\mathbf{1}$ is the unit symmetric sequence. Since the composition product preserves colimits and is symmetric monoidal in its left variable, it follows that there is a natural isomorphism
\[
\mathcal{O} \circ A \cong \mathrm{Sym}_{\mathcal{O}}(\mathbf{1}) \circ A \cong \mathrm{Sym}_{\mathcal{O}}(\mathbf{1} \circ A) \cong \mathrm{Sym}_{\mathcal{O}}(A)
\]
for $A \in \mathrm{SSeq}(\mathcal{C})$. More generally, if $\mathcal{M}$ denotes either the symmetric monoidal $\infty$-category $\mathrm{RMod}_{\mathcal{P}}(\mathrm{SSeq}(\mathcal{C}))$ or $\mathrm{RcoMod}_{\mathcal{Q}}(\mathrm{SSeq}(\mathcal{C}))$ of the Lemma \ref{lem:rightmodules}, we have
\[
\mathcal{O} \circ A \cong \mathrm{Sym}_{\mathcal{O}}(A)
\]
for $A \in \mathcal{M}$ and in this setting we can identify $\mathcal{O}$-algebras and left $\mathcal{O}$-modules. Dual comments apply to left $B\mathcal{O}$-comodules and $\mathrm{Sym}_{B\mathcal{O}}$-coalgebras.



\begin{theorem}
\label{thm:KDleftmodules}
Let $\mathcal{M}$ denote either the $\infty$-category $\mathrm{RMod}_{\mathcal{P}}(\mathrm{SSeq}(\mathcal{C}))$ of right modules for an operad $\mathcal{P}$ or the $\infty$-category $\mathrm{RcoMod}_{\mathcal{Q}}(\mathrm{SSeq}(\mathcal{C}))$ of right comodules for a cooperad $\mathcal{Q}$. Let $\mathcal{O}$ be an operad in $\mathcal{C}$. Then the adjoint functors
\[
\begin{tikzcd}[column sep = large]
\mathrm{LMod}_{\mathcal{O}}(\mathcal{M}) \ar[shift left]{r}{B(\mathbf{1}, \mathcal{O},-)} & \mathrm{LcoMod}_{B\mathcal{O}}(\mathcal{M}) \ar[shift left]{l}{C(\mathbf{1}, B\mathcal{O},-)} 
\end{tikzcd}
\]
are mutually inverse equivalences of $\infty$-categories.
\end{theorem}
\begin{proof}
By the discussion above, we may identify the adjoint pair of the theorem with
\[
\begin{tikzcd}
\mathrm{Alg}_{\mathcal{O}}(\mathcal{M}) \ar[shift left]{r}{\mathrm{indec}^{\mathrm{nil}}_{\mathcal{O}}} & \mathrm{coAlg}^{\mathrm{nil},\mathrm{dp}}_{B\mathcal{O}}(\mathcal{M}). \ar[shift left]{l}{\mathrm{prim}^{\mathrm{nil}}_{B\mathcal{O}}} 
\end{tikzcd}
\]
The symmetric monoidal $\infty$-category $\mathcal{M}$ is pronilpotent. Indeed, it can be written as the inverse limit over $n$ of the $\infty$-categories of left $\mathcal{P}$-modules (or left $\mathcal{Q}$-comodules) in the $\infty$-category of $n$-truncated symmetric sequences $\mathrm{SSeq}_{\leq n}(\mathcal{C})$. Therefore Theorem \ref{thm:pronilpotentKD} applies and the result follows.
\end{proof}

\appendix

\section{Coalgebras}
\label{app:coalgebras}

Let $\mathcal{Q}$ be a cooperad in a presentably symmetric monoidal stable $\infty$-category $\mathcal{C}$. The aim of this appendix is to give the precise construction of the $\infty$-category $\mathrm{coAlg}^{\mathrm{dp}}_{\mathcal{Q}}(\mathcal{C})$ of divided power $\mathcal{Q}$-coalgebras (Definition \ref{def:dpQcoalgebras}) and establish some of its fundamental properties. In particular, we will prove Theorem \ref{thm:coalgdpdecomposition} (the inductive decomposition of the $\infty$-category of divided power coalgebras) and Proposition \ref{prop:cofreetruncated} (the description of cofree divided power coalgebras for truncated cooperads). 


As we explained at the start of Section \ref{sec:coalgebras}, the main difficulty in describing divided power coalgebras is that the functor $\widehat{\mathrm{Sym}}_{\mathcal{Q}} = \prod_{n \geq 1} D_n^{\mathcal{Q}}$ is generally not a comonad on $\mathcal{C}$. There are (at least) two ways to address this issue:
\begin{itemize}
\item[(1)] One can embrace the fact that the functor $\widehat{\mathrm{Sym}}\colon \mathrm{SSeq}(\mathcal{C}) \to \mathrm{End}(\mathcal{C})$ is only lax monoidal. It is possible to define a notion of \emph{lax comonad}, of which $\widehat{\mathrm{Sym}}_{\mathcal{Q}}$ is an example, and a notion of coalgebra for such a lax comonad. This approach can be made precise by adapting the work of Anel \cite{anel} to the $\infty$-categorical setting.
\item[(2)] Write $\mathrm{Sym}^{\leq n}_A := \bigoplus_{k \leq n} D_n^A$, so that $\widehat{\mathrm{Sym}} \cong \varprojlim_n \mathrm{Sym}^{\leq n}$. The failure of $\widehat{\mathrm{Sym}}$ to be monoidal essentially boils down to the failure of the natural transformations 
\[
\widehat{\mathrm{Sym}}_A \circ \varprojlim_n \mathrm{Sym}^{\leq n}_B \to \varprojlim_n \widehat{\mathrm{Sym}}_A \circ \mathrm{Sym}^{\leq n}_B
\]
to be isomorphisms. This problem can be circumvented in a rather formal way, namely by redefining $\widehat{\mathrm{Sym}}$ to act on the $\infty$-category $\mathrm{Pro}(\mathcal{C})$ of pro-objects in $\mathcal{C}$.
\end{itemize}

In this appendix we will work out the second approach in detail. We recall some generalities on Pro-objects in Section \ref{sec:pro}, describing our main constructions and their properties in \ref{sec:appdpcoalgs}, and proving the necessary technical lemmas in \ref{subsec:limitscofree}.

\subsection{Pro-objects and completions of cofree coalgebras}
\label{sec:pro}

For the basics of pro-objects we will mostly rely on \cite[Section A.8.1]{SAG}. The $\infty$-category $\mathrm{Pro}(\mathcal{C})$ is the full subcategory of $\mathrm{Fun}(\mathcal{C},\mathcal{S})^{\mathrm{op}}$ on the functors that are accessible and preserve finite limits. The Yoneda embedding determines a functor we denote
\[
c\colon \mathcal{C} \to \mathrm{Pro}(\mathcal{C})\colon X \mapsto \mathrm{Map}_{\mathcal{C}}(X,-)
\]
with the $c$ standing for \emph{constant} pro-objects. This functor is fully faithful and preserves small colimits \cite[Proposition 5.1.3.2]{HTT} as well as finite limits \cite[Proposition 5.3.5.14]{HTT}. By \cite[Remark A.8.1.5]{SAG} a general object of $\mathrm{Pro}(\mathcal{C})$ may always be written as the limit of a filtered diagram $\{cX_{\alpha}\}$ of constant pro-objects. The functor $c$ admits a right adjoint $M$ (sometimes called \emph{materialization}) which takes such a `formal' filtered limit in $\mathrm{Pro}(\mathcal{C})$ to the actual limit $\varprojlim_{\alpha} X_{\alpha}$ in $\mathcal{C}$.

The $\infty$-category $\mathrm{Pro}(\mathcal{C})$ freely adjoins all small filtered limits to $\mathcal{C}$. More precisely, it is characterized by the following universal property. Let $\mathcal{D}$ be an $\infty$-category that admits small filtered limits. Then precomposing by $c$ defines an equivalence of $\infty$-categories
\[
c^*\colon \mathrm{Fun}^{\mathrm{filt}}(\mathrm{Pro}(\mathcal{C}), \mathcal{D}) \xrightarrow{\simeq} \mathrm{Fun}(\mathcal{C},\mathcal{D}),
\]
where the domain denotes the full subcategory of $\mathrm{Fun}(\mathrm{Pro}(\mathcal{C}), \mathcal{D})$ on the functors preserving small filtered limits \cite[Proposition A.8.1.6]{SAG}. In particular, setting $\mathcal{D} = \mathrm{Pro}(\mathcal{C})$ we find an equivalence of $\infty$-categories
\[
\mathrm{End}^{\mathrm{filt}}(\mathrm{Pro}(\mathcal{C})) \xrightarrow{\simeq} \mathrm{Fun}(\mathcal{C},\mathrm{Pro}(\mathcal{C})).
\]
Every functor $F \in \mathrm{End}(\mathcal{C})$ gives a functor $c \circ F\colon \mathcal{C} \to \mathrm{Pro}(\mathcal{C})$, which by the above equivalence defines a corresponding functor $pF \in \mathrm{End}^{\mathrm{filt}}(\mathrm{Pro}(\mathcal{C}))$ which we refer to as the \emph{prolongation} of $F$. This process of prolongation is easily seen to define a fully faithful monoidal functor
\[
\mathrm{End}(\mathcal{C}) \to \mathrm{End}^{\mathrm{filt}}(\mathrm{Pro}(\mathcal{C}))\colon F \mapsto pF,
\]
whose essential image is those endofunctors in $\mathrm{End}^{\mathrm{filt}}(\mathrm{Pro}(\mathcal{C}))$ that preserve constant pro-objects.

\begin{lemma}
\label{lem:prolongationadjoint}
The prolongation $p$ admits a right adjoint described by
\[
\mathrm{End}^{\mathrm{filt}}(\mathrm{Pro}(\mathcal{C})) \to \mathrm{End}(\mathcal{C})\colon G \to M \circ G \circ c.
\]
\end{lemma}
\begin{proof}
It suffices to observe that the functor
\[
c \circ - \colon \mathrm{End}(\mathcal{C}) \to \mathrm{Fun}(\mathcal{C}, \mathrm{Pro}(\mathcal{C})) 
\]
has right adjoint $M \circ -$.
\end{proof}


Now let $n \geq 1$ and consider the functor
\[
\mathrm{Sym}^{\leq n}\colon \mathrm{SSeq}(\mathcal{C}) \to \mathrm{End}(\mathcal{C})\colon A \mapsto \bigoplus_{k \leq n} D_k^A.
\]
It can be thought of as the composite of the lax monoidal $n$-truncation $\tau_n$ and the monoidal functor $\mathrm{Sym}$, so that $\mathrm{Sym}^{\leq n}$ is also naturally a lax monoidal functor. Since $p$ is monoidal, we find a lax monoidal functor
\[
p\mathrm{Sym}^{\leq n}\colon \mathrm{SSeq}(\mathcal{C}) \to \mathrm{End}^{\mathrm{filt}}(\mathrm{Pro}(\mathcal{C})). 
\]
Limits in the $\infty$-category of lax monoidal functors are computed on the `underlying' functors, so that 
\[
\widehat{p\mathrm{Sym}} := \varprojlim_n p\mathrm{Sym}^{\leq n}
\]
also acquires a lax monoidal structure.

\begin{remark}
\label{rmk:pSym}
Generally it is \emph{not} the case that for a symmetric sequence $A$, the functor $\widehat{p\mathrm{Sym}}_A$ is the prolongation of $\widehat{\mathrm{Sym}}_A$. However, for $X \in \mathcal{C}$ we do have natural isomorphisms $\widehat{\mathrm{Sym}}_A(X) \cong M\widehat{p\mathrm{Sym}}_A(cX)$.
\end{remark}

\begin{lemma}
\label{lem:pSymmon}
The functor $\widehat{p\mathrm{Sym}}\colon \mathrm{SSeq}(\mathcal{C}) \to \mathrm{End}^{\mathrm{filt}}(\mathrm{Pro}(\mathcal{C}))$ is monoidal.
\end{lemma}
\begin{proof}
We should check that for symmetric sequences $A, B \in \mathrm{SSeq}(\mathcal{C})$, the natural transformation
\[
\widehat{p\mathrm{Sym}}_A \circ \widehat{p\mathrm{Sym}}_B \to \widehat{p\mathrm{Sym}}_{A \circ B}
\]
arising from the lax monoidal structure of $\widehat{p\mathrm{Sym}}$ is an isomorphism. This map can be identified with the composite of the following natural isomorphisms:
\begin{eqnarray*}
\varprojlim_m p\mathrm{Sym}^{\leq m}_A \circ \varprojlim_n p\mathrm{Sym}^{\leq n}_B & \cong & \varprojlim_{m,n}p(\mathrm{Sym}^{\leq m}_A \circ \mathrm{Sym}^{\leq n}_B) \\
& \cong & \varprojlim_k p\mathrm{Sym}^{\leq k}_{A \circ B}.
\end{eqnarray*}
The first isomorphism relies on the fact that the functors $p\mathrm{Sym}^{\leq m}_A$ preserve filtered limits (by construction) and that $p$ is monoidal, the second isomorphism on the fact that the functor $\mathrm{Sym}\colon \mathrm{SSeq}(\mathcal{C}) \to \mathrm{End}(\mathcal{C})$ is monoidal.
\end{proof}

\subsection{Divided power coalgebras}
\label{sec:appdpcoalgs}

Let $\mathcal{Q}$ be a cooperad in $\mathcal{C}$. According to Lemma \ref{lem:pSymmon} it defines a corresponding comonad $\widehat{p\mathrm{Sym}}_{\mathcal{Q}}$ on $\mathrm{Pro}(\mathcal{C})$. 

\begin{definition}
\label{def:dpQcoalgebras}
We define the $\infty$-category of \emph{divided power $\mathcal{Q}$-coalgebras} to be the pullback in the following square:
\[
\begin{tikzcd}
\mathrm{coAlg}^{\mathrm{dp}}_{\mathcal{Q}}(\mathcal{C}) \ar{r}\ar{d} & \mathrm{coAlg}_{\widehat{p\mathrm{Sym}}_{\mathcal{Q}}}(\mathrm{Pro}(\mathcal{C})) \ar{d}{\mathrm{forget}} \\
\mathcal{C} \ar{r}{c} & \mathrm{Pro}(\mathcal{C}).
\end{tikzcd}
\]
\end{definition}

The comonad $\mathrm{Sym}_{\mathcal{Q}}$ gives a `prolonged' comonad $p\mathrm{Sym}_{\mathcal{Q}}$ on $\mathrm{Pro}(\mathcal{C})$ and by construction we have a natural map of comonads $p\mathrm{Sym}_{\mathcal{Q}} \to \widehat{p\mathrm{Sym}}_{\mathcal{Q}}$. These observations give rise to the following commutative diagram:
\[
\begin{tikzcd}
\mathrm{coAlg}_{\mathrm{Sym}_{\mathcal{Q}}}(\mathcal{C}) \ar{r}\ar{d}{\mathrm{forget}} & \mathrm{coAlg}_{p\mathrm{Sym}_{\mathcal{Q}}}(\mathrm{Pro}(\mathcal{C})) \ar{d}{\mathrm{forget}}\ar{r} & \mathrm{coAlg}_{\widehat{p\mathrm{Sym}}_{\mathcal{Q}}}(\mathrm{Pro}(\mathcal{C})) \ar{d}{\mathrm{forget}} \\
\mathcal{C} \ar{r}{c} & \mathrm{Pro}(\mathcal{C}) \ar[equal]{r} & \mathrm{Pro}(\mathcal{C}).
\end{tikzcd}
\]
Recall that in the body of this paper we have written $\mathrm{coAlg}^{\mathrm{nil},\mathrm{dp}}_{\mathcal{Q}}(\mathcal{C})$ for the $\infty$-category of $\mathrm{Sym}_{\mathcal{Q}}$-coalgebras. Hence the diagram above and the defining property of the $\infty$-category of divided power $\mathcal{Q}$-coalgebras now give a comparison functor
\[
\mathrm{coAlg}^{\mathrm{nil},\mathrm{dp}}_{\mathcal{Q}}(\mathcal{C}) \to \mathrm{coAlg}^{\mathrm{dp}}_{\mathcal{Q}}(\mathcal{C}).
\]

For the reader's convenience let us repeat the statement of Proposition \ref{prop:cofreetruncated} and then prove it:

\begin{proposition}
\label{cor:cofreetruncated}
If $\mathcal{Q} = \tau^n\mathcal{Q}$ is a truncated cooperad, then the comparison functor
\[
\mathrm{coAlg}^{\mathrm{nil,dp}}_{\mathcal{Q}}(\mathcal{C}) \to \mathrm{coAlg}^{\mathrm{dp}}_{\mathcal{Q}}(\mathcal{C})
\]
is an equivalence of $\infty$-categories.
\end{proposition}
\begin{proof}
In the truncated case, we have equivalences of comonads on $\mathrm{Pro}(\mathcal{C})$ as follows:
\[
\widehat{p\mathrm{Sym}_{\mathcal{Q}}} \cong \prod_{k \leq n} pD_k^{\mathcal{Q}} \cong p\mathrm{Sym}_{\mathcal{Q}}.
\]
In other words, the comonad $\widehat{p\mathrm{Sym}_{\mathcal{Q}}}$ is the prolongation of $\mathrm{Sym}_{\mathcal{Q}}$ and in particular preserves the full subcategory $\mathcal{C} \subseteq \mathrm{Pro}(\mathcal{C})$ of pro-constant objects. Hence the functor
\[
\mathrm{coAlg}_{\mathrm{Sym}_{\mathcal{Q}}}(\mathcal{C}) \to \mathcal{C} \times_{\mathrm{Pro}(\mathcal{C})}\mathrm{coAlg}_{\widehat{p\mathrm{Sym}_{\mathcal{Q}}}}(\mathrm{Pro}(\mathcal{C}))
\]
is an equivalence. Up to change of notation, this is the statement of the proposition.
\end{proof}

The aim of the remainder of this appendix is to prove Theorem \ref{thm:coalgdpdecomposition}. Let us recall the statement for the reader's convenience:

\begin{theorem}[Theorem \ref{thm:coalgdpdecomposition}]
For each $n \geq 2$, the evident commutative square of $\infty$-categories
\[
\begin{tikzcd}
\mathrm{coAlg}_{\varphi^n\mathcal{Q}}^{\mathrm{dp}}(\mathcal{C}) \ar{r}\ar{d} & \mathrm{coalg}_{D_n^{\mathcal{Q}}}(\mathcal{C}) \ar{d} \\
\mathrm{coAlg}_{\varphi^{n-1}\mathcal{Q}}^{\mathrm{dp}}(\mathcal{C}) \ar{r} & \mathrm{coalg}_{D_n^{\varphi^{n-1}\mathcal{Q}}}(\mathcal{C})
\end{tikzcd}
\]
is a pullback. Furthermore, the natural map
\begin{equation*}
\mathrm{coAlg}_{\mathcal{Q}}^{\mathrm{dp}}(\mathcal{C}) \rightarrow \varprojlim_n \mathrm{coAlg}_{\varphi^n\mathcal{Q}}^{\mathrm{dp}}(\mathcal{C})
\end{equation*}
is an equivalence of $\infty$-categories.
\end{theorem}

We split the proof into two lemmas, which we will deal with in Section \ref{subsec:limitscofree}:

\begin{lemma}
\label{lem:laxcoalgpb}
For every $n \geq 1$, the square
\[
\begin{tikzcd}
\mathrm{coAlg}^{\mathrm{dp}}_{\varphi^n\mathcal{Q}}(\mathcal{C}) \ar{r}\ar{d} & \mathrm{coAlg}^{\mathrm{dp}}_{\mathrm{cofree}(\mathcal{Q}(n))}(\mathcal{C}) \ar{d} \\
\mathrm{coAlg}^{\mathrm{dp}}_{\varphi^{n-1}\mathcal{Q}}(\mathcal{C}) \ar{r} & \mathrm{coAlg}^{\mathrm{dp}}_{\mathrm{cofree}(\varphi^{n-1}\mathcal{Q}(n))}(\mathcal{C})
\end{tikzcd}
\]
induced by the corresponding square of cooperads is a pullback of $\infty$-categories. Moreover, the evident functor $\mathrm{coAlg}^{\mathrm{dp}}_{\mathcal{Q}}(\mathcal{C}) \to \varprojlim_n \mathrm{coAlg}^{\mathrm{dp}}_{\varphi^n\mathcal{Q}}(\mathcal{C})$ is an equivalence.
\end{lemma}

\begin{lemma}
\label{lem:laxcoalgcofree}
For a symmetric sequence $A \in \mathrm{SSeq}(\mathcal{C})$ with $A(1) = 0$, generating a cofree cooperad $\mathcal{Q} = \mathrm{cofree}(A)$, the evident forgetful functor
\[
\mathrm{coAlg}^{\mathrm{dp}}_{\mathcal{Q}}(\mathcal{C}) \to \mathrm{coalg}_{\widehat{\mathrm{Sym}}_A}(\mathcal{C})
\]
is an equivalence of $\infty$-categories.
\end{lemma}

\subsection{Limits and cofree constructions}
\label{subsec:limitscofree}

In this section we investigate the extent to which the functor
\[
\mathrm{coOp}(\mathcal{C}) \to (\mathrm{Cat}_\infty)_{/\mathcal{C}}\colon \mathcal{Q} \mapsto (\mathrm{coAlg}^{\mathrm{dp}}_{\mathcal{Q}}(\mathcal{C}) \xrightarrow{\mathrm{forget}} \mathcal{C})
\]
preserves limits. This will lead to proofs of Lemmas \ref{lem:laxcoalgpb} and \ref{lem:laxcoalgcofree}. The following class of limits plays a special role in our arguments:

\begin{definition}
We say that a diagram $\cdots \to A_3 \to A_2 \to A_1$ in the $\infty$-category $\mathrm{SSeq}(\mathcal{C})$ is \emph{eventually constant} if for each $n \geq 1$ there exists a natural number $N(n)$ such that the maps $A_{k+1}(n) \to A_k(n)$ are isomorphisms for all $k \geq N(n)$.
\end{definition}

Note that for an eventually constant diagram $\{A_i\}_{i \geq 1}$ it need \emph{not} be true that the natural map $\mathrm{Sym}_{\varprojlim_i A_i} \to \varprojlim_i \mathrm{Sym}_{A_i}$ is an isomorphism, since the limit need not commute past the direct sum. However, this problem does not occur for the `completed' versions of $\mathrm{Sym}$ we are considering:

\begin{lemma}
\label{lem:pSymlimits}
The functor 
\[
\mathrm{SSeq}(\mathcal{C}) \to \mathrm{End}(\mathrm{Pro}(\mathcal{C}))\colon A \mapsto \widehat{p\mathrm{Sym}}_A
\]
preserves finite limits and eventually constant limits.
\end{lemma}
\begin{proof}
For finite limits this is clear from the fact that both $A \mapsto D_n^A$ and $p$ preserve finite limits. For an eventually constant diagram $\{A_i\}_{i \geq 1}$ of symmetric sequences, observe that
\[
pD_n^{\varprojlim_i A_i} \cong \varprojlim_i pD_n^{A_i}
\]
for each $n$, since the diagram $\{A_i(n)\}_{i \geq 1}$ has a constant coinitial subdiagram. Hence
\[
\widehat{p\mathrm{Sym}}_{\varprojlim_i A_i} \cong \prod_n \varprojlim_i pD_n^{A_i} \cong \varprojlim_i \widehat{p\mathrm{Sym}}_{A_i}
\]
as desired.
\end{proof}

A further useful property is the following:

\begin{lemma}
\label{lem:compositionevtconstdiag}
The composition product on $\mathrm{SSeq}(\mathcal{C})$ preserves limits of eventually constant diagrams in each variable separately.
\end{lemma}
\begin{proof}
A composition product $A \circ B$ can be written as the inverse limit of its truncations, $A \circ B = \varprojlim_n \tau_n(A \circ B)$. In turn, the $n$-truncation $\tau_n(A \circ B)$ only depends on the $n$-truncations of $A$ and $B$; more precisely, the truncation maps $A \to \tau_n A$ and $B \to \tau_n B$ induce an isomorphism
\[
\tau_n(A \circ B) \xrightarrow{\cong} \tau_n(\tau_n A \circ \tau_n B).
\]
If $\{A_i\}_{i \geq 1}$ is an eventually constant diagram, then $\{\tau_n A_i\}_{i \geq 1}$ has a constant coinitial subdiagram. Hence we find
\[
(\varprojlim_i A) \circ B \cong \varprojlim_n \tau_n(\tau_n(\varprojlim_i A_i) \circ B) \cong \varprojlim_{i,n} \tau_n(A_i \circ B) \cong \varprojlim_i (A_i \circ B)
\]
as desired. A similar argument applies for an eventually constant limit in the second variable $B$.
\end{proof}


We will now briefly review the construction of cofree objects, both in the context of cooperads and of comonads, following the discussion in \cite[Appendix B]{BCN}. Let $\mathcal{M}$ be a monoidal $\infty$-category with countable sequential limits (preserved by the tensor product in both variables) and finite products (preserved by the tensor product in the first variable). For an object $X$ of $\mathcal{M}$, define a sequence of objects $T^{(n)}(X)$ by 
\[
T^{(0)}(X) := \mathbf{1}, \quad\quad T^{(n+1)}(X) := \mathbf{1} \times (X \otimes T^{(n)}(X)).
\]
The projection $\mathbf{1} \times X \to \mathbf{1}$ defines, by recursion, a natural map $T^{(n+1)}(X) \to T^{(n)}(X)$. Then the dual of \cite[Theorem B.2]{BCN} shows that the cofree coalgebra on $X$ exists and is given by the inverse limit $\varprojlim_n T^{(n)}(X)$. Taking $\mathcal{M}$ to be the $\infty$-category $\mathrm{End}^{\mathrm{filt}}(\mathrm{Pro}(\mathcal{C}))$, this construction now produces the cofree comonad on a functor preserving filtered limits. If we take $\mathcal{M}$ to be $\mathrm{SSeq}(\mathcal{C})$, the reasoning above does not apply verbatim, since the composition product need not preserve countable sequential limits in each variable. However, an inspection of the argument of \cite{BCN} shows that one only needs to assume that the specific inverse limit $\varprojlim_n T^{(n)}(X)$ is preserved by the composition product in each variable. If $X$ is a symmetric sequence with $X(1) = 0$, then one easily verifies that the map $T^{(n+1)}(X) \to T^{(n)}(X)$ is an isomorphism in arities $\leq n+1$. Hence the diagram of symmetric sequences $\{T^{(n)}(X)\}_{n \geq 1}$ is eventually constant, so that Lemma \ref{lem:compositionevtconstdiag} implies that its limit commutes with the composition product in either variable.

\begin{proof}[Proof of Lemma \ref{lem:laxcoalgcofree}]
We have seen that the functor $\widehat{p\mathrm{Sym}}\colon \mathrm{SSeq}(\mathcal{C}) \to \mathrm{End}^{\mathrm{filt}}(\mathrm{Pro}(\mathcal{C}))$ is monoidal and preserves products. Hence it commutes with the construction $T^{(n)}$ described above; more precisely, for a given symmetric sequence $A$ we find natural isomorphisms
\[
\widehat{p\mathrm{Sym}}_{T^{(n)}(A)} \cong T^{(n)}(\widehat{p\mathrm{Sym}}_{A}).
\]
If we assume that $A(1) = 0$, then the remarks above the proof together with Lemma \ref{lem:pSymlimits} show that moreover $\widehat{p\mathrm{Sym}}$ commutes with the inverse limit over $n$, giving an isomorphism
\[
\widehat{p\mathrm{Sym}}_{\mathrm{cofree}(A)} \cong \mathrm{cofree}(\widehat{p\mathrm{Sym}}_{A}).
\]
In words, $\widehat{p\mathrm{Sym}}$ sends the cofree cooperad on $A$ to the cofree comonad on $\widehat{p\mathrm{Sym}}_A$. Now coalgebras for a comonad that is cofree on a functor $F$ are the same thing as coalgebras for the functor $F$. More precisely, we have an equivalence of $\infty$-categories
\[
\mathrm{coAlg}_{\mathrm{cofree}(\widehat{p\mathrm{Sym}}_A)}(\mathrm{Pro}(\mathcal{C})) \simeq \mathrm{coalg}_{\widehat{p\mathrm{Sym}}_A}(\mathrm{Pro}(\mathcal{C})), 
\]
which leads to
\begin{eqnarray*}
\mathrm{coAlg}^{\mathrm{dp}}_{\mathrm{cofree}(A)}(\mathcal{C}) & \simeq & \mathcal{C} \times_{\mathrm{Pro}(\mathcal{C})} \mathrm{coalg}_{\widehat{p\mathrm{Sym}}_A}(\mathrm{Pro}(\mathcal{C})) \\
& \simeq &  \mathrm{coalg}_{\widehat{\mathrm{Sym}}_A}(\mathcal{C}).
\end{eqnarray*}
The second equivalence is a straightforward consequence of Remark \ref{rmk:pSym}.
\end{proof}

Before we proceed we need some further observations on the commutation of the functor $\widehat{p\mathrm{Sym}}$ with limits. We write $\mathrm{coMnd}^{\mathrm{filt}}(\mathrm{Pro}(\mathcal{C}))$ for the $\infty$-category of comonads on $\mathrm{Pro}(\mathcal{C})$ of which the underlying functor commutes with filtered limits. As a consequence of Lemma \ref{lem:pSymmon}, the functor $\widehat{p\mathrm{Sym}}\colon \mathrm{SSeq}(\mathcal{C}) \to \mathrm{End}^{\mathrm{filt}}(\mathrm{Pro}(\mathcal{C}))$ also defines a functor from cooperads to comonads, which we use in the following statement:

\begin{lemma}
\label{lem:pSymsplitTot}
The functor $\widehat{p\mathrm{Sym}}\colon \mathrm{coOp}(\mathcal{C}) \to \mathrm{coMnd}^{\mathrm{filt}}(\mathrm{Pro}(\mathcal{C}))$ preserves $\mathrm{forget}$-split totalizations. More precisely, if $X^{\bullet}$ is a coaugmented cosimplicial object of $\mathrm{coOp}(\mathcal{C})$ such that the underlying coaugmented cosimplicial object $\mathrm{forget}(X^{\bullet})$ of $\mathrm{SSeq}(\mathcal{C})$ is split (cf. \cite[Definition 4.7.2.2]{higheralgebra}), then the map
\[
\widehat{p\mathrm{Sym}}_{X^{-1}} \to \mathrm{Tot}(\widehat{p\mathrm{Sym}}_{X^{\bullet}})
\]
is an isomorphism in $\mathrm{coMnd}^{\mathrm{filt}}(\mathrm{Pro}(\mathcal{C}))$.
\end{lemma}
\begin{proof}
The coaugmented cosimplicial object $\widehat{p\mathrm{Sym}}_{\mathrm{forget}(X^{\bullet})}$ of $\mathrm{End}^{\mathrm{filt}}(\mathrm{Pro}(\mathcal{C}))$ is split. This object underlies the diagram $\widehat{p\mathrm{Sym}}_{X^{\bullet}}$ in $\mathrm{coMnd}^{\mathrm{filt}}(\mathrm{Pro}(\mathcal{C}))$. Since the forgetful functor exhibits $\mathrm{coMnd}^{\mathrm{filt}}(\mathrm{Pro}(\mathcal{C}))$ as comonadic over $\mathrm{End}^{\mathrm{filt}}(\mathrm{Pro}(\mathcal{C}))$, the conclusion follows from the Barr--Beck--Lurie theorem \cite[Theorem 4.7.3.5]{higheralgebra}.
\end{proof}

Our final preparation for the proof of Lemma \ref{lem:laxcoalgpb} is the following lift of the statement of Lemma \ref{lem:pSymlimits} from endofunctors to comonads:

\begin{lemma}
\label{lem:comndProlimits}
The functor $\widehat{p\mathrm{Sym}}\colon \mathrm{coOp}(\mathcal{C}) \to \mathrm{coMnd}^{\mathrm{filt}}(\mathrm{Pro}(\mathcal{C}))$ preserves finite limits and eventually constant limits. (Recall that it is our convention that all cooperads $\mathcal{Q}$ are reduced and that $\mathcal{Q}(1) = \mathbf{1}$.)
\end{lemma}
\begin{proof}
We start with the case of finite limits. Every cooperad $\mathrm{coOp}(\mathcal{C})$ is naturally the totalization of a $\mathrm{forget}$-split cosimplicial object of cofree cooperads, namely the `cobar resolution' $\mathrm{cofree}^{\bullet + 1}(\mathcal{Q})$. (To be clear, here $\mathrm{cofree}(\mathcal{Q})$ denotes the cofree cooperad on the coaugmentation coideal of $\mathcal{Q}$, so that $\mathrm{cofree}(\mathcal{Q})$ still has the monoidal unit as its arity 1 term.) The functor $\widehat{p\mathrm{Sym}}$ preserves such totalizations by Lemma \ref{lem:pSymsplitTot}, so that we may reduce to checking that the composition $\widehat{p\mathrm{Sym}} \circ \mathrm{cofree}$ preserves finite limits. In the proof of Lemma \ref{lem:laxcoalgcofree} above we showed that the following diagram commutes:
\[
\begin{tikzcd}
\mathrm{coOp}(\mathcal{C}) \ar{r}{\widehat{p\mathrm{Sym}}} &  \mathrm{coMnd}^{\mathrm{filt}}(\mathrm{Pro}(\mathcal{C})) \\
\mathrm{SSeq}(\mathcal{C}) \ar{u}{\mathrm{cofree}}\ar{r}{\widehat{p\mathrm{Sym}}} & \mathrm{End}^{\mathrm{filt}}(\mathrm{Pro}(\mathcal{C})). \ar{u}{\mathrm{cofree}}
\end{tikzcd}
\]
Since the vertical functors are right adjoint, so in particular preserve finite limits, we can reduce to checking that the bottom horizontal arrow preserves finite limits. This is part of the conclusion of Lemma \ref{lem:pSymlimits}.

The reasoning for the case of eventually constant limits proceeds in entirely the same way; the only extra observation that is needed is that for an eventually constant diagram $\{\mathcal{Q}_i\}_{i \geq 1}$ of cooperads, the diagram $\{\mathrm{cofree}^{\circ k}\mathcal{Q}_i\}_{i \geq 1}$ is still eventually constant for every value of $k \geq 0$. This is clear from the explicit construction of cofree objects discussed above, specifically the observation that for a symmetric sequence $X$ the map $T^{(n+1)}(X) \to T^{(n)}(X)$ is an isomorphism in arities $\leq n+1$.
\end{proof}

We can now conclude this section with the promised proof.

\begin{proof}[Proof of Lemma \ref{lem:laxcoalgpb}]
We wish to show that the assignment $\mathcal{Q} \mapsto \mathrm{coAlg}^{\mathrm{dp}}_{\mathcal{Q}}(\mathcal{C})$ preserves the pullback square
\[
\begin{tikzcd}
\varphi^n\mathcal{Q} \ar{r}\ar{d} & \mathrm{cofree}(\mathcal{Q}(n)) \ar{d} \\
\varphi^{n-1}\mathcal{Q} \ar{r} & \mathrm{cofree}(\varphi^{n-1}\mathcal{Q}(n))
\end{tikzcd}
\]
and the inverse limit $\mathcal{Q} = \varprojlim_n \varphi^n\mathcal{Q}$. The functor
\[
\mathrm{coMnd}^{\mathrm{filt}}(\mathrm{Pro}(\mathcal{C})) \to (\mathrm{Cat}_{\infty})_{/\mathrm{Pro}(\mathcal{C})}\colon S \mapsto (\mathrm{coAlg}_S(\mathrm{Pro}(\mathcal{C})) \xrightarrow{\mathrm{forget}} \mathrm{Pro}(\mathcal{C}))
\]
preserves limits by a variation of the argument used in the proof of Theorem \ref{thm:algdecomposition}. Therefore it will suffice to see that the functor
\[
\widehat{p\mathrm{Sym}}\colon \mathrm{coOp}(\mathcal{C}) \to \mathrm{coMnd}^{\mathrm{filt}}(\mathrm{Pro}(\mathcal{C}))
\]
preserves the pullback square and the inverse limit described above. Note that the second is an eventually constant limit. Hence the conclusion now follows from Lemma \ref{lem:comndProlimits}.
\end{proof}

\section{Functoriality of the composition product}
\label{app:functoriality}



Let $\mathcal{C}$ and $\mathcal{D}$ be presentably symmetric monoidal $\infty$-categories. The aim of this appendix is to explain how, for an operad $\mathcal{O}$ in $\mathcal{C}$, a lax monoidal functor $f\colon \mathcal{C} \to \mathcal{D}$ induces a functor
\[
f_!\colon \mathrm{Alg}_{\mathcal{O}}(\mathcal{C}) \to \mathrm{Alg}_{f_!\mathcal{O}}(\mathcal{D}).
\]
Dually, for a colimit-preserving oplax monoidal functor $g\colon \mathcal{D} \to \mathcal{C}$ and a cooperad $\mathcal{Q}$ in $\mathcal{D}$, we will construct a corresponding functor
\[
g_!\colon \mathrm{coAlg}^{\mathrm{dp}}_{\mathcal{Q}}(\mathcal{D}) \to \mathrm{coAlg}^{\mathrm{dp}}_{g_!\mathcal{Q}}(\mathcal{C}),
\]
at least in the case where $\mathcal{C}$ and $\mathcal{D}$ are stable. (In the previous section we have only defined divided power coalgebras in a stable context.)

In Section \ref{subsec:laxmonoidal} we construct from a lax symmetric monoidal functor $f\colon \mathcal{C} \to \mathcal{D}$ a lax monoidal functor $f_!\colon \mathrm{SSeq}(\mathcal{C}) \to \mathrm{SSeq}(\mathcal{D})$ with respect to the composition product of symmetric sequences. It is not straightforward to do this from the description of the composition product we gave in Section \ref{sec:barcobar}. Instead, we will use an alternative construction of the composition product on $\mathrm{SSeq}(\mathcal{C})$ due to Haugseng \cite{haugsengsymseq}. The equivalence of the two approaches was shown by Arakawa \cite{arakawa}. In Section \ref{subsec:laxmonoidalalgcoalg} we will construct the functors $f_!$ and $g_!$ on algebra and divided power coalgebra categories alluded to above. Throughout this appendix, we allow our symmetric sequences $A$ to have a nontrivial term $A(0)$. (Recall that in the main text we worked with symmetric sequences indexed over \emph{non-empty} finite sets.) This extra flexibility is useful in order to discuss the action of $\mathrm{SSeq}(\mathcal{C})$ on $\mathcal{C}$ in Section \ref{subsec:laxmonoidalalgcoalg}.

\subsection{Lax monoidal functors and the composition product}
\label{subsec:laxmonoidal}

Let us fix some notation to describe Haugseng's construction (see also Arakawa's exposition in \cite[Section 5]{arakawa}). We recall the category $\mathbf{\Delta}_{\mathbb{F}}$ of \cite[Definition 4.1.1]{haugsengsymseq}. Its objects are finite strings of morphisms of finite sets
\[
S_0 \to S_1 \to \cdots \to S_n
\] 
with $n \geq 0$. A morphism from such a string to one of the form $(T_0 \to \cdots \to T_m)$ is a morphism $\alpha\colon {[n]} \to {[m]}$ in $\mathbf{\Delta}$ together with injections $S_i \to T_{\alpha(i)}$ for which each of the squares
\[
\begin{tikzcd}
S_i \ar{r}\ar{d} & S_{i+1} \ar{d} \\
T_{\alpha(i)} \ar{r} & T_{\alpha(i+1)}
\end{tikzcd}
\] 
commutes and is a pullback. We can think of an object of $\mathbf{\Delta}_{\mathbb{F}}$ as a layered forest with leaves $S_0$ and roots $S_n$. The edge set of the forest is $\coprod_{i=0}^n S_i$ and the maps in the string determine how these edges are attached to each other via vertices. To be precise, for a string $S_\bullet$ as above we define its pointed vertex set as $V(S_\bullet) := (\coprod_{i > 0} S_i)_+$. This determines a functor $V\colon \mathbf{\Delta}^{\mathrm{op}}_{\mathbb{F}} \to \mathrm{Fin}_*$. Also, observe that $\mathbf{\Delta}_{\mathbb{F}}$ admits an evident forgetful functor to $\mathbf{\Delta}$, for which we do not introduce separate notation.

Recall that a morphism in $\mathbf{\Delta}$ is called \emph{inert} if it is the inclusion of a subinterval, or in other words if its image is convex. Write $\Sigma$ for the full subcategory of the arrow category $\mathrm{Fun}(\Delta^1, \mathbf{\Delta}^{\mathrm{op}})$ on the inert morphisms. Define a category $\Sigma\mathbf{\Delta}_{\mathbb{F}}$ by the pullback square
\[
\begin{tikzcd}
\Sigma\mathbf{\Delta}^{\mathrm{op}}_{\mathbb{F}} \ar{r}\ar{d} & \Sigma \ar{d}{\mathrm{ev}_1} \\
\mathbf{\Delta}^{\mathrm{op}}_{\mathbb{F}} \ar{r} & \mathbf{\Delta}^{\mathrm{op}}.
\end{tikzcd}
\]
There is a projection $\pi_{\mathrm{Fin}_*}\colon \Sigma\mathbf{\Delta}^{\mathrm{op}}_{\mathbb{F}} \to \mathrm{Fin}_*$ defined as the composition of the functors 
\[
\Sigma\mathbf{\Delta}^{\mathrm{op}}_{\mathbb{F}} \to \mathbf{\Delta}^{\mathrm{op}}_{\mathbb{F}} \xrightarrow{V} \mathrm{Fin}_* 
\]
and a projection $\pi_{\mathbf{\Delta}^{\mathrm{op}}}$ to $\mathbf{\Delta}^{\mathrm{op}}$ via 
\[
\Sigma\mathbf{\Delta}^{\mathrm{op}}_{\mathbb{F}} \to \Sigma \xrightarrow{\mathrm{ev}_0} \mathbf{\Delta}^{\mathrm{op}}.
\]

\begin{definition}[cf. Definition 5.5 of \cite{arakawa}]
Let $\mathcal{C}^\otimes \to \mathrm{Fin}_*$ be an $\infty$-operad. Then we define a functor of $\infty$-categories $\widehat{\mathrm{SSeq}_H}(\mathcal{C})^{\circ} \to \mathbf{\Delta}^{\mathrm{op}}$ characterized by the universal property
\[
\mathrm{Fun}_{\mathbf{\Delta}^{\mathrm{op}}}(K,\widehat{\mathrm{SSeq}_H}(\mathcal{C})^{\circ}) \cong \mathrm{Fun}_{\mathrm{Fin}_*}(K \times_{\mathbf{\Delta}^{\mathrm{op}}} \Sigma\mathbf{\Delta}^{\mathrm{op}}_{\mathbb{F}}, \mathcal{C}^{\otimes}),
\]
naturally in maps of simplicial sets $K \to \mathbf{\Delta}^{\mathrm{op}}$. We define $\mathrm{SSeq}_H(\mathcal{C})^{\circ}$ to be the full subcategory of $\widehat{\mathrm{SSeq}_H}(\mathcal{C})^{\circ}$ spanned by functors 
\[
\{[n]\} \times_{\mathbf{\Delta}^{\mathrm{op}}} \Sigma\mathbf{\Delta}^{\mathrm{op}}_{\mathbb{F}} \to \mathcal{C}^{\otimes}
\]
sending every morphism to an inert morphism in $\mathcal{C}^\otimes$.
\end{definition}

The fiber of $\mathrm{SSeq}_H(\mathcal{C})^{\circ}$ at ${[1]} \in \mathbf{\Delta}^{\mathrm{op}}$ is easily identified with the $\infty$-category $\mathrm{SSeq}(\mathcal{C})$ of symmetric sequences in $\mathcal{C}$, see \cite[Remark 5.6]{arakawa}. In fact, if the $\infty$-operad $\mathcal{C}^{\otimes}$ is a symmetric monoidal $\infty$-category for which the tensor product preserves colimits over countable groupoids in each variable, then the functor $\mathrm{SSeq}_H(\mathcal{C})^{\circ} \to \mathbf{\Delta}^{\mathrm{op}}$ is a monoidal $\infty$-category by \cite[Corollary 5.7]{arakawa}. It does in fact encode the composition product, by virtue of the following:

\begin{proposition}[See Proposition 5.18 of \cite{arakawa}.]
\label{prop:SSeqcomparison}
Let $\mathcal{C}^\otimes$ be a symmetric monoidal $\infty$-category for which the tensor product preserves colimits over countable groupoids in each variable. Then there is an equivalence of monoidal $\infty$-categories
\[
(\mathrm{SSeq}(\mathcal{C}),\circ) \simeq \mathrm{SSeq}_H(\mathcal{C})^{\circ}.
\]
\end{proposition}

We will now investigate the functoriality of the construction $\mathrm{SSeq}_H(\mathcal{C})^{\circ}$ in $\mathcal{C}^{\otimes}$. Let $f\colon \mathcal{C}^{\otimes} \to \mathcal{D}^\otimes$ be a morphism of $\infty$-operads. Then by construction we obtain a functor
\[
\widehat{\mathrm{SSeq}_H}(\mathcal{C})^{\circ} \to \widehat{\mathrm{SSeq}_H}(\mathcal{D})^{\circ} 
\]
over $\mathbf{\Delta}^{\mathrm{op}}$. Since $f$ preserves inert morphisms by hypothesis, the functor above restricts to a functor
\[
f_!\colon \mathrm{SSeq}_H(\mathcal{C})^{\circ} \to \mathrm{SSeq}_H(\mathcal{D})^{\circ} 
\]
over $\mathbf{\Delta}^{\mathrm{op}}$.

\begin{lemma}
\label{lem:SSeqfunctoriality}
The functor $f_!$ constructed above preserves inert morphisms. In particular, if $\mathcal{C}^\otimes$ and $\mathcal{D}^\otimes$ are symmetric monoidal $\infty$-categories satisfying the hypothesis of Proposition \ref{prop:SSeqcomparison}, then $f_!$ is a lax monoidal functor between monoidal $\infty$-categories.
\end{lemma}
\begin{proof}
Let $\beta$ be a morphism of $\mathrm{SSeq}_H(\mathcal{C})^{\circ}$ lying over an inert map $\alpha$ in $\mathbf{\Delta}^{\mathrm{op}}$. Write $(\Sigma\mathbf{\Delta}^{\mathrm{op}}_{\mathbb{F}})_{\alpha}$ for the pullback of 
\[
\Delta^1 \xrightarrow{\alpha} \mathbf{\Delta}^{\mathrm{op}} \leftarrow \Sigma\mathbf{\Delta}^{\mathrm{op}}_{\mathbb{F}},
\]
so that $\beta$ is, by definition, determined by a map 
\[
\widehat{\beta}\colon (\Sigma\mathbf{\Delta}^{\mathrm{op}}_{\mathbb{F}})_{\alpha} \to \mathcal{C}^{\otimes}
\]
lying over $\mathrm{Fin}_*$. By \cite[Proposition 5.7]{arakawa}, the morphism $\beta$ is coCartesian over $\mathbf{\Delta}^{\mathrm{op}}$ if and only if $\widehat{\beta}$ carries every morphism lying over an inert morphism of $\mathrm{Fin}_*$ to an inert morphism of $\mathcal{C}^{\otimes}$. Clearly this property is preserved under composition with the morphism of $\infty$-operads $f\colon \mathcal{C}^{\otimes} \to \mathcal{D}^{\otimes}$, since $f$ preserves inerts by definition.
\end{proof}

\begin{remark}
    \label{rmk:laxmonoidalfunctors}
If $\mathcal{C}^\otimes$ and $\mathcal{D}^\otimes$ are symmetric monoidal $\infty$-categories, then a morphism of $\infty$-operads $f\colon \mathcal{C}^{\otimes} \to \mathcal{D}^{\otimes}$ is the same thing as a lax symmetric monoidal functor. In particular, we conclude that if $\mathcal{C}^{\otimes}$ and $\mathcal{D}^{\otimes}$ satisfy the hypothesis of Proposition \ref{prop:SSeqcomparison}, then $f$ induces a lax monoidal functor 
\[
f_!\colon (\mathrm{SSeq}(\mathcal{C}),\circ) \to (\mathrm{SSeq}(\mathcal{D}),\circ)
\]
by virtue of Proposition \ref{prop:SSeqcomparison} and Lemma \ref{lem:SSeqfunctoriality}. A lax monoidal functor induces a functor on algebra objects; in particular, we also obtain a functor $f_!\colon \mathrm{Op}(\mathcal{C}) \to \mathrm{Op}(\mathcal{D})$ sending operads in $\mathcal{C}$ to operads in $\mathcal{D}$.
\end{remark}

\subsection{Lax symmetric monoidal functors and (co)algebras}
\label{subsec:laxmonoidalalgcoalg}

For simplicity, we will assume that $\mathcal{C}^{\otimes}$ and $\mathcal{D}^{\otimes}$ are presentably symmetric monoidal $\infty$-categories throughout this section. In particular, they satisfy the hypothesis of Proposition \ref{prop:SSeqcomparison}. 

Note that we may identify $\mathcal{C}$ with the full subcategory of $\mathrm{SSeq}(\mathcal{C})$ on symmetric sequences concentrated in arity 0. Moreover, if $X$ is such a symmetric sequence, then for any other symmetric sequence $A$ the composition product $A \circ X$ is still concentrated in arity 0; in fact, it corresponds to the object of $\mathcal{C}$ that we have denoted $\mathrm{Sym}_A(X)$ in the main text. In this way the composition product gives $\mathcal{C}$ the structure of an $\infty$-category left tensored over the monoidal $\infty$-category $\mathrm{SSeq}(\mathcal{C})$ (see \cite[Section 4.2.2]{higheralgebra}). 

To be more precise, let $\mathcal{LM}^{\otimes}$ denote the $\infty$-operad on two colors $a$ and $m$ that parametrizes an associative algebra and a left module over it. Then there is an evident morphism of $\infty$-operads to the associative operad $p\colon \mathcal{LM}^{\otimes} \to \mathrm{Ass}^{\otimes}$ encoding the fact that any associative algebra may be viewed as a left module over itself (see \cite[Remark 4.2.1.5]{higheralgebra}). Encoding the monoidal $\infty$-category of symmetric sequences with composition product as a coCartesian fibration $\mathrm{SSeq}(\mathcal{C})^{\circ} \to \mathrm{Ass}^{\otimes}$, we form the pullback
\[
(\mathrm{SSeq}(\mathcal{C})^{\circ} \curvearrowright \mathrm{SSeq}(\mathcal{C})) := \mathcal{LM}^{\otimes} \times_{\mathrm{Ass}^{\otimes}} \mathrm{SSeq}(\mathcal{C})^{\circ}.
\]
This is a coCartesian fibration over $\mathcal{LM}^{\otimes}$ encoding the fact that $\mathrm{SSeq}(\mathcal{C})$ is left-tensored over $\mathrm{SSeq}(\mathcal{C})^{\circ}$. Now the evident full subcategory $(\mathrm{SSeq}(\mathcal{C})^{\circ} \curvearrowright \mathcal{C})$ whose fiber over $m$ is $\mathcal{C}$, viewed as symmetric sequences concentrated in degree zero, encodes the left tensoring of $\mathcal{C}$ over $\mathrm{SSeq}(\mathcal{C})^{\circ}$. Here we use the fact that $\mathcal{C}$ is a left ideal in $\mathrm{SSeq}(\mathcal{C})^{\circ}$ to see that this is still a coCartesian fibration over $\mathcal{LM}^{\otimes}$.

\begin{remark}
\label{rmk:laxmonoidalfunctors2}
Let us explicitly record the functoriality of the above construction: if $f\colon \mathcal{C}^{\otimes} \to \mathcal{D}^{\otimes}$ is a lax symmetric monoidal functor, then by Remark \ref{rmk:laxmonoidalfunctors} we obtain a lax monoidal functor $f_!\colon \mathrm{SSeq}(\mathcal{C})^{\circ} \to \mathrm{SSeq}(\mathcal{D})^{\circ}$, and the construction above now yields a morphism of $\infty$-operads
\[
f_!\colon (\mathrm{SSeq}(\mathcal{C})^{\circ} \curvearrowright \mathcal{C}) \to (\mathrm{SSeq}(\mathcal{D})^{\circ} \curvearrowright \mathcal{D}). 
\]
\end{remark}

For the following definition, we observe that an $\mathcal{LM}^{\otimes}$-algebra in $(\mathrm{SSeq}(\mathcal{C})^{\circ} \curvearrowright \mathcal{C})$ is precisely the data of an operad $\mathcal{O}$ and an $\mathcal{O}$-algebra $X$ in $\mathcal{C}$. The inclusion $\mathrm{Ass}^{\otimes} \to \mathcal{LM}^{\otimes}$ determines a forgetful functor from such $\mathcal{LM}^{\otimes}$-algebras to associative algebras in $\mathrm{SSeq}(\mathcal{C})$, i.e., operads in $\mathcal{C}$. For an operad $\mathcal{O}$ in $\mathcal{C}$, we have an equivalence
\[
\mathrm{Alg}_{\mathcal{O}}(\mathcal{C}) \simeq \{\mathcal{O}\} \times_{\mathrm{Op}(\mathcal{C})} \mathrm{Alg}_{\mathcal{LM}^{\otimes}}(\mathrm{SSeq}(\mathcal{C})^{\circ} \curvearrowright \mathcal{C}).
\]

\begin{definition}
    \label{def:alglaxfunctoriality}
Let $f\colon \mathcal{C}^{\otimes} \to \mathcal{D}^{\otimes}$ be a lax symmetric monoidal functor. Combining the functor $f_!$ of Remark \ref{rmk:laxmonoidalfunctors2} with the identification of $\mathrm{Alg}_{\mathcal{O}}(\mathcal{C})$ directly above, we find a functor
\[
\mathrm{Alg}_{\mathcal{O}}(\mathcal{C}) \to \mathrm{Alg}_{f_!\mathcal{O}}(\mathcal{D})
\]
for which we still write $f_!$.
\end{definition}

This establishes the functoriality of algebra $\infty$-categories $\mathrm{Alg}_{\mathcal{O}}(\mathcal{C})$ with respect to lax symmetric monoidal functors $\mathcal{C}^{\otimes} \to \mathcal{D}^{\otimes}$. The real purpose of this section is to explain the dual construction, namely how an oplax symmetric monoidal functor $g\colon \mathcal{D}^{\otimes} \to \mathcal{C}^\otimes$ determines a functor on divided power coalgebra $\infty$-categories $g_!\colon \mathrm{coAlg}^{\mathrm{dp}}_{\mathcal{Q}}(\mathcal{D}) \to \mathrm{coAlg}^{\mathrm{dp}}_{g_!\mathcal{Q}}(\mathcal{C})$.

To be more precise, for the remainder of this section we will assume that $g$ preserves colimits and that $\mathcal{C}$ and $\mathcal{D}$ are stable. By the adjoint functor theorem, $g$ then admits a right adjoint $f$, which acquires a lax symmetric monoidal structure by \cite[Corollary 3.4.8]{haugsenghebestreitlinskensnuiten}. By Remark \ref{rmk:laxmonoidalfunctors}, $f$ induces a lax monoidal functor $f_!\colon \mathrm{SSeq}(\mathcal{C}) \to \mathrm{SSeq}(\mathcal{D})$. This functor $f_!$ admits a left adjoint $g_!$ (simply applying $g$ levelwise to a symmetric sequence), which by another application of op. cit. now acquires an oplax monoidal structure.

\begin{remark}
Oplax monoidal functors preserve coalgebras, so that in the setting above we obtain a functor
\[
g_!\colon \mathrm{coOp}(\mathcal{D}) \to \mathrm{coOp}(\mathcal{C}). 
\]
In particular, this justifies our use of the notation $g_!\mathcal{Q}$ for a cooperad $\mathcal{Q}$ in $\mathcal{D}$.
\end{remark}

To get the desired functoriality of divided power coalgebras, we will need to upgrade the oplax monoidal functor $g_!\colon \mathrm{SSeq}(\mathcal{D}) \to \mathrm{SSeq}(\mathcal{C})$ to one that is compatible with the left tensorings of $\mathrm{Pro}(\mathcal{D})$ and $\mathrm{Pro}(\mathcal{C})$ over these monoidal $\infty$-categories via the `completed' symmetric algebra functors we introduced in Appendix \ref{app:coalgebras}. Again, to achieve this we first work at the level of the right adjoint $f$.

Recall from Section \ref{sec:pro} the action of $\mathrm{SSeq}(\mathcal{C})^{\circ}$ on $\mathrm{Pro}(\mathcal{C})$ via the functor $\widehat{p\mathrm{Sym}}$. As before, this can be packaged as an $\mathcal{LM}^{\otimes}$-monoidal $\infty$-category $(\mathrm{SSeq}(\mathcal{C})^{\circ} \curvearrowright \mathrm{Pro}(\mathcal{C}))$. Writing $pf\colon \mathrm{Pro}(\mathcal{C}) \to \mathrm{Pro}(\mathcal{D})$ for the prolongation of $f$ (which is right adjoint to the prolongation $pg$), we have for every symmetric sequence $A \in \mathrm{SSeq}(\mathcal{C})$ a natural transformation
\[
\widehat{p\mathrm{Sym}}_{f_!A} \circ pf \xrightarrow{\widehat{\gamma}_A} pf \circ \widehat{p\mathrm{Sym}}_{A},
\] 
described as follows. Since both sides preserve filtered limits, the natural transformation is uniquely determined by its restriction to constant pro-objects. For each $n \geq 0$ we have natural transformations
\[
\mathrm{Sym}^{\leq n}_{f_!A} \circ f \xrightarrow{\gamma_A} f \circ \mathrm{Sym}^{\leq n}_A
\]
from our construction of the morphism of left-tensored $\infty$-categories 
\[
f_!\colon (\mathrm{SSeq}(\mathcal{C})^{\circ} \curvearrowright \mathcal{C}) \to (\mathrm{SSeq}(\mathcal{D})^{\circ} \curvearrowright \mathcal{D})
\]
in Remark \ref{rmk:laxmonoidalfunctors2}, applied to the truncation $\tau_n A$. Taking prolongations and passing to the limit over $n$ gives the natural transformation $\widehat{\gamma}_A$ displayed before. Combining these with the lax monoidal structure of $f_!\colon \mathrm{SSeq}(\mathcal{C}) \to \mathrm{SSeq}(\mathcal{D})$ already constructed, this produces a morphism of $\infty$-operads
\[
f_!\colon (\mathrm{SSeq}(\mathcal{C})^{\circ} \curvearrowright \mathrm{Pro}(\mathcal{C})) \to (\mathrm{SSeq}(\mathcal{D})^{\circ} \curvearrowright \mathrm{Pro}(\mathcal{D}))
\]
over $\mathcal{LM}^{\otimes}$ or, in other words, a lax $\mathcal{LM}^{\otimes}$-monoidal functor. Applying \cite[Corollary 4.8]{haugsenghebestreitlinskensnuiten} yet again, but now for $\mathcal{LM}^\otimes$-monoidal structures, we find an oplax $\mathcal{LM}^\otimes$-monoidal left adjoint $g_!$ from $(\mathrm{SSeq}(\mathcal{D})^{\circ} \curvearrowright \mathrm{Pro}(\mathcal{D}))$ to $(\mathrm{SSeq}(\mathcal{C})^{\circ} \curvearrowright \mathrm{Pro}(\mathcal{C}))$.

In parallel with Definition \ref{def:alglaxfunctoriality}, we now find that for a cooperad $\mathcal{Q}$ in $\mathcal{D}$, the functor $g$ induces a functor
\[
g_!\colon \mathrm{coAlg}_{\widehat{p\mathrm{Sym}}_{\mathcal{Q}}}(\mathrm{Pro}(\mathcal{D})) \to \mathrm{coAlg}_{\widehat{p\mathrm{Sym}}_{g_!\mathcal{Q}}}(\mathrm{Pro}(\mathcal{C})).
\]
By construction, this functor is compatible with the forgetful functors to $\mathrm{Pro}(\mathcal{D})$ and $\mathrm{Pro}(\mathcal{C})$ and the prolongation $pg\colon \mathrm{Pro}(\mathcal{D}) \to \mathrm{Pro}(\mathcal{C})$ of $g$ between those. In particular, the functor $g_!$ above preserves coalgebras whose underlying object is pro-constant. In conclusion, this $g_!$ also restricts to a functor
\[
\mathrm{coAlg}_{\mathcal{Q}}^{\mathrm{dp}}(\mathcal{D}) \to \mathrm{coAlg}_{g_!\mathcal{Q}}^{\mathrm{dp}}(\mathcal{C}), 
\]
as desired.





\bibliographystyle{plain}
\bibliography{biblio}

\begin{thebibliography}{10}

\bibitem{stacksproject}
Stacks project.
\newblock https://stacks.math.columbia.edu.

\bibitem{amabel}
A.~Amabel.
\newblock {Poincar\'e/Koszul Duality for General Operads}.
\newblock arXiv preprint arXiv:1910.09076.

\bibitem{anel}
M.~Anel.
\newblock Cofree coalgebras over operads and representative functions.
\newblock {\em arXiv preprint arXiv:1409.4688}, 2014.

\bibitem{arakawa}
Kensuke Arakawa.
\newblock On the equivalence of {Brantner's and Chu--Haugseng's} approaches to
  enriched $\infty$-operads.
\newblock {\em arXiv preprint arXiv:2603.23019}, 2026.

\bibitem{ayalafrancis}
D.~Ayala and J.~Francis.
\newblock Zero-pointed manifolds.
\newblock {\em Journal of the Institute of Mathematics of Jussieu},
  20(3):785--858, 2021.

\bibitem{basterramandell}
M.~Basterra and M.A. Mandell.
\newblock Homology and cohomology of ${E}_{\infty}$ ring spectra.
\newblock {\em Mathematische Zeitschrift}, 249(4):903--944, 2005.

\bibitem{behrensrezk}
M.~Behrens and C.~Rezk.
\newblock The {B}ousfield-{K}uhn functor and topological {A}ndr\'{e}-{Q}uillen
  cohomology.
\newblock preprint, 2012.

\bibitem{blansblom}
Max Blans and Thomas Blom.
\newblock On the chain rule in goodwillie calculus.
\newblock {\em arXiv preprint arXiv:2410.20504}, 2024.

\bibitem{bousfieldkan}
Aldridge~Knight Bousfield and Daniel~Marinus Kan.
\newblock {\em Homotopy limits, completions and localizations}, volume 304.
\newblock Springer Science \& Business Media, 1972.

\bibitem{brantnerthesis}
D.L.B. Brantner.
\newblock {\em {The Lubin-Tate theory of spectral Lie algebras}}.
\newblock PhD thesis, Harvard University, 2017.

\bibitem{BCN}
D.L.B. Brantner, R.~Campos, and J.~Nuiten.
\newblock Pd operads and explicit partition lie algebras.
\newblock {\em arXiv preprint arXiv:2104.03870}, 2021.

\bibitem{brantnerheutsKD}
D.L.B. Brantner and G.S.K.S. Heuts.
\newblock Symmetric sequences, $\infty$-operads, and the universal {L}ie
  algebra.
\newblock {\em to appear}.

\bibitem{brantnermathew}
D.L.B. Brantner and A.~Mathew.
\newblock Deformation theory and partition lie algebras.
\newblock {\em arXiv preprint arXiv:1904.07352}, 2019.

\bibitem{calaquecamposnuiten}
D.~Calaque, R.~Campos, and J.~Nuiten.
\newblock Moduli problems for operadic algebras.
\newblock {\em Journal of the London Mathematical Society}, 106(4):3450--3544,
  2022.

\bibitem{chingbar}
M.~Ching.
\newblock Bar-cobar duality for operads in stable homotopy theory.
\newblock {\em Journal of {T}opology}, 2012.

\bibitem{chingharper}
M.~Ching and J.E. Harper.
\newblock Derived koszul duality and tq-homology completion of structured ring
  spectra.
\newblock {\em Advances in Mathematics}, 341:118--187, 2019.

\bibitem{francisgaitsgory}
J.~Francis and D.~Gaitsgory.
\newblock Chiral {K}oszul duality.
\newblock {\em Selecta Mathematica}, 18(1):27--87, 2012.

\bibitem{gaitsgoryrozenblyum}
D.~Gaitsgory and N.~Rozenblyum.
\newblock A study in derived algebraic geometry: Volume ii: Deformations, lie
  theory and formal geometry.
\newblock {\em American Mathematical Society}, 2017.

\bibitem{getzlerjones}
E.~Getzler and J.D.S. Jones.
\newblock Operads, homotopy algebra and iterated integrals for double loop
  spaces.
\newblock Available online at arXiv hep-th/9305013.

\bibitem{ginzburgkapranov}
V.~Ginzburg and M.~Kapranov.
\newblock Koszul duality for operads.
\newblock {\em Duke {M}athematical {J}ournal}, 76(1):203--272, 1994.

\bibitem{haugsengsymseq}
R.~Haugseng.
\newblock $\infty$-operads via symmetric sequences.
\newblock arXiv preprint arXiv:1708.09632.

\bibitem{haugsenghebestreitlinskensnuiten}
Rune Haugseng, Fabian Hebestreit, Sil Linskens, and Joost Nuiten.
\newblock Lax monoidal adjunctions, two-variable fibrations and the calculus of
  mates.
\newblock {\em Proceedings of the London Mathematical Society},
  127(4):889--957, 2023.

\bibitem{heine}
H.~Heine.
\newblock A monadicity theorem for higher algebraic structures.
\newblock arXiv preprint arXiv:1712.00555.

\bibitem{heutsgoodwillie}
G.S.K.S. Heuts.
\newblock Goodwillie approximations to higher categories.
\newblock Available online at arXiv:1510.03304, 2015.

\bibitem{kuhnmccord}
N.~J. Kuhn.
\newblock {The McCord model for the tensor product of a space and a commutative
  ring spectrum}.
\newblock In {\em {Categorical Decomposition Techniques in Algebraic
  Topology}}, pages 213--235. Springer, 2003.

\bibitem{lodayvallette}
J.-L. Loday and B.~Vallette.
\newblock {\em Algebraic operads}, volume 346.
\newblock Springer Science \& Business Media, 2012.

\bibitem{HTT}
J.~Lurie.
\newblock {\em Higher topos theory}, volume 170.
\newblock Princeton University Press, 2009.

\bibitem{dagx}
J.~Lurie.
\newblock Derived {A}lgebraic {G}eometry {X}: {F}ormal {M}oduli {P}roblems.
\newblock Available online at math.harvard.edu/~lurie/, 2011.

\bibitem{higheralgebra}
J.~Lurie.
\newblock Higher algebra.
\newblock Available online at math.harvard.edu/~lurie/, 2014.

\bibitem{lurieRot}
J.~Lurie.
\newblock Rotation invariance in algebraic {K}-theory.
\newblock Available online at math.harvard.edu/~lurie/, 2015.

\bibitem{SAG}
J.~Lurie.
\newblock Spectral algebraic geometry.
\newblock {\em preprint}, 2018.

\bibitem{moore}
J.C. Moore.
\newblock Differential homological algebra.
\newblock In {\em Actes du Congres International des Math{\'e}maticiens (Nice,
  1970)}, volume~1, pages 335--339, 1970.

\bibitem{pridham}
J.P. Pridham.
\newblock Unifying derived deformation theories.
\newblock {\em Advances in Mathematics}, 224(3):772--826, 2010.

\bibitem{rationalhomotopy}
D.~G. Quillen.
\newblock Rational homotopy theory.
\newblock {\em Annals of Mathematics}, pages 205--295, 1969.

\end{thebibliography}

\end{document}